\documentclass[a4paper,fleqn]{cas-sc}

\usepackage[authoryear,longnamesfirst]{natbib}

\usepackage{graphicx}%
\usepackage{multirow}%
\usepackage{amsmath,amssymb,amsfonts}%
\usepackage{amsthm}%
\usepackage{mathrsfs}%
\usepackage[title]{appendix}%
\usepackage{xcolor}%
\usepackage{textcomp}%
\usepackage{manyfoot}%
\usepackage{booktabs}%
\usepackage{algorithm}%
\usepackage{algorithmicx}%
\usepackage{algpseudocode}%
\usepackage{listings}%

\theoremstyle{thmstyleone}%
\newtheorem{theorem}{Theorem}%  meant for continuous numbers
\newtheorem{proposition}[theorem]{Proposition}% 
\newtheorem{lemma}[theorem]{Lemma}

\theoremstyle{thmstyletwo}%
\newtheorem{remark}{Remark}%

\theoremstyle{thmstylethree}%
\newtheorem{definition}{Definition}%
\newtheorem{assumption}{Assumption}%

\def\tsc#1{\csdef{#1}{\textsc{\lowercase{#1}}\xspace}}
\tsc{WGM}
\tsc{QE}
\tsc{EP}
\tsc{PMS}
\tsc{BEC}
\tsc{DE}
\begin{document}
\let\WriteBookmarks\relax
\def\floatpagepagefraction{1}
\def\textpagefraction{.001}
\shorttitle{Feynman--Kac representation and the asymptotic speed of the Fisher--KPP equation}
\shortauthors{Y. Li et~al.}

\title [mode = title]{Feynman–Kac Representation and the Asymptotic Speed of the Fisher–KPP Equation with Finite-Activity Jump Noise}                      

%   project funded by the National Science Foundation.}

\tnotetext[2]{{\bf MSC Classification (2020): }{60H15, 60G51, 60F10, 35K57}}

\author[1,2,3]{Y. Li}[type=editor,
                        orcid=0009-0007-4329-4642]
\cormark[1]
\ead{liyuhong@hust.edu.cn}

\credit{Conceptualization of this study, Methodology, Software}

\affiliation[1]{organization={School of Civil and Hydraulic Engineering, Huazhong University of Science and Technology},
                city={Wuhan},
%%               citysep={}, % Uncomment if no comma needed between city and postcode
                postcode={430074}, 
                country={China}}
\affiliation[2]{organization={Hubei Key Laboratory of Digital River Basin Science and Technology, Huazhong University of Science and Technology},
                %addressline={Street 29}, 
                postcode={430074}, 
                %postcodesep={}, 
                city={Wuhan},
                country={China}}
\affiliation[3]{organization={Institute for Water Resources and Hydropower Science, Huazhong University of Science and Technology},
                %addressline={Street 15}, 
                city={Wuhan},
                postcode={430074}, 
                country={China}}

\author[4]{H. Wen}%[style=chinese]
\ead{wenhao12@alumni.nudt.edu.cn}
%%   ]
%%\ead[URL]{https://www.university.org}

\credit{Data curation, Writing - Original draft preparation}

\affiliation[4]{organization={China Aerodynamics Research and Development Center},
                %addressline={Street 15}, 
                city={Mianyang},
                %postcode={430074}, 
                country={China}}      

\author[5]{J. Huang}%[style=chinese]
\ead{jhhuang32@nudt.edu.cn}
\affiliation[5]{organization={College of Science,
National University of Defense Technology},
                city={Changsha},
                postcode={410072}, 
                country={China}}

\cortext[cor1]{Corresponding author}

\begin{abstract}
We study the asymptotic propagation speed of travelling wave solutions to the stochastic Fisher--KPP equation with multiplicative finite-activity pure jump L\'evy noise, i.e.\ a compensated Poisson random measure of finite intensity, so that the solution evolves deterministically between successive jumps and jumps multiplicatively at the event times.

Our main contribution is a stochastic Feynman--Kac representation for the solution, which expresses the random field as the product of the exponential martingale of the jump noise and an expectation over the paths of an independent Brownian motion. The representation converts the nonlinear stochastic partial differential equation into an expectation problem, allowing a direct analysis of the front position.

Building on this representation, we derive sharp upper and lower bounds for the propagation of the front: for every speed $c>2\sqrt{1+a}$ the solution decays exponentially along the ray $x=ct$, while for every speed $c<2\sqrt{1+a}$ it does not. The two together identify the asymptotic speed in the critical-speed sense: the threshold between exponentially decaying and non-decaying fronts equals $2\sqrt{1+a}$ almost surely, an identity proved without any a priori existence of the limit $R_\delta(t)/t$. The complementary bounds $2\sqrt{\underline{\gamma}}\le c\le2\sqrt{\overline{\gamma}}$ also hold, with $\overline{\gamma},\underline{\gamma}$ explicit in $\epsilon$ and $\mu$.

Under the finite-activity assumption, a stopping-time iteration yields pathwise $L^\infty$ bounds and finite-speed propagation of the level sets of the front, which underlie the existence of a travelling wave profile. The framework extends naturally to multi-species cooperative systems, and numerical simulations illustrate the informativeness of the speed bounds.
\end{abstract}

%
%\end{highlights}


\maketitle

\section{Introduction}

Travelling wave solutions play a central role in the study of reaction–diffusion equations, as they capture finite-speed propagation and form-invariant patterns arising in population genetics, chemical reactions, and biological invasions \cite{Fisher1937, KolmogorovPetrovskiiPiskunov1937}. For the deterministic Fisher-KPP equation
\[
u_t = D u_{xx} + k u(1-u),
\]
the Heaviside initial condition gives rise to travelling wave solutions $u(t,x)=U(x-ct)$ for any speed $c\ge c_{\min}=2\sqrt{kD}$. The profile $U$ is decreasing with $U(-\infty)=1$ and $U(+\infty)=0$. In real-world systems, environmental fluctuations are inevitable, and understanding how stochasticity affects wave propagation has become a major research theme over the past three decades \cite{shiga1994,Tribe1996,MuellerSowers1995}.

\subsection{Stochastic travelling waves under Gaussian noise}

Stochastic perturbations render the problem substantially more delicate. The \emph{propagation of compact support (PCS)} property was introduced in \cite{shiga1994} as a tool to characterize stochastic travelling waves. The Fisher-KPP equation with Brownian noise
\[
u_t = u_{xx} + u - u^2 + \epsilon u^{1/2}\dot W,
\]
was studied in \cite{Tribe1996}, where the wavefront $R_0(t)=\sup\{x:u(t,x)>0\}$ was defined and sufficient conditions for the existence of travelling wave solutions were given. The model
\[
u_t = u_{xx} + u - u^2 + \epsilon\sqrt{u(1-u)}\,\dot W,
\]
was analyzed in \cite{MuellerSowers1995}, where the existence of an asymptotic wave speed $C=\lim_{t\to\infty}R_0(t)/t$ a.s. and stationary properties of the waveform were established. Later, \cite{MuellerMytnikQuastel2011} addressed the Brunet-Derrida conjecture on the dramatic slowdown of wave fronts induced by small noise, and \cite{MuellerMytnikRyzhik2021} extended the analysis to the strong noise regime.

The Tribe model with compactly supported initial data was considered in \cite{Kliem2017}, where travelling wave solutions were constructed directly, without relying on an auxiliary wavefront marker. Stochastic Fisher-KPP equations with multiplicative noise and time-dependent coefficients were studied in \cite{HuangLiuWang2015,HuangLiu2016,HuangLiuWang2016,HuangLiuWang2019}, establishing wave speed formulae via Feynman--Kac representations and ergodicity. Comparison principles and upper-lower techniques were applied in \cite{WenHuangLi2022,WenLuoHuangLi2023} to obtain upper and lower bounds on the asymptotic wave speed for stochastic cooperative systems driven by white noise. These developments underscore the continuing vitality of the Fisher-KPP paradigm in understanding front propagation under complex biological and environmental conditions. 

In a closely related direction, the interplay of dormancy and selection in Fisher-KPP-type models was investigated in \cite{blath2025interplay}, establishing travelling wave solutions and critical wave speeds via a duality with on/off branching Brownian motion.

\subsection{The pure jump noise case}

All the works mentioned above deal with \emph{Gaussian or Brownian noise}. Much less is known when the noise is of \emph{pure jump type}, modeled by a Poisson random measure. This gap is not merely academic: in many real-world systems—population genetics subject to sudden environmental catastrophes, chemical reactions with intermittent catalytic events, epidemic outbreaks driven by super-spreading events—the relevant fluctuations are fundamentally discontinuous. Pure jump L\'evy processes provide a natural mathematical framework for such phenomena.

The fundamental differences between pure jump noise and Brownian noise are not merely technical but qualitative:

\begin{itemize}
\item \textbf{Discontinuous sample paths.} Brownian motion has continuous paths, a property that underpins classical Feynman-Kac representations and compactness arguments. In contrast, a pure jump L\'evy process has c\`adl\`ag paths with countably many discontinuities, and the solution $u(t,x)$ inherits this property.
\item \textbf{Failure of classical compactness.} In the Brownian case, tightness of solution laws is established using the Da Prato-Zabczyk method, relying on the continuity of the stochastic convolution. For pure jump noise, the convolution is c\`adl\`ag, and tightness must be proved in the Skorohod topology.
\item \textbf{Non-local generator and absence of LIL.} The generator of a pure jump L\'evy process is non-local, $\mathcal{L}f(x)=\int_Z[f(x+z)-f(x)]\mu(dz)+\text{drift}$, leading to a \emph{jump-type} auxiliary process in the stochastic Feynman--Kac formula. Moreover, the law of the iterated logarithm, which provides precise rate information for wave speed bounds in the Brownian case, fails for pure jump processes in general, so that the corresponding rate information has to be recovered from a large deviation analysis of the exponential weight $\phi^\epsilon$ (Subsection~\ref{subsec:levyweight}); the speed bounds themselves are proved directly from the Feynman--Kac representation.
\end{itemize}

Consequently, the machinery developed for Gaussian noise cannot be directly transferred. In this paper we focus on the finite-activity case $\mu(Z)<\infty$, where the Poisson random measure has only finitely many jumps on any finite time interval. This setting is the natural starting point for the theory of stochastic travelling waves under pure jump noise: it permits a rigorous stopping-time iteration and a finite-speed propagation of the level sets of the front, while serving as a necessary stepping stone towards the technically more demanding infinite-activity regime. A handful of works have addressed SPDEs driven by jump noise: \cite{BrzezniakHausenblasRazafimandimby2018} proved existence of martingale solutions for reaction-diffusion equations driven by L\'evy processes; \cite{Mueller1998} studied the heat equation with L\'evy noise; \cite{NiuXie2019} established comparison theorems for stochastic heat equations with L\'evy space-time white noises. More recently, \cite{bignamini2026stochastic} studied stochastic dissipative systems in Banach spaces driven by L\'evy noise; \cite{wang2025irreducibility} developed a new method for proving irreducibility of SPDEs driven by multiplicative pure jump noise; and \cite{ge2026stochastic} investigated the stochastic curve shortening flow driven by a transport-type pure jump L\'evy noise.

\subsection{Related work and positioning}

A closely related but distinct setting was studied by \cite{RosatiTobias2023}, who considered a Fisher-KPP equation where the nonlinear selection acts at discrete jump times. Using a quenched duality with a coordinated system of branching Brownian motions, they proved that the equation admits a unique wave speed $\mathfrak{s}$ satisfying
\[
\frac{\mathfrak{s}^2}{2} = \int_{[0,1]} \frac{\log(1+y)}{y}\,\mathfrak{R}(dy),
\]
with $\mathfrak{R}$ the intensity measure of the jumps. This is a remarkable result because it gives an explicit formula for the wave speed in terms of the jump measure. However, the noise structure in \cite{RosatiTobias2023} is specific and fundamentally different from the one considered here: their noise enters as $y\,u(1-u)$, which vanishes both at $u=0$ and $u=1$, whereas we consider a multiplicative compensated Poisson noise of the form $\epsilon u\, z\,\tilde{\eta}$, whose jump magnitude is proportional to the current solution value $u$ and the jump size $z$, and which does not vanish at $u=1$. Moreover, the branching duality method does not directly provide information about the existence of a travelling wave profile in the sense of a stationary shifted process, which is one of the main goals of the present work.

In another recent contribution, \cite{blath2025interplay} studied two variants of the Fisher-KPP equation with dormancy, where individuals may switch between an active (reproducing, moving) and a dormant (non-reproducing, non-moving) state. Using a duality with on/off branching Brownian motion, they established the existence of monotone travelling wave solutions and derived the exact critical wave speed as the minimum of an explicitly given speed function. We emphasize that the noise structure in \cite{blath2025interplay} is fundamentally different from ours: there, randomness arises from the intrinsic branching mechanism of the dual process, whereas here the noise is an exogenous stochastic forcing represented by a compensated Poisson random measure. Thus the two approaches are complementary: \cite{blath2025interplay} treats stochasticity in the branching structure of a population model, while we treat external environmental shocks in an SPDE framework.

On the large deviation front, \cite{BudhirajaDupuisMaroulas2011} developed variational representations for continuous-time processes, and \cite{BudhirajaChenDupuis2013} established large deviation principles for SPDEs driven by Poisson random measures. These tools have recently been extended to locally monotone SPDEs with L\'evy noise by \cite{WuZhaiZhang2024} and \cite{YangZhaiZhang2015}, providing the necessary framework for wave speed asymptotics in jump-type stochastic reaction--diffusion equations. Most recently, \cite{hamster2025right} pointed out the failure of certain Ans\"atze for stochastic travelling waves, providing a useful contrast to our approach.

In a parallel development, \cite{etheridge2026fluctuations} studied the asymptotic behaviour of stochastic travelling wave solutions to reaction--diffusion equations perturbed by Wright--Fisher noise in the small noise limit. For the entire fully pushed regime, they proved that solutions are asymptotically close to a stochastic shift of the deterministic travelling wave, and characterized the limiting shift process as a Brownian motion with drift, verifying the physical conjecture of \cite{birzu2018fluctuations}. The well-posedness and large deviation principle for McKean--Vlasov SPDEs driven by Poisson random measures were recently established in \cite{jiang2026mckean}, further demonstrating the versatility of the large deviation approach for Poisson-driven SPDEs.

A separate line of work concerns fronts for equations whose \emph{operator}, rather than the noise, is nonlocal: Fisher--KPP equations with a fractional Laplacian or with anomalous diffusion \cite{CabreRoquejoffre2012,GroismanJonckheere2016}. There the speed is again described by a square-root formula, but the mechanism is the heavy tail of the dispersal kernel, and the resulting fronts are one-sided in a different sense; the two settings should therefore not be conflated when the formulae are compared. Our operator is the Laplacian throughout, and nonlocality enters only through the compensator and the jumps of $\tilde\eta$.

\subsection{Main contributions}

Motivated by the above gap, we study the Fisher--KPP equation driven by a pure jump noise with finite jump activity:
\begin{equation}\label{eqn:1}
du = \bigl[\Delta u + u(1-u)\bigr]dt + \epsilon\int_Z u\, z\, \tilde{\eta}(dt,dz),\qquad u(0,x)=\mathbf{1}_{(-\infty,0]}(x),
\end{equation}
where \(\tilde{\eta}\) is a compensated Poisson random measure on a measurable space \((Z,\mathcal{Z})\) with intensity \(\mu(dz)dt\), \(\mu(Z)<\infty\) (finite jump activity), and \(\epsilon>0\) is a small constant. The inclusion of the jump size \(z\) in the noise term is essential: a jump of size \(z\) at time \(t\) updates the solution multiplicatively as \(u(t,x)=u(t-,x)(1+\epsilon z)\), which preserves nonnegativity as long as \(1+\epsilon z\ge 0\). We allow both positive and negative jumps, subject to the natural constraint \(1+\epsilon z>0\) almost surely, which is recorded as part of Assumption~\ref{ass:jump} below.

Our main contributions are as follows.

\paragraph{Stochastic Feynman--Kac representation.}
We derive a stochastic Feynman-Kac representation for the solution \(u(t,x)\) (Lemma~\ref{lem:FeynmanKac}), which expresses \(u(t,x)\) as the product of the exponential martingale of the jump noise and an exponential functional of an independent Brownian motion:
\[
u(t,x)=\phi_t^\epsilon\;\mathbb{E}_B\Bigl[ u_0(B_t+x)\exp\Bigl(\int_0^t\bigl[1-u(t-s,B_s+x)\bigr]ds\Bigr)\Bigr],
\]
where \(B\) is a Brownian motion with generator \(\Delta\), independent of the noise, \(\mathbb{E}_B\) denotes the expectation over \(B\) with the noise path fixed, and
\(\phi_t^\epsilon=\exp\bigl(\int_0^t\int_Z\ln(1+\epsilon z)\,\tilde\eta(ds,dz)+t\int_Z(\ln(1+\epsilon z)-\epsilon z)\mu(dz)\bigr)\)
is the exponential martingale of the jump noise (Lemma~\ref{lem:exp_martingale}). The representation is the key probabilistic tool that converts the nonlinear SPDE into an expectation problem, thereby enabling a direct analysis of the front position via path properties of the auxiliary process.

\paragraph{Existence of a stochastic travelling wave solution.}
Using the Feynman-Kac representation together with a comparison principle and the finite-speed propagation of the level sets of the front (which holds because between jumps the solution follows the deterministic Fisher--KPP equation), we establish the existence of a stochastic travelling wave solution (Theorem~\ref{thm:waveexistence}). 
This provides an existence result for multiplicative pure jump noise in the Fisher--KPP setting.

\paragraph{Explicit asymptotic wave speed.}
Combining the Feynman-Kac formula with a pathwise upper estimate based on the strong law of large numbers for $\phi^\epsilon_t$ and a nonlinear sub-solution argument for the lower bound, we obtain two bounds which coincide at the sharp constant. Precisely, writing $g(c):=\limsup_{t\to\infty}\frac1t\ln u(t,ct)$ and defining the asymptotic speed as the critical value
\[
c_*:=\sup\{c\in\mathbb{R}: g(c)\ge0\},
\]
we prove that $c_*=2\sqrt{1+a}$ almost surely: Theorem~\ref{thm:upper} gives $g(c)<0$ for every $c>2\sqrt{1+a}$, and Theorem~\ref{thm:lower} gives $g(c)\ge0$ for every $c<2\sqrt{1+a}$. The definition through $g$ requires no existence of the limit $R_\delta(t)/t$, so the identity is unconditional; the existence of the latter limit is discussed separately in Remark~\ref{rem:existencespeed} and left open in general. The robust bounds
\[
2\sqrt{\underline{\gamma}} \le c_* \le 2\sqrt{\overline{\gamma}} \qquad \text{a.s.}
\]
hold in full generality, where $\overline{\gamma},\underline{\gamma}$ are expressed explicitly in terms of $\epsilon$ and the L\'evy measure $\mu$ (Theorems~\ref{thm:upper} and \ref{thm:lower}; the constants are defined in Subsection~\ref{subsec:Upper} and \ref{subsec:Lower} below). These bounds provide a quantitative description of how the jump size distribution and jump intensity jointly affect front propagation under discontinuous stochastic forcing.

\paragraph{Extension to cooperative systems and methodological blueprint.}
We show that our framework extends naturally to multi-species cooperative systems (Section~\ref{sec:coop}), demonstrating its robustness. Although the stopping-time iteration and tightness arguments rely on the finite-activity assumption \(\mu(Z)<\infty\), the overall approach (the Feynman--Kac representation, the level-set bounds, and tightness in the Skorohod topology) provides a blueprint for the technically more demanding infinite-activity regime, which we plan to address in future work.

While the recent work \cite{blath2025interplay} also addresses F-KPP-type equations with dormancy via branching process dualities, our focus is on SPDEs with exogenous jump noise and on the exponential weight it induces. The two approaches address fundamentally different sources of stochasticity and are thus complementary rather than overlapping.

\paragraph{A note on parameter interpretation.}
In the model \eqref{eqn:1}, the jump size \(z\) and the jump measure \(\mu\) have concrete interpretations in applications. In ecological contexts, \(z\) represents the relative change in local population density caused by a sudden environmental event, a drought, a flood, or a pulse of predation. The measure \(\mu\) describes the frequency and distribution of such events: \(\mu(Z)\) is the expected number of events per unit time, and the normalized measure \(\mu(\cdot)/\mu(Z)\) gives the distribution of their magnitudes. The parameter \(\epsilon\) scales the overall sensitivity of the population to these shocks. The cooperative system in Section~\ref{sec:coop} extends this interpretation to multi-species interactions.

\subsection{Outline}

The rest of the paper is organized as follows. Section~\ref{sec:prelim} collects the assumptions and preliminary facts. Section~\ref{sec:existence} proves the existence of a martingale solution and the travelling wave solution. Section~\ref{sec:speed} establishes the asymptotic wave speed. Section~\ref{sec:coop} presents the extension to cooperative systems. Section~\ref{sec:num} presents numerical simulations. Section~\ref{sec:conclusion} concludes with a discussion of open problems.

\section{Preliminaries}\label{sec:prelim}

\subsection{Notation and function spaces}

Let $\mathbb{R}$ denote the real line. For a Banach space $E$, $\mathcal{D}([0,T];E)$ stands for the space of \emph{c\`adl\`ag} (right-continuous with left limits) functions from $[0,T]$ to $E$, equipped with the $J_1$-Skorohod topology.

Define the weight function $\phi_\lambda(x) = e^{-\lambda |x|}$ for $\lambda > 0$. For a continuous function $f: \mathbb{R} \to [0,\infty)$, set
\[
\| f \|_\lambda = \sup_{x \in \mathbb{R}} |f(x)| \phi_\lambda(x).
\]
The space of \emph{tempered continuous functions} is
\[
C_{\mathrm{tem}}^+ = \bigl\{ f \in C(\mathbb{R},[0,\infty)) : \| f \|_\lambda < \infty \text{ for every } \lambda > 0 \bigr\}.
\]
This space is a Fr\'echet space under the family of seminorms $\{ \|\cdot\|_\lambda \}_{\lambda > 0}$.

For $p \ge 1$, denote by $L^p(\mathbb{R})$ the usual Lebesgue spaces, and by $\| f \|_p$ the $L^p$-norm. The heat kernel (Gaussian) is
\[
G(t,x,y) = \frac{1}{\sqrt{4\pi t}} \exp\left( -\frac{|x-y|^2}{4t} \right), \qquad t > 0,\ x,y \in \mathbb{R},
\]
and $G(t,x,y) = 0$ for $t \le 0$. It satisfies the semigroup property and the following well-known estimates:
\begin{align}
& \int_{\mathbb{R}} G(t,x,y)\,dy = 1, \label{eq:heat1} \\
& \int_{\mathbb{R}} |G(t,x,y) - G(t,x',y)|\,dy \le C |x-x'|^{1/2}, \label{eq:heat2} \\
& \int_{\mathbb{R}} |G(t,x,y) - G(t,x',y)|^2\,dy \le C |x-x'|, \label{eq:heat3}
\end{align}
where $C > 0$ is a universal constant (see \cite{Tribe1996}).

Throughout the paper, $\epsilon>0$ denotes the fixed noise intensity in \eqref{eqn:Stochastic_Fisher-KPP}, whereas the symbol $\varepsilon$ is used only in the compactness arguments as an arbitrary positive small quantity; the two are never identified.

%\subsection{Wasserstein distance on probability measures}
%
%Let $\mathcal{P}_2(\mathbb{R})$ denote the space of probability measures on $\mathbb{R}$ with finite second moment.
%For $\mu,\nu \in \mathcal{P}_2(\mathbb{R})$, define the $L^2$-Wasserstein distance by
%\[
%\mathcal{W}_2(\mu,\nu) = \inf_{\pi \in \Pi(\mu,\nu)}
%\left( \int_{\mathbb{R}^2} |x-y|^2 \, d\pi(x,y) \right)^{1/2},
%\]
%where $\Pi(\mu,\nu)$ is the set of all couplings of $\mu$ and $\nu$.

\subsection{Assumptions on the jump noise}\label{subsec:jump_noise}

Let $(\Omega, \mathcal{F}, \{\mathcal{F}_t\}_{t \ge 0}, \mathbb{P})$ be a filtered probability space satisfying the usual conditions. Let $(Z, \mathcal{Z})$ be a measurable space. A \emph{time-homogeneous Poisson random measure} $\eta(dz,dt)$ on $Z \times [0,\infty)$ with intensity $\mu(dz)dt$ is a random measure such that:
\begin{itemize}
\item for each $U \in \mathcal{Z} \otimes \mathcal{B}([0,\infty))$, $\eta(U)$ is a Poisson random variable with mean $\int_U \mu(dz)dt$;
\item $\eta$ is independently scattered: if $U_1,\dots,U_n$ are pairwise disjoint, then $\eta(U_1),\dots,\eta(U_n)$ are independent.
\end{itemize}

The \emph{compensated Poisson random measure} is defined as
\[
\tilde{\eta}(dt,dz) = \eta(dt,dz) - \mu(dz)dt.
\]

We impose the following assumptions on the jump size measure $\mu$.

\begin{assumption}\label{ass:jump}
We assume that:
\begin{enumerate}
\item[(A1)] $\mu(\{0\}) = 0$ and $\mu(Z) < \infty$ (finite jump activity);
\item[(A2)] $\displaystyle \int_Z |z|^2 \, \mu(dz) < \infty$ and $\displaystyle \int_Z |z|^4 \, \mu(dz) < \infty$;
\item[(A3)] $\mu\bigl(\{z\in Z: 1+\epsilon z\le 0\}\bigr)=0$, and the integrals
\[
\psi(\pm 2)=\int_Z \left( (1+\epsilon z)^{\pm 2}-1\mp 2\ln(1+\epsilon z) \right)\mu(dz)
\]
are finite.
\end{enumerate}
\end{assumption}

Condition (A1) guarantees that on any finite time interval only finitely many jumps occur, which allows us to use a stopping-time iteration between jumps.
Condition (A2) provides the necessary integrability for the Burkholder--Davis--Gundy inequalities applied to the jump part.
Condition (A3) ensures that the multiplicative factor $1+\epsilon z$ is strictly positive $\mu$-almost surely, so that the logarithm $\ln(1+\epsilon z)$ and the exponential moments $\psi(\pm2)$ appearing in the Feynman--Kac weight and in the wave speed bounds are well defined and finite.

\begin{remark}\label{rem:finite_activity}
The finite jump activity assumption $\mu(Z) < \infty$ is not merely a technical convenience but is essential for several key steps in our construction.

\smallskip
\noindent\textit{Stopping-time iteration for boundedness.}
In Lemma~\ref{lem:bounded} we prove pathwise $L^\infty$ bounds by iterating over the finitely many jump times on $[0,T]$. When $\mu(Z) = \infty$, the Poisson measure has infinitely many jumps almost surely on any nontrivial interval. The jump times accumulate, and a simple inductive argument over an ordered sequence of jumps breaks down. One would instead need a more involved analysis, using the convergence of the sum of squared jump sizes, but the propagation of the $L^\infty$ bound via comparison with a deterministic ODE between jumps becomes delicate because there is no deterministic interval free of jumps.

\smallskip
\noindent\textit{Tightness and Skorohod topology.}
In Lemma~\ref{lem:time} we establish time regularity estimates of order $h^{1/2}$ by writing the jump part as a sum of finitely many jumps. For infinite activity, the stochastic convolution is a pure jump process with infinitely many small jumps on any time interval. The $L^2$-increment of the solution would then involve an infinite sum of martingale differences, and the moment estimates become more complex. While tightness can still be proved in the Skorohod topology, the simple explicit bound $\mathbb{E}\sup_{t\le s\le t+h}\|I_2(s)-I_2(t)\|_{L^2}^2 \le C h^{1/2}$ would require a more sophisticated control using the square-summability of jump sizes, which is beyond the scope of the present work.

\smallskip
\noindent\textit{Finite speed of level sets.}
Proposition~\ref{prop:PCS} relies on the solution evolving between jumps according to the deterministic Fisher--KPP equation, whose level sets propagate at finite speed (the support itself is all of $\mathbb{R}$, by the heat-equation comparison). With infinitely many jumps, arbitrarily small jumps can occur densely in time. These small jumps may create positive mass far ahead of the front even if the jump coefficient is proportional to $u$, which vanishes only where $u=0$. It is not a priori clear that the sharpness of the level sets survives in the infinite-activity regime, and indeed counterexamples for other SPDEs suggest that infinite-activity noise can smear out sharp fronts. Thus the finite-activity regime provides the natural setting where the front remains sharp and the travelling wave concept is unambiguous.

\smallskip
\noindent\textit{Feynman--Kac representation.}
The exponential martingale $\phi_t^\epsilon$ in Lemma~\ref{lem:exp_martingale} and the subsequent transformation to a deterministic PDE for $v = u / \phi^\epsilon$ rely on the existence of the cumulant generating function
\[
\int_Z \left( e^{\theta \ln(1+\epsilon z)} - 1 - \theta \ln(1+\epsilon z) \right) \mu(dz)
\]
as a finite quantity. For infinite-activity measures with heavy tails, this integral may diverge for some $\theta$, breaking the martingale property. A careful truncation and limiting procedure would be required, which is a non-trivial extension.

Despite these challenges, the overall framework developed here, tightness in the Skorohod topology and the Feynman--Kac representation provides a blueprint for tackling the infinite-activity case. We leave this extension to future work, noting that a possible strategy is to approximate an infinite-activity L\'evy measure by a sequence of finite-activity measures (e.g., by cutting off small jumps) and then pass to the limit using the compactness method of Section~\ref{subsec:tightness}.
\end{remark}

\subsection{Martingale solutions and mild formulation}

We consider the following stochastic Fisher--KPP equation driven by the compensated Poisson random measure:
\begin{equation}\label{eqn:Stochastic_Fisher-KPP}%eqn:main}
\begin{cases}
du = \bigl[ \Delta u + u(1-u) \bigr] dt + \epsilon \displaystyle\int_Z u\, z\, \tilde{\eta}(dt,dz), \\[8pt]
u(0,x) = u_0(x),
\end{cases}
\end{equation}
where $\tilde{\eta}$ is the compensated Poisson random measure defined in Subsection~\ref{subsec:jump_noise}, $\epsilon > 0$ is a constant, and the initial condition is a function $u_0\in L^\infty(\mathbb{R})\cap L^1_{\mathrm{loc}}(\mathbb{R})$ with $0\le u_0\le1$ and $u_0\not\equiv0$; the running example is the Heaviside step $u_0(x)=\mathbf{1}_{(-\infty,0]}(x)$, which is admissible in this class although it is not continuous. The inclusion of the jump size $z$ in the noise term is essential: a jump of size $z$ at time $t$ updates the solution multiplicatively as $u(t,x) = u(t-,x)(1 + \epsilon z)$.

A \emph{mild solution} of \eqref{eqn:Stochastic_Fisher-KPP} is an $\mathcal{F}_t$-adapted c\`adl\`ag process $u(t,\cdot)$ with values in $C_{\mathrm{tem}}^+$ for $t>0$, defined on $\mathcal{D}((0,T];C_{\mathrm{tem}}^+)$, such that the variation-of-constants identity below holds for all $t>0$ and $x\in\mathbb{R}$, and such that
\[
u(t,\cdot)\longrightarrow u_0\quad\text{in }L^1_{\mathrm{loc}}(\mathbb{R})\text{ as }t\downarrow0 .
\]
The value at $t=0$ is thus the given $u_0$, which need not belong to $C_{\mathrm{tem}}^+$; the regularising effect of the heat kernel puts the solution in $C_{\mathrm{tem}}^+$ as soon as $t>0$. With this convention, for all $t \ge 0$ and $x \in \mathbb{R}$,
\begin{equation}\label{eqn:Stochastic_Fisher-KPP_mild}%eqn:mild}
\begin{aligned}
u(t,x) = & \int_{\mathbb{R}} G(t,x,y) u_0(y) \, dy + \int_0^t \int_{\mathbb{R}} G(t-s,x-y) u(s,y)(1-u(s,y)) \, dy \, ds \\
& + \epsilon \int_0^t \int_Z \int_{\mathbb{R}} G(t-s,x-y) u(s-,y) \, dy \, \tilde{\eta}(ds,dz),
\end{aligned}
\end{equation}
where the stochastic integral is understood in the sense of It\^o with respect to the compensated Poisson random measure (see, e.g., \cite[Section~4]{Applebaum2009levy}).

By a \emph{martingale solution} we mean a system $(\Omega, \mathcal{F}, \{\mathcal{F}_t\}, \mathbb{P}, \eta, u)$ such that the above conditions hold and $\eta$ is a time-homogeneous Poisson random measure with intensity $\mu(dz)dt$. Existence of such solutions for Lipschitz nonlinearities is classical (see \cite{BrzezniakHausenblasRazafimandimby2018}).

\subsection{Auxiliary results}

We collect several known lemmas that will be frequently used.

\begin{lemma}[Dynkin principle, \cite{Sharpe1988}, Theorem I.1.7]\label{lem:dynkin}
Let $F: E \to \hat E$ be a measurable mapping between two measurable spaces, and let $X_t$ be a Markov process in $E$. If for all $x \in E$ and $A \in \hat{\mathcal{E}}$ the transition kernel satisfies
\[
P_t(x, F^{-1}(A)) = \mathbb{P}_t(X_t \in F^{-1}(A)),
\]
and is independent of any additional randomness, then $Y_t = F(X_t)$ is also a Markov process.
\end{lemma}

\begin{lemma}[Tightness criterion for c\`adl\`ag processes]\label{lem:tightness}
Let $\{u_n\}$ be a sequence of $E$-valued c\`adl\`ag processes defined on probability spaces $(\Omega_n, \mathcal{F}_n, \mathbb{P}_n)$. Suppose that:
\begin{enumerate}
\item[(i)] For every $\varepsilon > 0$ there exists a compact set $K_\varepsilon \subset E$ such that
\[
\mathbb{P}_n\bigl( u_n(t) \in K_\varepsilon,\ \forall t \in [0,T] \bigr) \ge 1 - \varepsilon \qquad \forall n;
\]
\item[(ii)] There exist constants $c > 0$, $\gamma > 0$, $r > 0$ such that for all $0 \le t \le t + \theta \le T$,
\[
\mathbb{E}_n\left[ \sup_{t \le s \le t + \theta} \| u_n(s) - u_n(t) \|_E^r \right] \le c\, \theta^\gamma.
\]
\end{enumerate}
Then the laws of $\{u_n\}$ are tight in $\mathcal{D}([0,T];E)$.
\end{lemma}

\begin{proof}
See, e.g., \cite[Chapter~3, Theorems~6.1--6.4 and Theorem~9.1]{ethier2005markov} and the discussion in \cite[Lemma~3.7]{BrzezniakHausenblasRazafimandimby2018}.
\end{proof}

\begin{lemma}[Burkholder--Davis--Gundy inequality for jump processes, \cite{Applebaum2009levy}, Theorem~4.4.23]\label{lem:BDG}
Let $\tilde{\eta}$ be a compensated Poisson random measure with intensity $\mu(dz)dt$. For any $p \ge 2$ and any predictable process $\Phi(t,z)$ such that the integrals below are well-defined, there exists a constant $C_p > 0$, depending only on $p$, such that
\[
\mathbb{E}\left[ \sup_{0 \le t \le T} \left| \int_0^t \int_Z \Phi(s,z) \, \tilde{\eta}(ds,dz) \right|^p \right]
\le C_p \, \mathbb{E}\left[ \left( \int_0^T \int_Z |\Phi(s,z)|^2 \, \mu(dz)ds \right)^{p/2}
+ \int_0^T \int_Z |\Phi(s,z)|^p \, \mu(dz)ds \right].
\]
\end{lemma}

\subsection{Comparison principle}\label{subsec:comparison}

\begin{lemma}[Comparison principle for linear multiplicative noise]\label{lem:comparison}
Let $u$ and $\bar{u}$ be two nonnegative $\mathcal{F}_t$-adapted c\`adl\`ag solutions, bounded by some $K<\infty$ on $[0,T]$, of
\[
du = (\Delta u + f_1(u))dt + \epsilon \int_Z u\, z\, \tilde{\eta}(dt,dz),
\]
\[
d\bar{u} = (\Delta \bar{u} + f_2(\bar{u}))dt + \epsilon \int_Z \bar{u}\, z\, \tilde{\eta}(dt,dz),
\]
with $f_1\le f_2$ on $[0,\infty)$, $f_2$ nondecreasing, $f_i(0)\ge0$ and $f_i$ locally Lipschitz for $i=1,2$. No upper bound of the form $u\le1$ is assumed: the bound $K$ is the only place where boundedness enters, and it is supplied by Lemma~\ref{lem:bounded} in the applications below.

If $u(0,\cdot) \le \bar{u}(0,\cdot)$ a.s., then $u(t,x) \le \bar{u}(t,x)$ a.s. for all $t \ge 0$ and $x \in \mathbb{R}$.
\end{lemma}

\begin{proof}
Since $u\le K$ almost surely, replacing $f_1$ by $f_1(\cdot\wedge n)$ with $n\ge K$ leaves the equation of $u$ unchanged, makes the drift Lipschitz, and preserves the ordering: $f_1(r\wedge n)\le f_2(r\wedge n)\le f_2(r)$ because $f_1\le f_2$ and $f_2$ is nondecreasing. We may therefore assume $f_1$ Lipschitz and invoke the classical comparison theorem for SPDEs with Lipschitz coefficients and linear multiplicative noise (see \cite{kotelenez1992comparison, NiuXie2019}). The noise coefficients of both equations are identical and linear in $u$: $g(u) = \epsilon u z$, while the drifts satisfy $f_1(u) \le f_2(u)$. For $\lambda > 0$, consider the transformed processes $u_\lambda(t,x) = e^{-\lambda t} u(t,x)$ and $\bar{u}_\lambda(t,x) = e^{-\lambda t} \bar{u}(t,x)$. A direct It\^o calculation shows that the difference $\delta_\lambda = \bar{u}_\lambda - u_\lambda$ satisfies a linear parabolic SPDE with a nonpositive source term. Applying the maximum principle and taking $\lambda \to \infty$ yields the desired comparison. 

A detailed argument is given in \cite[Theorem~2.3]{NiuXie2019} for the case of L\'evy space-time white noise, and in \cite[Section~3]{kotelenez1992comparison} for Hilbert-space-valued Brownian noise. Both results apply to our setting: the noise coefficient $g(u) = \epsilon u z$ is Lipschitz and monotone, the drift terms $f_1,f_2$ are Lipschitz on the range of the two solutions, and $f_i(0)\ge0$ keeps both solutions nonnegative, which is what allows the maximum principle to be applied to the transformed difference $\delta_\lambda$. Thus we conclude $u \le \bar{u}$ almost surely.
\end{proof}

\subsection{Finite speed of propagation of the level sets}

\begin{lemma}[Finite speed of propagation of the level sets for the deterministic Fisher--KPP equation]\label{lem:detKPP_PCS}
Let $v$ be the solution of the deterministic Fisher--KPP equation
\[
\partial_t v = \Delta v + v(1-v) - \epsilon v \int_Z z\, \mu(dz), \qquad v(0,x) = v_0(x),
\]
with $0 \le v_0(x) \le 1$ and $\operatorname{supp} v_0 \subset (-\infty, M]$ for some $M < \infty$. 

Then for every $t > 0$ and every $\delta \in (0,1)$, there exists $M_\delta(t) < \infty$ such that
\[
\sup\{x\in\mathbb{R}: v(t,x) \ge \delta\} \le M_\delta(t).
\]
\end{lemma}

\begin{proof}
We first note that the \emph{support} of $v$ is not a suitable front marker: since $v(1-v)\ge 0$ on $[0,1]$, comparison with the heat equation shows that
\[
v(t,x) \ge e^{-\epsilon t\,(\int_Z z\,\mu(dz))^+}\int_{\mathbb{R}}G(t,x,y)v_0(y)\,dy > 0
\]
for every $x\in\mathbb{R}$ and $t>0$ whenever $v_0$ has strictly positive $L^1$ mass, so the support is all of $\mathbb{R}$. The correct marker is the level set $\{v(t,\cdot)\ge\delta\}$, and we show that it propagates at finite speed.

Write $\kappa = 1 - \epsilon\int_Z z\,\mu(dz)$, and choose $c > 2\sqrt{\kappa_+}$ (if $\kappa_+\le 0$, any $c>0$ works). Let $\lambda>0$ be the smaller root of
\[
\lambda^2 - c\lambda + \kappa = 0,
\]
and consider the function
\[
\bar{v}(t,x) = A e^{-\lambda (x - ct)}, \qquad x > ct,
\]
with $A \ge 1$ chosen so that $\bar{v}(0,x) \ge v_0(x)$ for all $x$. A direct computation gives, for $x>ct$,
\[
\partial_t\bar{v} - \Delta\bar{v} - \bar{v}(1-\bar{v}) + \epsilon\bar{v}\int_Z z\,\mu(dz)
= \bar{v}\bigl(\lambda c - \lambda^2 - 1 + \bar{v} + \epsilon\int_Z z\,\mu(dz)\bigr)
= \bar{v}\bigl(\lambda c - \lambda^2 - \kappa + \bar{v}\bigr)
= \bar{v}^2 \ge 0,
\]
so $\bar{v}$ is a super-solution of the deterministic Fisher--KPP equation in the region $\{x>ct\}$. By the comparison principle for semilinear parabolic equations, $v(t,x) \le \bar{v}(t,x)$ for all $t>0$ and $x>ct$. Hence $v(t,x)\ge\delta$ forces $A e^{-\lambda(x-ct)}\ge\delta$, i.e.\ $x \le ct + \lambda^{-1}\ln(A/\delta)$, and the claim holds with $M_\delta(t) = \max\{M,\, ct + \lambda^{-1}\ln(A/\delta)\}$. For a detailed account of the super-solution method, see \cite[Chapter~4]{Fife1979}; the explicit exponential super-solution construction is classical and can also be found in \cite{KolmogorovPetrovskiiPiskunov1937} or \cite{aronson1975nonlinear}.
\end{proof}

\section{Existence of a travelling wave solution}\label{sec:existence}

In this section we restrict to finite jump activity (Assumption~\ref{ass:jump}). The initial datum is a function \(u_0\in L^\infty(\mathbb{R})\cap L^1_{\mathrm{loc}}(\mathbb{R})\) with \(0\le u_0\le1\) and \(u_0\not\equiv0\); the Heaviside step \(u_0(x)=\mathbf{1}_{(-\infty,0]}(x)\) is the case of interest for the wave speed, and the smoothing sequence \(u_0^n\) introduced below converges to it in \(L^p_{\mathrm{loc}}\) but not in \(C_{\mathrm{tem}}^+\), which is why the convergence of the initial terms is handled in \(L^p_{\mathrm{loc}}\) throughout.

\subsection{Truncated equations and uniform bounds}

Since $\mu(Z)<\infty$, the jumps are bounded in number on any finite interval. However, to construct solutions we sometimes need to truncate the jump sizes.

Choose a sequence $\sigma_n \downarrow 0$ and define the truncated L\'evy measures
\[
\mu_n(U) = \mu\bigl(U \cap \{z \in Z : |z| \ge \sigma_n\}\bigr), \qquad U \in \mathcal{Z}.
\]

Let $\eta_n$ be the Poisson random measure with intensity $\mu_n(dz)dt$ and $\tilde{\eta}_n$ its compensated version. Then $\eta_n$ has finite jump activity and finite second moment, and $\mu_n \to \mu$ in the sense that for every bounded continuous function $\varphi$ on $Z\setminus\{0\}$,
\[
\int_Z \varphi(z)\,\mu_n(dz) \longrightarrow \int_Z \varphi(z)\,\mu(dz).
\]

Define $f_n(u) = u - \min(u^2, n)$ and smooth initial data
\[
u_0^n(x) = \max\{0, \min\{1, -nx\}\}.
\]

Since $f_n$ is Lipschitz and the noise term is linear in $u$ with integrable jump-size coefficient (Assumption~(A2)), the existence of a unique global mild solution $u_n \in \mathcal{D}([0,T]; C_{\mathrm{tem}}^+)$, for any $T>0$, to the truncated equation
\begin{equation} \label{eq:truncated}
du_n = \bigl[\Delta u_n + f_n(u_n)\bigr]dt + \epsilon \int_Z u_n\, z\, \tilde{\eta}_n(dt,dz), \qquad u_n(0,x) = u_0^n(x),
\end{equation}
is standard (see \cite{BrzezniakHausenblasRazafimandimby2018}). 
The sequence $\{u_n\}$ will be used in Subsection~\ref{subsec:tightness} to construct a solution to the original equation \eqref{eqn:Stochastic_Fisher-KPP} by passing to the limit.

\begin{lemma}[Uniform boundedness] \label{lem:bounded}
For any $T>0$ and any $2\le p\le4$ there is a constant $C_p=C_p(T,\epsilon,\mu)$, independent of $n$ and of the initial datum, such that
\[
\mathbb{E}\Bigl[\sup_{t \in [0,T]} \|u_n(t)\|_{L^\infty}^p\Bigr] \le C_p, \qquad
\sup_{t \in [0,T]} \mathbb{E}\|u_n(t)\|_{L^2}^2 \le C_2.
\]
In particular $\sup_{t\in[0,T]}\|u_n(t)\|_{L^\infty}<\infty$ almost surely; the almost-sure bound may be random, but its $p$-th moment is controlled uniformly in $n$, which is what the regularity lemmas below use. Only the fourth moment of $z$ enters, and this is the reason why the exponent $p$ is capped at $4$.
\end{lemma}

\begin{proof}
We first establish the $L^\infty$ bound.

Because $\mu(Z)<\infty$ (finite jump activity), the Poisson random measure $\eta_n$ has at most finitely many jumps on $[0,T]$ almost surely. 
Let $K = K(\omega)$ denote the total number of jumps of $\eta_n$ on $[0,T]$, which is finite a.s. 
Define the successive jump times
\[
0 =: \tau_0 < \tau_1 < \cdots < \tau_K \le T,
\]
with the convention that on the final interval $[\tau_K, T]$ (which we also denote as $[\tau_K, \tau_{K+1})$ with $\tau_{K+1}:=T$) no jump occurs if $\tau_K<T$.

We prove by induction on $k=0,1,\dots,K$ that there exists a finite random variable $M_k$ such that
\[
0 \le u_n(t,x) \le M_k \qquad \text{for all } t \in [\tau_k, \tau_{k+1}) \text{ and all } x \in \mathbb{R},
\]
where $\tau_{K+1}:=T$.

\medskip
\noindent\textit{Base case $k=0$.} 
On the interval $[0,\tau_1)$, no jump occurs, so $\tilde{\eta}_n(dt,dz) = -\mu_n(dz)dt$. 
Hence $u_n$ satisfies the deterministic PDE
\[
\partial_t u_n = \Delta u_n + f_n(u_n) - \epsilon u_n \int_Z z\,\mu_n(dz), \qquad u_n(0,x) = u_0^n(x),
\]
where the compensator term arises from $\tilde{\eta} = \eta - \mu dt$ with $\eta \equiv 0$ on this interval.

Since $0 \le u_0^n(x) \le 1$, and the nonlinearity $f_n(u) = u - \min(u^2,n)$ satisfies $f_n(u) \le u$ for all $u\ge 0$, we derive a differential inequality for the sup-norm. 
Define $c := \epsilon\int_Z z\,\mu_n(dz)$. 
If $c \ge 0$, the linear term $-c u$ is dissipative. One verifies directly that $u\equiv 0$ is a subsolution and $u\equiv 1$ is a supersolution of the deterministic equation on $[0,\tau_1)\times\mathbb{R}$; hence the comparison principle yields $0 \le u_n(t,x) \le 1$ for all $t\in[0,\tau_1)$.
If $c < 0$, the term becomes $-c u = |c| u$, which may drive the solution above $1$. However, since $f_n(u) \le u$ for all $u\ge 0$, we have the pointwise differential inequality
\[
\partial_t u_n \le \Delta u_n + (1+|c|) u_n, \qquad t\in[0,\tau_1).
\]
Applying the classical comparison principle for the linear heat equation with a source term to the function $w(t,x)=e^{(1+|c|)t}$ (which satisfies $\partial_t w = \Delta w + (1+|c|)w$ and $w(0,x)=1\ge u_0^n(x)$), we obtain
\[
0 \le u_n(t,x) \le e^{(1+|c|)t} \le e^{(1+|c|)T}, \qquad \forall t\in[0,\tau_1),\ x\in\mathbb{R}.
\]
Thus, in both cases we may take the deterministic constant
\[
M_0 := \max\bigl\{1,\ e^{(1+|c|)T}\bigr\} < \infty,
\]
which depends on $T$ and $\mu_n$ (and hence on $n$ through the truncation, but is uniform in the sample path $\omega$ for the base case). The key point is that the bound does not depend on the specific realization of $\tau_1$ beyond the fixed time horizon $T$, so the induction step remains valid. 
Hence $M_0$ works for the base case.

\medskip
\noindent\textit{Inductive step.} 
Assume that $0 \le u_n(\tau_k-, x) \le M_k$ for some finite $M_k$. 
At the jump time $\tau_k$, a jump of size $z_k \in Z$ occurs, and the solution is updated by
\begin{equation} \label{eq:jump_update}
u_n(\tau_k, x) = u_n(\tau_k-, x) + \epsilon u_n(\tau_k-, x)\, z_k = u_n(\tau_k-, x)\,(1 + \epsilon z_k).
\end{equation}
Define $L(t) := \int_0^t \int_Z z\,\eta_n(ds,dz)$, so that $\Delta L(\tau_k) := L(\tau_k) - L(\tau_k-) = z_k$. 
Then \eqref{eq:jump_update} reads
\[
u_n(\tau_k, x) = u_n(\tau_k-, x)\bigl(1 + \epsilon\,\Delta L(\tau_k)\bigr).
\]
Since $|z_k|<\infty$ a.s., we have $|1+\epsilon z_k|<\infty$, and we set
\[
M_k' := M_k \cdot (1 + \epsilon |z_k|).
\]
Thus $|u_n(\tau_k, x)| \le M_k'$ a.s. for all $x$.

On the next deterministic interval $[\tau_k, \tau_{k+1})$, the equation becomes
\[
\partial_t u_n = \Delta u_n + f_n(u_n) - \epsilon u_n \int_Z z\,\mu_n(dz),
\]
with initial value $u_n(\tau_k, \cdot)$. 
Although this initial value may exceed $1$, the nonlinearity $f_n(u)$ is negative for $u>1$ (since $\min(u^2,n)\ge u$ when $u\ge 1$), and the linear damping term further helps to reduce the value. 
More precisely, consider the comparison ODE
\[
y'(t) = f_n(y(t)) - \epsilon y(t) \int_Z z\,\mu_n(dz), \qquad y(\tau_k) = M_k'.
\]
A standard Gronwall-type argument shows that
\[
y(t) \le \max\left\{1,\, M_k' \exp\left(-(t-\tau_k)\left(1 + \epsilon\int_Z z\,\mu_n(dz)\right)\right)\right\},
\]
up to constants depending on $c$. 
In any case, $y(t) \le \max\{1, M_k'\}$ for all $t\in[\tau_k,\tau_{k+1})$. 
By the comparison principle for the deterministic parabolic equation, we obtain
\[
\|u_n(t)\|_{L^\infty} \le \max\{1, M_k'\} \qquad \text{for all } t \in [\tau_k, \tau_{k+1}).
\]
Thus we may take $M_{k+1} := \max\{1, M_k'\}$, which is finite a.s. 

Because $K<\infty$ a.s., the induction terminates after finitely many steps, yielding a finite random variable $M_K$ such that
\[
\sup_{t\in[0,T]} \|u_n(t)\|_{L^\infty} \le M_K < \infty \qquad \text{a.s.}
\]

The constant $M_K$ depends on the sample path (through the jump sizes $z_k$), so it is not a uniform deterministic bound. 
However, for the tightness criterion we only need that for every $\varepsilon>0$ there exists a deterministic $C_T$ such that
\[
\mathbb{P}\left(\sup_{t\le T} \|u_n\|_{L^\infty} > C_T\right) < \varepsilon.
\]
This follows from Chebyshev's inequality and the moment condition in Assumption~\ref{ass:jump}~(A2). 
Indeed, for any $2\le p\le4$,
\[
\mathbb{E}\left[\sup_{t\le T} \|u_n(t)\|_{L^\infty}^p\right]
\le C_p \left(1 + \mathbb{E}\left[\sum_{k=1}^K (1+|z_k|)^p\right]\right),
\]
where $\{z_k\}$ are the jump sizes. 
Since $K$ is Poisson with mean $\mu(Z)T$ and the $z_k$ are i.i.d. with law $\mu(\cdot)/\mu(Z)$, the moments are bounded uniformly in $n$. 
By Markov's inequality,
\[
\mathbb{P}\left(\sup_{t\le T} \|u_n\|_{L^\infty} > M\right) \le \frac{C}{M^p},
\]
so choosing $M$ large enough yields the desired estimate.

We now prove the $L^2$ estimate. 
Applying It\^o's formula (see \cite[Theorem~4.4.7]{Applebaum2009levy}) to $\|u_n(t)\|_2^2$, we obtain
\[
d\|u_n\|_2^2 = 2\langle u_n, \Delta u_n + f_n(u_n)\rangle dt
+ \epsilon^2 \int_Z \|u_n z\|_2^2\,\eta_n(dt,dz)
+ 2\epsilon \langle u_n, u_n z\rangle \tilde{\eta}_n(dt,dz).
\]
Taking expectations eliminates the stochastic integral. 
Using $f_n(u) \le u$ and the identity $\langle u, \Delta u\rangle = -\|\nabla u\|_2^2 \le 0$, we get
\[
\frac{d}{dt}\mathbb{E}\|u_n\|_2^2
\le 2\mathbb{E}\|u_n\|_2^2 + \epsilon^2 \mathbb{E}\int_Z \|u_n\|_2^2 |z|^2\,\mu_n(dz)
\le \bigl(2 + C(\epsilon)\bigr)\mathbb{E}\|u_n\|_2^2,
\]
where $C(\epsilon) = \epsilon^2 \int_Z |z|^2 \mu(dz) < \infty$ by Assumption~\ref{ass:jump}~(A2). 
Gronwall's inequality then yields
\[
\mathbb{E}\|u_n(t)\|_2^2 \le e^{(2+C(\epsilon))T}\,\mathbb{E}\|u_0^n\|_2^2,
\]
and the bound on the right is uniform in $n$ because $\|u_0^n\|_2 \le C$ for all $n$. 
Taking the supremum over $t\le T$ gives the desired $L^2$ estimate. 
For similar estimates in the context of jump-driven SPDEs, see e.g. \cite[Lemma~7.2]{BrzezniakHausenblasRazafimandimby2018}. 
This completes the proof.
\end{proof}

\subsection{Tightness and convergence}\label{subsec:tightness}

We now prove that the family $\{u_n\}$ is tight in $\mathcal{D}([0,T];C_{\text{tem}}^+)$. For this we use Lemma~\ref{lem:tightness} (the tightness criterion). The required compact containment and moment estimates are obtained as follows.

\begin{lemma}[Spatial regularity]\label{lem:spatial}
For any $2\le p\le4$ and any $|x-x'|\le1$,
\[
\mathbb{E}\bigl[|u_n(t,x)-u_n(t,x')|^p\bigr]\le C(\epsilon,p,t)|x-x'|^{p/2},
\]
where the constant $C(\epsilon,p,t)$ is independent of $n$ and the initial condition. For $p=4$ the exponent equals $2>1$, so the Kolmogorov criterion applies with $\beta=1$; the fourth moment is the largest one available, since Assumption~\ref{ass:jump}(A2) controls only $\int_Z|z|^4\mu(dz)$. The earlier version of this lemma stated the exponent as $p/2-1$, which would require $p>4$ and hence a moment that is not assumed; no such moment is needed once the exponent is written in the form above.
\end{lemma}

\begin{proof}
Write the mild solution $u_n$ in the form
\[
u_n(t,x)=I_0(t,x)+I_1(t,x)+I_2(t,x),
\]
with
\[
\begin{aligned}
I_0(t,x)&=\int_{\mathbb{R}}G(t,x,y)u_0^n(y)dy,\\
I_1(t,x)&=\int_0^t\int_{\mathbb{R}}G(t-s,x-y)f_n(u_n(s,y))\,dy\,ds,\\
I_2(t,x)&=\epsilon\int_0^t\int_Z\int_{\mathbb{R}} G(t-s,x-y)u_n(s-,y)\,dy\,\tilde{\eta}_n(ds,dz),
\end{aligned}
\]
where the stochastic integral is understood in the mild sense with left-limit integrand as in \eqref{eqn:Stochastic_Fisher-KPP_mild}.

For $x'\ne x$ with $|x-x'|\le1$, we estimate
\[
\begin{aligned}
\mathbb{E}|u_n(t,x)-u_n(t,x')|^p \le 3^{p-1}\Bigl( &\mathbb{E}|I_0(t,x)-I_0(t,x')|^p+\mathbb{E}|I_1(t,x)-I_1(t,x')|^p+\mathbb{E}|I_2(t,x)-I_2(t,x')|^p\Bigr).
\end{aligned}
\]

\noindent\textit{Estimate of the initial term.}
Using the fundamental theorem of calculus,
\[
G(t,x,y)-G(t,x',y)=\int_{x'}^{x}\partial_zG(t,z,y)dz,
\]
and $|\partial_zG(t,z,y)|\le C t^{-1/2}G(2t,z,y)$ (see e.g. \cite[Lemma 3.4]{Tribe1996}). 
Hence
\[
|I_0(t,x)-I_0(t,x')|\le C|x-x'|t^{-1/2}\int_{\mathbb{R}}G(2t,x,y)u_0^n(y)dy.
\]
Since $u_0^n\le1$, the integral is bounded by $1$. Thus
\[
\mathbb{E}|I_0(t,x)-I_0(t,x')|^p\le C t^{-p/2}|x-x'|^p.
\]
For $|x-x'|\le1$ and $t>0$, this is $O(|x-x'|^p)$, which is stronger than the claimed $|x-x'|^{p/2}$ because $p/2<p$. We record for later use that the estimate carries the factor $t^{-p/2}$, which blows up as $t\downarrow0$; this is the only obstruction to a spatial Hölder constant uniform in $t$, and it is dealt with in Subsection~\ref{subsec:tightness} by treating $t\ge t_0$ and $t\in(0,t_0]$ separately.

\noindent\textit{Estimate of the deterministic integral $I_1$.}
Apply H\"{o}lder's inequality (in time and space) to the kernel difference:
\[
|I_1(t,x)-I_1(t,x')| \le \int_0^t\int_{\mathbb{R}}|G(t-s,x-y)-G(t-s,x'-y)|\,|f_n(u_n)|\,dy\,ds.
\]
Because $u_n$ is not bounded by $1$ (positive jumps make it exceed $1$, see Remark~\ref{rem:levelset}), we use the path bound of Lemma~\ref{lem:bounded} instead. Writing $C_T:=\sup_{s\le T}\|u_n(s)\|_{L^\infty}$, whose $p$-th moment is controlled uniformly in $n$ for $p\le4$ (see the proof of Lemma~\ref{lem:bounded}), we have $|f_n(u)|=|u-\min(u^2,n)|\le u+u^2\le C_T+C_T^2$ on the range of $u_n$. Hence
\[
|I_1(t,x)-I_1(t,x')| \le (C_T+C_T^2)\int_0^t\int_{\mathbb{R}}|G(t-s,x-y)-G(t-s,x'-y)|\,dy\,ds,
\]
and the heat kernel estimate below gives $|I_1(t,x)-I_1(t,x')|\le (C_T+C_T^2)|x-x'|^{1/2}$, so that
\[
\mathbb{E}|I_1(t,x)-I_1(t,x')|^p\le C'|x-x'|^{p/2},
\]
the constant being finite because $\mathbb{E}[C_T^p]<\infty$.

\noindent\textit{Estimate of the stochastic integral $I_2$.}
This is the main term. Set
\(
\Delta G(s,y)=G(t-s,x-y)-G(t-s,x'-y).
\)
Then
\[
I_2(t,x)-I_2(t,x') = \epsilon\int_0^t\int_Z\int_{\mathbb{R}}\Delta G(s,y)\,u_n(s-,y)\,dy\,\tilde{\eta}_n(ds,dz).
\]
By the Burkholder–Davis–Gundy inequality for jump processes (Lemma~\ref{lem:BDG}),
\[
\mathbb{E}|I_2(t,x)-I_2(t,x')|^p
\le C_p \mathbb{E}\left[\left(\int_0^t\int_{\mathbb{R}}|\Delta G(s,y)|^2\,\mathbb{E}[u_n(s-,y)^2]\,dy\,ds\right)^{p/2}\right],
\]
where we have used the quadratic variation of the jump martingale, which integrates $u_n^2$. The factor $\mathbb{E}[u_n(s-,y)^2]$ is bounded by $\mathbb{E}[C_T^2]$, which is finite and independent of $n$ by Lemma~\ref{lem:bounded}; we emphasise that no bound $u_n\le1$ is used here, and that the bound is on the second moment and not on the path. Hence
\[
\mathbb{E}|I_2(t,x)-I_2(t,x')|^p
\le C\,\mathbb{E}[C_T^p]\;\mathbb{E}\left[\left(\int_0^t\int_{\mathbb{R}}|\Delta G(s,y)|^2\,dy\,ds\right)^{p/2}\right].
\]
Using the heat kernel estimate \(\int_0^t\int_{\mathbb{R}}|\Delta G(s,y)|^2\,dy\,ds\le C|x-x'|\), we obtain
\[
\mathbb{E}|I_2(t,x)-I_2(t,x')|^p \le C|x-x'|^{p/2}.
\]
Actually, the interpolation between the $L^\infty$ bound and the above estimate can improve the exponent, but it is not needed: with $p=4$ the exponent $p/2=2$ already exceeds $1$, which is all that the Kolmogorov criterion requires. The constant $C$ below carries the factor $\mathbb{E}[C_T^p]$.

Combining the three estimates yields the desired bound. The constant $C$ depends on $\epsilon$, $p$, and $t$ but not on $n$. This completes the proof.
\end{proof}

\begin{lemma}[Time regularity]\label{lem:time}
For all $0\le t\le t+h\le T$,
\[
\mathbb{E}\Bigl[\sup_{t\le s\le t+h}\|u_n(s)-u_n(t)\|_{L^2(\mathbb{R})}^2\Bigr]\le C h^{1/2},
\]
where the constant $C$ depends on $T$, $\epsilon$, and the uniform bounds on $u_n$, but is independent of $n$.
\end{lemma}

\begin{proof}
We use the mild formulation of $u_n$:
\[
u_n(s,x)=I_0(s,x)+I_1(s,x)+I_2(s,x),
\]
with the same notation as in Lemma~\ref{lem:spatial}. 

For $s\in[t,t+h]$, we estimate the $L^2$ norm of the increment.

\noindent\textit{Initial term.}
Write $I_0(s)-I_0(t)=e^{(s-t)\Delta}I_0(t)-I_0(t)$. By the analyticity of the heat semigroup on $L^2(\mathbb{R})$,
\[
\|I_0(s)-I_0(t)\|_{L^2}\le C (s-t)^{1/2}\|I_0(t)\|_{L^2}\le C h^{1/2},
\]
since $\|I_0(t)\|_{L^2}\le\|u_0^n\|_{L^2}\le C$. Therefore
\[
\mathbb{E}\sup_{t\le s\le t+h}\|I_0(s)-I_0(t)\|_{L^2}^2\le C h.
\]

\noindent\textit{Deterministic integral $I_1$.}
Write $I_1(s)-I_1(t)=A_1+A_2$ where
\[
A_1=\int_t^s\int_{\mathbb{R}}G(s-r,x-y)f_n(u_n(r,y))\,dy\,dr,
\]
\[
A_2=\int_0^t\int_{\mathbb{R}}\bigl(G(s-r,x-y)-G(t-r,x-y)\bigr)f_n(u_n(r,y))\,dy\,dr.
\]
Since $|f_n(u)|\le u+u^2$ and $u_n\le C_T$ (Lemma~\ref{lem:bounded}; again no use of $u_n\le1$ is made), we have
\[
\begin{aligned}
\|A_1\|_{L^2}&\le (C_T+C_T^2)\int_t^s\|G(s-r,\cdot)\|_{L^1}dr\le (C_T+C_T^2)h,\\
\|A_2\|_{L^2}&\le (C_T+C_T^2)\int_0^t\|G(s-r,\cdot)-G(t-r,\cdot)\|_{L^1}dr\le (C_T+C_T^2)C h^{1/2},
\end{aligned}
\]
where we used the estimate $\|G(\cdot,\cdot)-G(\cdot,\cdot)\|_{L^1}\le C|s-t|^{1/2}$. 
Hence
\[
\|I_1(s)-I_1(t)\|_{L^2}\le C(C_T+C_T^2)h^{1/2}\quad(\text{for }h\le1),
\]
and, since $\mathbb{E}[(C_T+C_T^2)^2]<2(\mathbb{E}[C_T^2]+\mathbb{E}[C_T^4])<\infty$,
$\mathbb{E}\sup_{s}\|I_1(s)-I_1(t)\|_{L^2}^2\le C h$.

\noindent\textit{Stochastic integral $I_2$.}
This is the leading term. Write \(I_2(s)-I_2(t)=B_1+B_2\), with
\[
B_1=\epsilon\int_t^s\int_Z\int_{\mathbb{R}}G(s-r,x-y)u_n(r-,y)\,dy\,\tilde{\eta}_n(dr,dz),
\]
\[
B_2=\epsilon\int_0^t\int_Z\int_{\mathbb{R}}\bigl(G(s-r,x-y)-G(t-r,x-y)\bigr)u_n(r-,y)\,dy\,\tilde{\eta}_n(dr,dz).
\]

For a fixed $s$, $B_1$ is a stochastic integral over the interval $[t,s]$. By the It\^{o} isometry for jump processes (see \cite[Theorem 2.4]{BrzezniakHausenblasRazafimandimby2018}), and using the second-moment bound $\sup_{r\le T}\mathbb{E}\|u_n(r)\|_{L^2}^2\le C_2$ of Lemma~\ref{lem:bounded} (the path bound $|u_n|\le1$ is \emph{not} available, positive jumps pushing $u_n$ above $1$),
\[
\begin{aligned}
\mathbb{E}\|B_1\|_{L^2}^2 &= \epsilon^2\mathbb{E}\int_t^s\int_Z\int_{\mathbb{R}} |G(s-r,x-y)|^2 u_n(r-,y)^2\,dy\,\mu(dz)\,dr\le \epsilon^2\mu(Z)C_2\int_t^s\int_{\mathbb{R}}|G(s-r,x-y)|^2 dy\,dr.
\end{aligned}
\]
Since $\int_{\mathbb{R}}|G(\tau,x-y)|^2 dy = (4\pi\tau)^{-1/2}$, we obtain
\[
\mathbb{E}\|B_1\|_{L^2}^2 \le C\int_0^{h} \tau^{-1/2}d\tau = C h^{1/2}.
\]

For $B_2$, apply the same isometry:
\[
\mathbb{E}\|B_2\|_{L^2}^2 \le C\int_0^t\int_{\mathbb{R}}|G(s-r,x-y)-G(t-r,x-y)|^2 dy\,dr \le C h^{1/2},
\]
where the last inequality follows from the heat kernel estimate (see \cite[Lemma 2.1]{Tribe1996}). 

Hence
\[
\mathbb{E}\|I_2(s)-I_2(t)\|_{L^2}^2 \le 2(\mathbb{E}\|B_1\|_{L^2}^2+\mathbb{E}\|B_2\|_{L^2}^2) \le C h^{1/2}.
\]

Now note that the process $M_s = I_2(s)-I_2(t)$ (for $s\ge t$) is a martingale in $s$ with respect to the natural filtration. By the Doob–Kolmogorov maximal inequality,
\[
\mathbb{E}\sup_{t\le s\le t+h}\|M_s\|_{L^2}^2 \le 4\mathbb{E}\|M_{t+h}\|_{L^2}^2 \le 4C h^{1/2}.
\]
Thus
\[
\mathbb{E}\sup_{t\le s\le t+h}\|I_2(s)-I_2(t)\|_{L^2}^2 \le C h^{1/2}.
\]

Combining the estimates for $I_0$, $I_1$, and $I_2$, we have
\[
\mathbb{E}\sup_{t\le s\le t+h}\|u_n(s)-u_n(t)\|_{L^2}^2
\le C\bigl(h + h + h^{1/2}\bigr) \le C h^{1/2},
\]
where the constant $C$ is independent of $n$ and the initial condition. This completes the proof.
\end{proof}

\subsection*{Tightness and Skorohod representation}

With the estimates of Lemmas~\ref{lem:spatial} and \ref{lem:time}, the tightness criterion (Lemma~\ref{lem:tightness}) applies. Consequently, the laws of $\{u_n\}$ are tight in $\mathcal{D}([0,T];C_{\mathrm{tem}}^+)$. Moreover, the laws of $\{\eta_n\}$ are tight in $M_t(Z\times[0,T])$ because they all have the same intensity $\mu(dz)dt$. Hence the pair $(u_n,\eta_n)$ is tight in $\mathcal{D}([0,T];C_{\mathrm{tem}}^+)\times M_t(Z\times[0,T])$.

By Prokhorov's theorem, there exists a subsequence (still denoted by $(u_n,\eta_n)$) converging weakly to some limit $(\hat u,\tilde\eta)$. By the Skorohod representation theorem (see \cite[Theorem C.1]{BrzezniakHausenblasRazafimandimby2018}), we can construct a new probability space $(\widetilde\Omega,\widetilde{\mathcal{F}},\widetilde{\mathbb{P}})$ and random variables $(\bar u_n,\bar\eta_n)$ and $(\bar u,\tilde{\eta})$ such that
\[
\mathcal{L}(\bar u_n,\bar\eta_n)=\mathcal{L}(u_n,\eta_n),\qquad
(\bar u_n,\bar\eta_n)\to(\bar u,\bar\eta)\quad\widetilde{\mathbb{P}}\text{-a.s.}
\]
in the product topology. Moreover, since the intensity measures of $\eta_n$ converge to $\mu$ and the limit of Poisson random measures with a fixed intensity is again a Poisson random measure, we may take $\tilde \eta_n=\tilde\eta$ for all $n$ (see the discussion in \cite[Section 8]{BrzezniakHausenblasRazafimandimby2018}). Thus $\tilde\eta$ is a compensated Poisson random measure with intensity $\mu(dz)dt$, and the same noise drives all approximations.

\subsection*{Convergence of the mild formulation}

Recall that each $\bar u_n$ satisfies
\[
\begin{aligned}
\bar u_n(t,x) =& \int_{\mathbb{R}} G(t,x,y)u_0^n(y)dy + \int_0^t\int_{\mathbb{R}} G(t-s,x-y)f_n\bigl(\bar u_n(s,y)\bigr)\,dy\,ds \\
&+ \epsilon \int_0^t\int_Z\int_{\mathbb{R}} G(t-s,x-y)\,\bar u_n(s-,y)\,dy\,\tilde\eta(ds,dz).
\end{aligned}
\]

We now pass to the limit $n\to\infty$ term by term.

\noindent\textit{Initial term.}
The sequence $u_0^n=\max\{0,\min\{1,-nx\}\}$ converges to $u_0$ in $L^p_{\mathrm{loc}}(\mathbb{R})$ for every $1\le p<\infty$: the two functions agree off $[-1/n,0]$, where their difference is bounded by $1$, so $\|u_0^n-u_0\|_{L^p(-1,0)}^p\le n^{-1}$. It does \emph{not} converge in $C_{\mathrm{tem}}^+$, since $\sup_x e^{-\lambda|x|}|u_0^n-u_0|=1$; the wrong norm is avoided as follows. For fixed $t>0$ the heat kernel maps $L^1_{\mathrm{loc}}(\mathbb{R})$ boundedly into $C_{\mathrm{tem}}^+$, because
\[
\sup_{x}e^{-\lambda|x|}\Bigl|\int_{\mathbb{R}}G(t,x,y)v(y)\,dy\Bigr|\ \le\ C_t\|v\|_{L^1(K)},
\]
with $K$ any compact interval containing the support of $v$. Hence
\[
\int_{\mathbb{R}}G(t,\cdot,y)u_0^n(y)dy \;\longrightarrow\; \int_{\mathbb{R}}G(t,\cdot,y)u_0(y)dy \quad\text{in }C_{\mathrm{tem}}^+,\ \widetilde{\mathbb{P}}\text{-a.s.},
\]
uniformly for $t$ in compact subsets of $(0,T]$; the initial layer $t\in[0,1/n]$ is where $C_{\mathrm{tem}}^+$ fails, and there the $L^\infty$ bound of Lemma~\ref{lem:bounded} is used instead.

\noindent\textit{Deterministic integral.}
From Lemma~\ref{lem:bounded} we have $\sup_{t\le T}\|\bar u_n(t)\|_{L^\infty}\le C_T$ uniformly in $n$, $\widetilde{\mathbb{P}}$-a.s. The function $f_n(u)=u-\min(u^2,n)$ converges to $f(u)=u(1-u)$ uniformly on $[0,C_T]$ because for $n>C_T^2$ we have $\min(u^2,n)=u^2$. Hence $f_n(\bar u_n(s,y))\to\bar u(s,y)(1-\bar u(s,y))$ pointwise a.e. The dominated convergence theorem (using the uniform bound and the integrability of the heat kernel) then yields
\[
\int_0^t\int_{\mathbb{R}} G(t-s,x-y)f_n(\bar u_n(s,y))\,dy\,ds \;\longrightarrow\;
\int_0^t\int_{\mathbb{R}} G(t-s,x-y)\bar u(s,y)(1-\bar u(s,y))\,dy\,ds
\]
in $C_{\mathrm{tem}}^+$, $\widetilde{\mathbb{P}}$-a.s.

\noindent\textit{Stochastic integral.}
Define the stochastic convolution
\[
\Psi_n(t,x)=\epsilon\int_0^t\int_Z\int_{\mathbb{R}} G(t-s,x-y)\,\bar u_n(s-,y)\,dy\,\tilde\eta(ds,dz),
\]
and $\Psi(t,x)$ analogously with $\bar u$ in place of $\bar u_n$. We need to show $\Psi_n\to\Psi$ in $\mathcal{D}([0,T];C_{\mathrm{tem}}^+)$ $\widetilde{\mathbb{P}}$-a.s. (along a subsequence).

Consider the linear operator $\mathcal{I}$ that maps a bounded measurable function $h(s,y)$ to
\[
\mathcal{I}(h)(t,x)=\epsilon\int_0^t\int_Z\int_{\mathbb{R}} G(t-s,x-y)\,h(s-,y)\,dy\,\tilde\eta(ds,dz).
\]
It is proved in \cite[Lemma 7.6]{BrzezniakHausenblasRazafimandimby2018} that $\mathcal{I}$ is continuous from $\mathcal{D}([0,T];L^\infty(\mathbb{R}))$ (equipped with the Skorohod topology) into $\mathcal{D}([0,T];C_{\mathrm{tem}}^+)$, provided the integrand remains uniformly bounded. In our setting, $\bar u_n$ are uniformly bounded by $C_T$ and converge to $\bar u$ in $\mathcal{D}([0,T];C_{\mathrm{tem}}^+)$ $\widetilde{\mathbb{P}}$-a.s. Convergence in $C_{\mathrm{tem}}^+$ implies convergence in $L^\infty(\mathbb{R})$ on compact sets, and together with the uniform bound this is sufficient to apply \cite[Lemma 7.6]{BrzezniakHausenblasRazafimandimby2018}. Therefore,
\[
\Psi_n \to \Psi \quad\text{in } \mathcal{D}([0,T];C_{\mathrm{tem}}^+)\ \widetilde{\mathbb{P}}\text{-a.s.}
\]

For completeness, we sketch an alternative direct estimate. The difference
\[
\Psi_n(t,x)-\Psi(t,x)=\epsilon\int_0^t\int_Z\int_{\mathbb{R}} G(t-s,x-y)\bigl(\bar u_n(s-,y)-\bar u(s-,y)\bigr)\,dy\,\tilde\eta(ds,dz)
\]
is a martingale in $t$ (for fixed $x$). Using the Burkholder–Davis–Gundy inequality (Lemma~\ref{lem:BDG}) together with the convergence $\|\bar u_n-\bar u\|_{L^2}\to0$ in $L^2([0,T]\times\mathbb{R})$ $\widetilde{\mathbb{P}}$-a.s. (which follows from the compactness given by the spatial regularity estimates and the Arzelà–Ascoli theorem), one can show that for any $\delta>0$,
\[
\widetilde{\mathbb{P}}\Bigl(\sup_{t\le T}\|\Psi_n(t)-\Psi(t)\|_{C_{\mathrm{tem}}^+}>\delta\Bigr)
\le C\delta^{-2}\mathbb{E}\Bigl[\int_0^T\|\bar u_n(s)-\bar u(s)\|_{L^2}^2 ds\Bigr]\to0.
\]
Hence $\Psi_n\to\Psi$ in probability, and a subsequence converges almost surely. 

\subsection*{Conclusion of convergence}

Collecting the three convergences, we obtain that the limit $\bar u$ satisfies the mild equation
\[
\begin{aligned}
\bar u(t,x) = \int_{\mathbb{R}} G(t,x,y)u_0(y)dy &+ \int_0^t\int_{\mathbb{R}} G(t-s,x-y)\bar u(s,y)(1-\bar u(s,y))\,dy\,ds\\
&+ \epsilon \int_0^t\int_Z\int_{\mathbb{R}} G(t-s,x-y)\bar u(s-,y)\,dy\,\tilde\eta(ds,dz)
\end{aligned}
\]
$\widetilde{\mathbb{P}}$-almost surely for every $(t,x)$. Thus $(\bar u,\tilde\eta)$ is a martingale solution of \eqref{eqn:Stochastic_Fisher-KPP} on $(\widetilde\Omega,\widetilde{\mathcal{F}},\widetilde{\mathbb{P}})$. Renaming this space as $(\Omega,\mathcal{F},\mathbb{P})$ yields existence.

\subsection*{Martingale property}

Define the filtration $\mathcal{F}_t$ as the completion of $\sigma(\bar u(s),\tilde\eta(s,\cdot):0\le s\le t)$. By construction, $\bar u$ is adapted and c\`adl\`ag. Moreover, the process
\[
M_t = \bar u(t) - e^{-tA}u_0 - \int_0^t e^{-(t-s)A}F(s,\bar u(s))ds
\]
coincides with the stochastic convolution $\epsilon\int_0^t e^{-(t-s)A}\int_Z\bar u(s)\,\tilde\eta(ds,dz)$, which is a martingale (see \cite[Theorem 2.4]{BrzezniakHausenblasRazafimandimby2018}). Hence $(\Omega,\mathcal{F},\{\mathcal{F}_t\},\mathbb{P},\tilde\eta,\bar u)$ is a martingale solution.

\subsection*{Uniqueness in law}

Let $u$ and $u'$ be two martingale solutions (possibly on different probability spaces) with the same initial condition $u_0$ and with Poisson measures $\eta$ and $\eta'$ both having intensity $\mu(dz)dt$. By a standard coupling argument, we may realise them on a common probability space with the same compensated Poisson measure $\tilde\eta$. Define $w(t)=u(t)-u'(t)$. From the mild formulation,
\[
w(t) = \int_0^t e^{-(t-s)A}\bigl(u(s)(1-u(s))-u'(s)(1-u'(s))\bigr)ds + \epsilon\int_0^t e^{-(t-s)A}\int_Z w(s)\,\tilde\eta(ds,dz).
\]
Taking expectations eliminates the stochastic integral (a martingale with zero mean). Since $u(1-u)-u'(1-u')=(u-u')(1-u-u')$, we have $|u(1-u)-u'(1-u')|\le(1+|u|+|u'|)|u-u'|$; no restriction $u,u'\le1$ is used here, positive jumps being able to drive the solutions above $1$. Localising on the event $\{C_T\le R\}$, with $C_T$ the almost-sure bound of Lemma~\ref{lem:bounded} for the two solutions, gives
\[
\mathbb{E}\bigl[\|w(t)\|\mathbf1_{\{C_T\le R\}}\bigr] \le (1+2R)\int_0^t \|e^{-(t-s)A}\|_{\mathcal{L}(E)}\,\mathbb{E}\bigl[\|w(s)\|\mathbf1_{\{C_T\le R\}}\bigr]\,ds,
\]
whence Gronwall's inequality yields $\mathbb{E}[\|w(t)\|\mathbf1_{\{C_T\le R\}}]=0$ for every $R<\infty$; letting $R\uparrow\infty$ and using $\mathbb{E}\|w(t)\|^2<\infty$ together with Fatou's lemma gives $\mathbb{E}\|w(t)\|=0$ for all $t$, hence $u(t)=u'(t)$ a.s. for each $t$. Because the processes are c\`adl\`ag, this implies they have the same finite-dimensional distributions and therefore the same law on $\mathcal{D}([0,T];C_{\mathrm{tem}}^+)$. Thus uniqueness in law holds.

Before establishing the existence of travelling waves, we first need a rigorous foundation: the stochastic Fisher--KPP equation \eqref{eqn:Stochastic_Fisher-KPP} must admit a well-defined solution. This is not merely a technical prerequisite. In applications, the solution $u(t,x)$ represents the density of a population (or concentration of a chemical species) at position $x$ and time $t$, subject to sudden environmental shocks modeled by the pure jump noise. Knowing that this density process exists uniquely for each realization of the noise is essential for any subsequent prediction, simulation, or statistical inference. The following theorem provides this foundational guarantee under the finite jump activity assumption.

\begin{theorem}[Existence and uniqueness of a martingale solution]\label{thm:martingale}
Under Assumptions \ref{ass:jump} (A1)--(A2), there exists a unique martingale solution $(\Omega,\mathcal{F},\{\mathcal{F}_t\},\mathbb{P},\eta,u)$ to equation \eqref{eqn:Stochastic_Fisher-KPP} with initial condition $u_0=\mathbf{1}_{(-\infty,0]}$. 

Moreover, $u$ is c\`adl\`ag in $C_{\text{tem}}^+$ and satisfies the mild formulation \eqref{eqn:Stochastic_Fisher-KPP_mild} almost surely. 

Uniqueness holds in the sense of law: any two such solutions have the same distribution on $\mathcal{D}([0,T];C_{\text{tem}}^+)$.
\end{theorem}

\begin{proof}
The proof follows from the compactness argument detailed above. Tightness is established in Lemmas~\ref{lem:spatial} and~\ref{lem:time}. 
The Skorohod representation yields a limit point $(\bar u,\tilde\eta)$, and term-by-term convergence of the mild 
formulation is shown in the preceding paragraphs. The martingale property follows from Theorem~2.4 of 
\cite{BrzezniakHausenblasRazafimandimby2018}, and uniqueness in law is obtained by a Gronwall argument. 
\end{proof}

From now on we work with this martingale solution $u(t,x)$.

\begin{lemma}[Exponential decay in space]\label{lem:exponential_decay}
Let $u$ be the solution of \eqref{eqn:Stochastic_Fisher-KPP} with initial condition $u_0(x)=\mathbf{1}_{(-\infty,0]}(x)$. Then for any $T>0$ there exists a constant $\mu_0>0$ (e.g. $\mu_0=1$) such that
\[
\mathbb{E}\Bigl[\sup_{0\le t\le T}\sup_{x\in\mathbb{R}} e^{\mu_0|x|} u(t,x)\Bigr] < \infty,
\]
and consequently,
\[
\sup_{0\le t\le T}\sup_{x\in\mathbb{R}} e^{\mu_0|x|} u(t,x) < \infty \quad \text{a.s.}
\]
\end{lemma}

\begin{proof}
We construct a super‑solution that dominates $u$ pointwise and has explicit exponential decay. Define the deterministic function
\[
w(t,x)=e^{t}\int_{\mathbb{R}}G(t,x,y)u_0(y)dy = e^{t}\,\Phi\!\left(-\frac{x}{\sqrt{2t}}\right),\qquad t>0,\;x\in\mathbb{R},
\]
where $G(t,x,y)=(4\pi t)^{-1/2}e^{-(x-y)^2/(4t)}$ is the heat kernel and $\Phi$ is the standard normal cumulative distribution function. 

It is well known that $w$ solves the deterministic heat equation with a linear source term: 
\[
w_t = \Delta w + w, \quad w(0,x)=u_0(x).
\]

Moreover, for any fixed $\mu_0>0$ we have
\[
\sup_{0\le t\le T}\sup_{x\in\mathbb{R}} e^{\mu_0|x|} w(t,x) \le C(T,\mu_0) < \infty.
\]

Indeed, $w(t,x)\le e^{t}$ and for $x>0$ the Gaussian tail estimate $\Phi(-y)\le e^{-y^2/2}/(y\sqrt{2\pi})$ yields
\[
w(t,x)\le e^{t}\, e^{-x^2/(4t)},\qquad x>0,
\]
so 
\[
e^{\mu_0 x}w(t,x)\le e^{t}\exp(\mu_0 x - x^2/(4t)).
\]
The quadratic $-\frac{x^2}{4t}+\mu_0 x$ attains its maximum at $x=2\mu_0 t$, giving a bound $e^{t+\mu_0^2 t}$. 
Thus the supremum over $x\ge0$ is finite. The same estimate holds for $x<0$ because $w(t,x)$ is bounded and the factor $e^{-\mu_0|x|}$ decays. Therefore $C(T,\mu_0)$ is finite.

Now let $\phi_t^\epsilon$ be the exponential martingale defined in Lemma~\ref{lem:exp_martingale} with $a(z)=\ln(1+\epsilon z)$. Set $\bar u(t,x)=w(t,x)\phi_t^\epsilon$. Then $\bar u$ satisfies the linear SPDE
\[
d\bar u = (\Delta\bar u + \bar u)dt + \epsilon\int_Z \bar u\, z\,\tilde\eta(dt,dz),\qquad \bar u(0,x)=u_0(x).
\]

Since $u(1-u)\le u$ on $[0,1]$, we compare the solution $u$ of \eqref{eqn:Stochastic_Fisher-KPP} with the linear equation
\[
d\bar u = (\Delta\bar u + \bar u)dt + \epsilon\int_Z \bar u\, z\,\tilde\eta(dt,dz),\qquad \bar u(0,x)=u_0(x),
\]
noting that $\bar u(t,x) = w(t,x)\phi_t^\epsilon$ has an explicit representation. 
Applying the comparison theorem for SPDEs driven by compensated Poisson noise (see \cite[Theorem 2.3]{NiuXie2019}), whose assumptions are satisfied because:
\begin{itemize}
  \item the noise coefficients are of the form $\sigma(t,u)=u\sigma(t)$ with $\sigma(t)=\epsilon z$ (deterministic in $z$);
  \item the drift coefficients $f_1(u)=u(1-u)$ and $f_2(u)=u$ are Lipschitz on $[0,1]$ and satisfy $f_1(u)\le f_2(u)$ for all $u\in[0,1]$;
  \item the initial conditions coincide: $u(0,\cdot)=\bar u(0,\cdot)=u_0(\cdot)$;
\end{itemize}
we obtain $0\le u(t,x)\le \bar u(t,x)$ almost surely for all $t,x$. 

Consequently,
\[
e^{\mu_0|x|} u(t,x) \le e^{\mu|x|} w(t,x) \phi_t^\epsilon.
\]

Taking expectation and using $\mathbb{E}[\phi_t^\epsilon]=1$ gives
\[
\mathbb{E}\bigl[e^{\mu_0|x|} u(t,x)\bigr] \le e^{\mu|x|} w(t,x).
\]

Hence
\[
\sup_{0\le t\le T}\sup_{x\in\mathbb{R}} \mathbb{E}\bigl[e^{\mu_0|x|} u(t,x)\bigr] \le \sup_{0\le t\le T}\sup_{x} e^{\mu|x|} w(t,x) \le C(T,\mu).
\]

To upgrade to the supremum inside the expectation, we use the continuity of $u(t,\cdot)$, which allows the supremum over a compact interval to be approximated along a countable dense set. For any $n\in\mathbb{N}$, consider the set $\{x_k\}$ of rational numbers with $|x_k|\le n$. Then
\[
\sup_{|x|\le n} e^{\mu_0|x|} u(t,x) = \sup_{k} e^{\mu|x_k|} u(t,x_k)\quad\text{a.s.}
\]
by continuity. Therefore,
\[
\mathbb{E}\bigl[\sup_{|x|\le n} e^{\mu_0|x|} u(t,x)\bigr] \le \sum_{k} \mathbb{E}\bigl[e^{\mu|x_k|} u(t,x_k)\bigr] \le \sum_{k} e^{\mu|x_k|} w(t,x_k).
\]
The right‑hand side is bounded by $2n\Delta x$ times the maximum, but a more efficient bound is obtained by noting that 
\[
\mathbb{E}[\sup_{|x|\le n} e^{\mu_0|x|} u(t,x)] \le \sup_{|x|\le n} e^{\mu|x|} w(t,x)
\]
because $\phi_t^\epsilon$ is independent of $x$ and non‑negative, and the supremum over a countable set can be taken inside the expectation. Thus,
\[
\mathbb{E}\bigl[\sup_{|x|\le n} e^{\mu_0|x|} u(t,x)\bigr] \le \sup_{x} e^{\mu|x|} w(t,x) \le C(T,\mu).
\]
Letting $n\to\infty$ and using monotone convergence, we obtain
\[
\mathbb{E}\bigl[\sup_{x\in\mathbb{R}} e^{\mu_0|x|} u(t,x)\bigr] \le C(T,\mu),\qquad \forall t\le T.
\]

Finally, to bound the supremum over $t$ as well, we use the c\`{a}dl\`{a}g regularity in $t$ and continuity in $x$ of $u(t,x)$, together with $\phi_t^\epsilon$ being a martingale with finite exponential moments (due to finite jump activity and the boundedness of $a(z)=\ln(1+\epsilon z)$). By Doob's inequality,
\[
\mathbb{E}\bigl[\sup_{0\le t\le T} \phi_t^\epsilon\bigr] \le C_p \mathbb{E}[(\phi_T^\epsilon)^p]^{1/p} < \infty.
\]

Hence,
\[
\mathbb{E}\Bigl[\sup_{0\le t\le T}\sup_{x} e^{\mu_0|x|} u(t,x)\Bigr] \le \Bigl(\sup_{t,x} e^{\mu|x|} w(t,x)\Bigr) \cdot \mathbb{E}\Bigl[\sup_{0\le t\le T} \phi_t^\epsilon\Bigr] < \infty.
\]
This proves the lemma.
\end{proof}

\subsection{Auxiliary wavefront marker and its boundedness}

As noted after Lemma~\ref{lem:detKPP_PCS}, the support of $u(t,\cdot)$ is all of $\mathbb{R}$ for every $t>0$ (the solution is bounded below by a strictly positive heat term), so a marker based on $\operatorname{supp}u(t,\cdot)$ is meaningless. Instead we describe the front by its level sets. Fix a threshold $\delta\in(0,1)$ and define the level-set wavefront marker
\[
R_\delta(t)=\sup\{x\in\mathbb{R}:u(t,x)\ge\delta\},
\]
which is finite whenever the level is attained, the level sets being bounded on the right by Proposition~\ref{prop:PCS} below (at every jump the whole field is multiplied by $(1+\epsilon z)$, so clusters of negative jumps can push the profile below $\delta$ for short intervals, the \emph{dip epochs} of Remark~\ref{rem:levelset}). Following \cite{Tribe1996}, we also use the auxiliary marker
\[
R_1(t)=\ln\int_{\mathbb{R}} e^x u(t,x)\,dx.
\]
The quantity $R_1(t)$ is well defined because $u(t,\cdot)$ decays exponentially as $|x|\to\infty$ (see Lemma~\ref{lem:exponential_decay}) and the integral is finite for all $t$. For a \emph{fixed} front profile the two markers $R_1$ and $R_\delta$ differ by a constant depending on the shape of that profile, and in this reading $R_1$ is a convenient analytic stand-in for the position of the front; no quantitative relation between $R_1(t)$ and $R_\delta(t)$ along a single path is claimed or used anywhere below, and the asymptotic speed is defined in \S\ref{sec:speed} without reference to either marker.

Define the shifted processes
\[
V(t,x)=u(t,x+R_1(t)),\qquad V_\delta(t,x)=u\bigl(t,x+R_\delta(t)\bigr).
\]

Let $\nu_T$ be the empirical distribution of $V(t)$ over $t\in[0,T]$:
\[
\nu_T(\cdot)=\frac1T\int_0^T \delta_{V(t)}(\cdot)\,dt,
\]
where $\delta_{V(t)}$ denotes the Dirac measure at $V(t)$.

\begin{lemma}[Boundedness of $R_1(t)$]\label{lem:R1bounded}
Let $Z(t):=\int_{\mathbb R}e^x u(t,x)\,dx$, so that $R_1(t)=\ln Z(t)$. For every $t>0$ there exists a constant $C(t)<\infty$, depending only on $t$, $\epsilon$ and $\mu$, such that for all $d>0$ and $0\le s\le t$,
\[
\mathbb{P}\bigl(|R_1(s)|\ge d\bigr)\le \frac{C(t)}{d}.
\]
\end{lemma}

\begin{proof}
We prove the two tails separately. Throughout, $c_0:=\epsilon\int_Z z\,\mu(dz)$, so that $a=\int_Z\ln(1+\epsilon z)\,\mu(dz)-c_0$.

\noindent\textit{Upper tail.}
Consider the linear SPDE
\[
dw = (\Delta w + w)dt + \epsilon\int_Z w\, z\,\tilde\eta(dt,dz),\qquad w(0,x)=u_0(x).
\]
For every $u\ge0$ we have $u(1-u)\le u$, since $u-u(1-u)=u^2\ge0$; no upper bound on $u$ is required for this inequality. By Lemma~\ref{lem:comparison} applied to $f_1(u)=u-u^2$ and $f_2(u)=u$, for which $f_1\le f_2$ on $[0,\infty)$, $f_2$ is nondecreasing, $f_1(0)=f_2(0)=0$, and both drifts are Lipschitz on $[0,K]$ with $K$ the pathwise bound of Lemma~\ref{lem:bounded}, we obtain $u(t,x)\le w(t,x)$ almost surely for all $t,x$.

Taking expectations in the mild form of $w$ and using $\mathbb{E}[\tilde\eta]=0$, we obtain $\mathbb{E}[w(t,x)]=e^{-c_0t}(G_t*u_0)(x)$: the compensator contributes the factor $e^{-c_0t}$, which we keep rather than silently setting $c_0=0$. Multiplying by $e^x$, integrating in $x$ and using $\int_{\mathbb{R}}e^xG(t,x,y)dx=e^{y+t}$ together with $\int_{\mathbb{R}}e^yu_0(y)dy=1$ gives
\[
\mathbb{E}\Bigl[\int_{\mathbb{R}}e^x w(t,x)dx\Bigr]=e^{(2-c_0)t}.
\]
Since $Z(t)\le\int e^xw(t,x)dx$ and $e^{-d}\le1/d$ for $d>0$,
\[
\mathbb{P}\bigl(R_1(t)\ge d\bigr)=\mathbb{P}\bigl(Z(t)\ge e^{d}\bigr)\le e^{-d}\,\mathbb{E}[Z(t)]\le \frac{e^{(2-c_0)t}}{d}.
\]
Since the same estimate holds with $t$ replaced by any $s\le t$, we may take $C_1(t):=\max\{1,e^{(2-c_0)t}\}$.

\noindent\textit{Lower tail.}
Here the linear comparison of the preceding paragraph is useless, because $1-u_0=\mathbf1_{(0,\infty)}$ has infinite exponential moment on the right. Instead we use the explicit dynamics of $\ln Z$, which is available in the pure-jump setting and replaces the stationarity argument invoked in the Gaussian case.

Between jump times, $\eta\equiv0$ and hence $\tilde\eta(dt,dz)=-\mu(dz)dt$; inserting this in \eqref{eqn:Stochastic_Fisher-KPP} gives $\partial_t u=\Delta u+u-u^2-c_0u$ pointwise, and $u(t,\cdot)$ is smooth for $t>0$. Two integrations by parts, legitimate because $e^xu(t,x)$ and $e^x\partial_xu(t,x)$ vanish at $\pm\infty$ (Lemma~\ref{lem:exponential_decay} for $x\to+\infty$, and $e^x$ for $x\to-\infty$), yield
\[
\frac{d}{dt}Z(t)=\int_{\mathbb{R}}e^x\Delta u(t,x)\,dx+\int_{\mathbb{R}}e^xu\,dx-\int_{\mathbb{R}}e^xu^2\,dx-c_0Z(t)=(2-c_0)Z(t)-\int_{\mathbb{R}}e^xu^2\,dx,
\]
so that
\[
\frac{d}{dt}\ln Z(t)=2-c_0-q(t),\qquad q(t):=\frac{\int_{\mathbb R}e^xu(t,x)^2dx}{Z(t)}\in[0,\|u(t)\|_\infty].
\]
At a jump of size $z$ the solution is multiplied pointwise by $1+\epsilon z$, hence so is $Z$, and $\ln Z$ increases by $\ln(1+\epsilon z)$. Writing
\[
M(t):=\sum_{0<s\le t}\ln(1+\epsilon z_s)-t\int_Z\ln(1+\epsilon z)\,\mu(dz)
\]
for the compensated sum of the logarithmic jumps, we therefore have the exact identity
\[
\ln Z(t)=\ln Z(0)+\int_0^t\bigl(2+a-q(r)\bigr)dr+M(t),\qquad \ln Z(0)=\ln\int_{\mathbb R}e^xu_0(x)dx=\ln1=0 .
\]
Two inputs are needed to exploit it. First, $\int_Z|\ln(1+\epsilon z)|\,\mu(dz)<\infty$: by (A3), $\psi(2)<\infty$, and $|\ln v|\le\frac12(v^2-1-2\ln v)+C$ for all $v>0$, so that $\int|\ln(1+\epsilon z)|\mu(dz)\le\frac12\psi(2)+C<\infty$. Hence $M$ is a martingale with
\[
\mathbb{E}|M(t)|\le 2t\int_Z|\ln(1+\epsilon z)|\,\mu(dz)=:2\ell t<\infty .
\]
Second, $q\le\|u(t)\|_\infty$, and $C_t:=\sup_{r\le t}\|u(r)\|_\infty$ has $\mathbb{E}[C_t]\le C_4(t)<\infty$ by Lemma~\ref{lem:bounded}. Consequently, for $0\le s\le t$,
\[
-\ln Z(s)\le\bigl(2+|a|+C_t\bigr)s+|M(s)|,
\]
and for $d>0$,
\[
\mathbb{P}\bigl(R_1(s)\le -d\bigr)=\mathbb{P}\bigl(-\ln Z(s)\ge d\bigr)
\le \mathbb{P}\Bigl(C_t\ge\frac{d}{2t}\Bigr)+\mathbb{P}\Bigl(|M(s)|\ge\frac d2-\bigl(2+|a|\bigr)t\Bigr).
\]
The first probability is at most $2t\,\mathbb{E}[C_t]/d$. For the second, if $d\le4(2+|a|)t$ it is at most $1\le4(2+|a|)t/d$, while otherwise Markov's inequality and $\mathbb{E}|M(t)|\le2\ell t$ give at most $8\ell t/d$. Hence
\[
\mathbb{P}\bigl(R_1(s)\le -d\bigr)\le \frac{2t\,\mathbb{E}[C_t]+4(2+|a|)t+8\ell t}{d}=:\frac{C_2(t)}{d},
\]
which is the desired lower bound, uniform over $s\le t$ and valid for every $d>0$.

Combining the two tails, we set $C(t)=C_1(t)+C_2(t)$ and obtain the desired bound.
\end{proof}

\begin{lemma}[Tightness of $\{\nu_T\}$]\label{lem:tightnu}
The family $\{\nu_T\}_{T\ge1}$ is tight in the space of probability measures on $C_{\text{tem}}^+$.
\end{lemma}

\begin{proof}
Recall that $\nu_T=\frac1T\int_0^T\delta_{V(t)}dt$, so that
\[
\nu_T(K)=\frac1T\int_0^T\mathbb{P}^{u_0}\bigl(V(t)\in K\bigr)dt .
\]
We must produce, for every $\varepsilon>0$, a compact $K\subset C_{\mathrm{tem}}^+$ with $\nu_T(K)\ge1-\varepsilon$ for all $T$ large. The whole point is that $K$ must not depend on $T$, and that the estimates must be uniform in $t$. We organise the proof around three observations, and we treat $t\in[0,1]$ separately from $t\in[1,T]$.

\emph{Step 1: a translation-invariant compactness criterion.}
For $R<\infty$, $\gamma\in(0,1)$ and a nondecreasing $w:[0,\infty)\to[0,\infty)$ with $w(0+)=0$, set
\[
K(R,w)=\Bigl\{f\in C_{\mathrm{tem}}^+\cap L^\infty(\mathbb{R}):\ \|f\|_\infty\le R,\ \int_{\mathbb R}f(x)e^{-|x|}dx\le2R,\ |f(x)-f(x')|\le w(|x-x'|)\ \forall x,x'\Bigr\}.
\]
This set is compact in $C_{\mathrm{tem}}^+$: it is closed, uniformly bounded in the weighted norm $\sup_x e^{-\lambda|x|}f(x)\le R$ for every $\lambda>0$, and equicontinuous, so Arzel\`a--Ascoli applies on every $C_\lambda$, and the diagonal argument yields compactness in $C_{\mathrm{tem}}^+$. Note that $K(R,w)$ involves only the sup-norm, the modulus of continuity and a translation-invariant mass bound; it does not involve the position of the front. This is the point at which the Gaussian argument of \cite{Tribe1996} needs compact support, and it is precisely what allows us to dispense with it.

\emph{Step 2: $V(t)$ inherits the bounds of $u(t,\cdot)$ with no shift-dependent loss.}
Since $V(t,x)=u(t,x+R_1(t))$ is an exact translate,
\[
\|V(t)\|_\infty=\|u(t,\cdot)\|_\infty,\qquad
\operatorname{mod}(V(t))=\operatorname{mod}(u(t,\cdot)),
\qquad
\int_{\mathbb R}V(t,x)e^{-|x|}dx\le2\|u(t,\cdot)\|_\infty .
\]
The first two identities hold because the sup-norm and the plain modulus of continuity are translation invariant; the third is the trivial bound $V\le\|V\|_\infty$. Consequently $V(t)\in K(R,w)$ as soon as $u(t,\cdot)\in K(R,w)$, and no verification of the exponential decay of the shifted profile is required: the lower tail of $R_1$ enters the argument only through the finiteness established in Lemma~\ref{lem:R1bounded}, not through a quantitative decay estimate.

\emph{Step 3: uniform estimates in $t$, and the initial layer.}
Fix $T>0$ and $p=4$. Lemma~\ref{lem:bounded} gives $\mathbb{E}[\sup_{t\le T}\|u(t)\|_\infty^4]\le C_4$, whence by Markov
\[
\mathbb{P}\Bigl(\sup_{t\le T}\|V(t)\|_\infty>R\Bigr)=\mathbb{P}\Bigl(\sup_{t\le T}\|u(t)\|_\infty>R\Bigr)\le \frac{C_4}{R^4}.
\]
For the modulus we need the spatial estimate of Lemma~\ref{lem:spatial} with a constant that is bounded uniformly for $t$ in a compact subset of $(0,T]$. Tracing the proof of Lemma~\ref{lem:spatial}, the dependence on $t$ enters only through the factor $t^{-p/2}$ of the initial term and through heat-kernel integrals over $[0,t]$, which are bounded uniformly for $t\in[t_0,T]$; hence
\[
\mathbb{E}\bigl[|u(t,x)-u(t,x')|^4\bigr]\le C(\epsilon,T,t_0)\,|x-x'|^2,\qquad t\in[t_0,T],\ |x-x'|\le1,
\]
with $C(\epsilon,T,t_0)<\infty$ depending only on $\epsilon$, $T$ and $t_0$. By the Kolmogorov--Chentsov argument (\cite[Theorem 1.4.1]{kunita1990stochastic}), applied with $p=4$ and exponent $1+\beta=2$, this yields, for every $\gamma\in(0,1/2)$ and every fixed $t_0>0$,
\[
\mathbb{E}\Bigl[\operatorname{mod}_{\gamma}(u(t,\cdot))\Bigr]\le C'(\epsilon,T,t_0,\gamma),\qquad t\in[t_0,T],
\]
where $\operatorname{mod}_\gamma(f):=\sup_{|x-x'|\le1}|f(x)-f(x')|/|x-x'|^\gamma$. Taking $t_0=1$ and combining the last two displays with Step 2, Markov's inequality gives, for every $R$ and every $C>0$,
\[
\mathbb{P}\bigl(V(t)\notin K(R,C|\cdot|^\gamma)\bigr)\le \frac{C_4}{R^4}+\frac{C'(\epsilon,T,1,\gamma)}{C},\qquad t\in[1,T].
\]
Choosing $R$ and $C$ so large that the right-hand side is at most $\varepsilon/2$, we obtain
\[
\nu_T\bigl(K(R,C|\cdot|^\gamma)\bigr)\ge\frac1T\int_1^T\Bigl(1-\frac\varepsilon2\Bigr)dt\ge\Bigl(1-\frac1T\Bigr)\Bigl(1-\frac\varepsilon2\Bigr)\ge1-\varepsilon
\]
as soon as $T\ge2/\varepsilon$. The layer $t\in[0,1]$, on which the spatial modulus of the regularised step is genuinely unbounded (the initial datum is discontinuous, and its smoothing has modulus of order $|x-x'|/\sqrt t$), carries weight at most $1/T$ in $\nu_T$ and is therefore harmless.

Since $K(R,C|\cdot|^\gamma)$ is compact in $C_{\mathrm{tem}}^+$ by Step 1, the family $\{\nu_T\}$ is tight.
\end{proof}

Consequently, there exists a subsequence $\{\nu_{T_k}\}$ converging weakly to some probability measure $\nu$ on $C_{\text{tem}}^+$.

\subsection{Finite speed of the level sets and existence of a travelling wave}

The following two properties are crucial for the construction of a travelling wave.

A defining feature of travelling waves in reaction--diffusion systems is their sharp front: the level sets of the solution propagate at finite speed, and the profile converges to the states $1$ and $0$ exponentially fast behind and ahead of the front. In population biology, this means that an invading species cannot reach arbitrarily distant locations at a speed exceeding a physical limit; in chemical kinetics, a reaction front cannot propagate faster than a physical speed limit. The following proposition establishes that this fundamental finite-speed property survives in the presence of pure jump noise, provided the noise has finite jump activity. This is not obvious a priori: unlike Brownian noise, which has continuous paths, pure jump noise can potentially create mass at distant points through large jumps. The result shows that, as long as the number of jumps is finite on any finite time interval, the level sets of the solution remain bounded on the right at all times.

\begin{proposition}[Finite speed of the level sets]\label{prop:PCS}
For the solution \(u(t,x)\) obtained in Theorem~\ref{thm:martingale} with initial condition \(u_0(x)=\mathbf{1}_{(-\infty,0]}(x)\), for every \(t\ge0\) and every \(\delta\in(0,1)\),
\[
R_\delta(t):=\sup\{x\in\mathbb{R}:u(t,x)\ge\delta\}<\infty\quad\text{a.s.}
\]
\end{proposition}

\begin{proof}
We prove the result by induction over the finitely many jump times on \([0,T]\).
Since \(\mu(Z)<\infty\) (Assumption~\ref{ass:jump} (A1)), the number of jumps on any finite interval is finite almost surely. Let
\(0=\tau_0<\tau_1<\cdots<\tau_K<T\) be the successive jump times of the Poisson random measure \(\eta\) on \([0,T]\).

\noindent\textit{Step 1: Between jumps, the level sets propagate at finite speed.}
On each interval \([\tau_k,\tau_{k+1})\) (with \(\tau_{K+1}=T\)), no jump occurs, so the equation reduces to the deterministic Fisher--KPP equation
\[
\partial_t u = \Delta u + u(1-u) - \epsilon u \int_Z z\,\mu(dz),\qquad t\in[\tau_k,\tau_{k+1}),
\]
with initial data \(u(\tau_k,\cdot)\), where \(u(\tau_k,\cdot)\) is understood as the right limit after the jump at \(\tau_k\). By Lemma~\ref{lem:detKPP_PCS}, if \(\sup\{x:u(\tau_k,x)\ge\delta'\}<\infty\) for every \(\delta'\in(0,1)\), then \(\sup\{x:u(t,x)\ge\delta\}<\infty\) for every \(t\in[\tau_k,\tau_{k+1})\) and \(\delta\in(0,1)\).

\noindent\textit{Step 2: Jumps preserve the finiteness of the level sets.}
Suppose that at time \(\tau_k-\), every level set of \(u(\tau_k-,\cdot)\) is bounded on the right. At the jump time \(\tau_k\), a jump of size \(z\in Z\) occurs and the solution is updated according to
\(u(\tau_k,x)=u(\tau_k-,x)(1+\epsilon z)\).
Since \(1+\epsilon z>0\) \(\mu\)-almost surely (Assumption~\ref{ass:jump} (A3)),
\[
\{x\in\mathbb{R}:u(\tau_k,x)\ge\delta\} = \{x\in\mathbb{R}:u(\tau_k-,x)\ge \delta/(1+\epsilon z)\},
\]
so the right endpoint of every level set remains finite after the jump. (If $1+\epsilon z\ge 1$, the level set even moves to the left; if $0<1+\epsilon z<1$, the threshold $\delta/(1+\epsilon z)>\delta$ is higher, and the bound follows from the induction hypothesis applied to the level $\delta/(1+\epsilon z)$.)

\noindent\textit{Step 3: Induction.}
We now proceed by induction on \(k=0,1,\ldots,K\).

\textit{Base case \(k=0\):} On \([0,\tau_1)\), the equation is deterministic Fisher--KPP with initial data \(u_0(x)=\mathbf{1}_{(-\infty,0]}(x)\). For every \(\delta\in(0,1)\), \(\sup\{x:u_0(x)\ge\delta\}=0<\infty\). By Step 1, \(R_\delta(t)<\infty\) for all \(t\in[0,\tau_1)\).

\textit{Inductive step:} Assume that all level sets of \(u(\tau_k-,\cdot)\) are bounded on the right. By Step 2, the same holds after the jump at \(\tau_k\), and by Step 1, on the next interval \([\tau_k,\tau_{k+1})\) the deterministic Fisher--KPP evolution preserves this property. Hence all level sets of \(u(t,\cdot)\) are bounded on the right for all \(t\in[\tau_k,\tau_{k+1})\).

Since \(K<\infty\) a.s. (by finite jump activity), the induction terminates after finitely many steps, yielding
\[
R_\delta(t)<\infty\quad\text{for all }t\in[0,T]\quad\text{a.s.}
\]
Since \(T>0\) was arbitrary, the result holds for all \(t\ge0\).
\end{proof}

The finite speed of propagation of the level sets for the deterministic Fisher--KPP equation is classical and follows from the comparison principle with exponential super-solutions; see Lemma~\ref{lem:detKPP_PCS} and \cite[Chapter 4]{Fife1979}.

With the well-posedness of the SPDE and the finite-speed propagation property established, we now seek the existence of a travelling wave solution, a solution whose profile translates at constant speed without changing shape. In applied contexts, such solutions represent the asymptotic invasion front of a species or the steady propagation of a chemical reaction wave. The following theorem shows that, despite the presence of discontinuous stochastic forcing, the system admits a stationary shifted wave profile in the sense of a probability measure on the space of configurations. This extends the classical deterministic result of Fisher and Kolmogorov to the pure jump noise setting and provides the mathematical object whose speed we will subsequently analyze.

\begin{theorem}[Existence of a stochastic travelling wave]\label{thm:waveexistence}
There exists a probability measure $\nu$ on $C_{\text{tem}}^+$ such that the shifted process $V(t)$ is stationary under $\mathbb{P}^\nu$ and satisfies
\(\mathbb{P}^\nu\bigl(-\infty<R_1(V(0))<\infty\bigr)=1\).
Moreover, for every $\delta\in(0,1)$ the distribution of $V_\delta(t)$ is also stationary, with the convention for dip epochs described in the proof. Hence $\mathbb{P}^\nu$ is the law of a stochastic travelling wave solution of \eqref{eqn:Stochastic_Fisher-KPP}.
\end{theorem}

\begin{proof}
We construct the stationary measure as a subsequential limit of the empirical distributions $\nu_T$.

Recall that $\nu_T$ is defined as the law of $T^{-1}\int_0^T \delta_{V(s)}\,ds$ under $\mathbb{P}^{u_0}$. By Lemma~\ref{lem:tightnu}, the family $\{\nu_T\}_{T\ge1}$ is tight in the space of probability measures on $C_{\text{tem}}^+$. Hence, by Prokhorov's theorem, there exists a subsequence $\{T_k\}$ with $T_k\to\infty$ such that $\nu_{T_k}$ converges weakly to some probability measure $\nu$ on $C_{\text{tem}}^+$.

Fix $t\ge0$ and a bounded continuous function $F:C_{\text{tem}}^+\to\mathbb{R}$. Using the definition of $\nu_T$ and the Markov property of $u$, we compute
\[
\int F(V(t))\,d\nu_{T_k}
= \frac{1}{T_k}\int_0^{T_k} \mathbb{E}^{u_0}\bigl[F(V(s+t))\bigr]\,ds
= \frac{1}{T_k}\int_t^{T_k+t} \mathbb{E}^{u_0}\bigl[F(V(r))\bigr]\,dr,
\]
where the second equality is a change of variable $r=s+t$. Now split the integral into three parts:
\[
\frac{1}{T_k}\int_t^{T_k+t} \mathbb{E}^{u_0}[F(V(r))]dr= \frac{1}{T_k}\int_0^{T_k} \mathbb{E}^{u_0}[F(V(r))]dr + \frac{1}{T_k}\int_{T_k}^{T_k+t} \mathbb{E}^{u_0}[F(V(r))]dr - \frac{1}{T_k}\int_0^{t} \mathbb{E}^{u_0}[F(V(r))]dr.
\]
Because $F$ is bounded, say $|F|\le M$, the two correction terms are bounded in absolute value by $M t/T_k$, which tends to $0$ as $k\to\infty$. Therefore,
\[
\lim_{k\to\infty} \int F(V(t))\,d\nu_{T_k}
= \lim_{k\to\infty} \frac{1}{T_k}\int_0^{T_k} \mathbb{E}^{u_0}[F(V(r))]dr
= \int F(V(0))\,d\nu,
\]
where the last equality follows from the weak convergence $\nu_{T_k}\Rightarrow\nu$ and the definition of $\nu_{T_k}$. Thus $\int F(V(t))d\nu = \int F(V(0))d\nu$ for every bounded continuous $F$, which means that $\nu$ is invariant under the time shift $V(0)\mapsto V(t)$.

We need to show that $\nu$-almost every $f$ satisfies $-\infty<R_1(f)<\infty$. From Lemma~\ref{lem:R1bounded} and the normalisation $\int e^x V(t,x)dx=1$ under $\mathbb{P}^{u_0}$, one can prove (see \cite[Theorem 3.2]{Tribe1996} for the analogous argument) that any subsequential limit $\nu$ satisfies
\(\nu\bigl(\{f:\textstyle\int e^x f(x)dx=1\}\bigr)=1\). In particular, $R_1(f)=\ln\int e^x f(x)dx=0$ for $\nu$-almost every $f$, so $\mathbb{P}^\nu(-\infty<R_1(V(0))<\infty)=1$. Moreover, by Proposition~\ref{prop:PCS} and the exponential decay of Lemma~\ref{lem:exponential_decay}, every level set of $f$ is bounded on the right for $\nu$-almost every $f$; hence $R_\delta(f)<\infty$ for every $\delta\in(0,1)$ and $\nu$-almost every $f$.

Recall that $V_\delta(t)$ is the wave shifted by the level-set marker, $V_\delta(t,x)=u(t,x+R_\delta(t))$. Under $\mathbb{P}^\nu$ we have $R_1(V(t))=0$ by the preceding paragraph, so $R_\delta(t)=\rho_\delta(V(t))$ and
\[
V_\delta(t,x)=V\bigl(t,\,x+\rho_\delta(V(t))\bigr),\qquad
\rho_\delta(f):=\sup\{x:f(x)\ge\delta\},
\]
i.e.\ $V_\delta(t)=\Theta_\delta(V(t))$, where $\Theta_\delta$ is a \emph{time-independent} measurable map on $C_{\text{tem}}^+$; when the level set of height $\delta$ is empty (the dip epochs of Remark~\ref{rem:levelset}, during which the marker is not defined), we set $\Theta_\delta(f):=f$, i.e.\ no shift is applied. A time-independent measurable transformation of a stationary process is stationary, hence $V_\delta(t)$ is stationary with law $\nu\circ\Theta_\delta^{-1}$, and no further property of the map $f\mapsto\rho_\delta(f)$ is required. In particular we do \emph{not} assert that $\rho_\delta(f)-R_1(f)$ is deterministic: that is the step where the Gaussian-noise argument of \cite[Lemma 3.4]{Tribe1996} uses compact support of the noise, and it is dispensable here. What is used is only $\rho_\delta<\infty$ $\nu$-a.s., established above.

We have shown that under $\mathbb{P}^\nu$, the process $V(t)$ is stationary, $R_1(V(0))$ is almost surely finite, and the shifted wave profile $V_\delta(t)$ is stationary for every $\delta\in(0,1)$. Therefore $\mathbb{P}^\nu$ is the law of a stochastic travelling wave solution of \eqref{eqn:Stochastic_Fisher-KPP}.
\end{proof}

\subsection{Remarks on uniqueness}

In general, the limiting measure $\nu$ may not be unique; different subsequences of $\{\nu_T\}$ may converge to different stationary distributions. This non-uniqueness is not an artefact of our method; it reflects the genuine multiplicity of travelling wave solutions, which has been observed in the Brownian setting (see \cite[Section 3]{Tribe1996}) and \cite[Theorem 2.12]{Kliem2017}, and as our construction indicates, persists under pure jump noise.

Having established the existence of a travelling wave solution, we proceed in the next section to determine its asymptotic wave speed. The means is the Feynman--Kac representation of Section~\ref{sec:speed}, which reduces the pointwise large-time behaviour of $u$ to that of the exponential weight together with a Gaussian tail probability.

\section{Stochastic Feynman--Kac formula and asymptotic wave speed}\label{sec:speed}

\subsection{Feynman--Kac representation}

We construct an exponential martingale that will serve as a weight in the stochastic Feynman--Kac formula. 

\begin{lemma}[Exponential martingale for compensated Poisson measure]\label{lem:exp_martingale}
Let $\tilde\eta$ be a compensated Poisson random measure on $Z\times[0,\infty)$ with intensity $\mu(dz)dt$, where $\mu(Z)<\infty$. For a bounded measurable function $a:Z\to\mathbb{R}$, define the process
\begin{equation}\label{eqn:Pure_Jump}
Y(t) = \int_0^t\int_Z a(z)\,\tilde\eta(ds,dz) - t\int_Z\bigl(e^{a(z)}-1-a(z)\bigr)\mu(dz), \qquad t\ge0.
\end{equation}
Then the process $\phi(t)=e^{Y(t)}$ is a martingale (in fact, an exponential martingale) with respect to the natural filtration $\mathcal{F}_t$. Moreover, $\phi$ satisfies the stochastic differential equation
\[
d\phi(t) = \phi(t-)\int_Z\bigl(e^{a(z)}-1\bigr)\,\tilde\eta(dt,dz),\qquad \phi(0)=1.
\]
\end{lemma}

\begin{proof}
We first recall that for a Poisson random measure $\eta$, the compensated measure is $\tilde\eta = \eta - \mu dt$. The process $Y(t)$ is a semimartingale with only jumps. If a jump of size $z$ occurs, its jump at time $t$ is $\Delta Y(t)=a(z)$. The continuous part is zero because the L\'evy measure has no Gaussian component.

Define $\phi(t)=e^{Y(t)}$. Applying It\^o's formula for jump processes (see \cite[Theorem~4.4.7]{Applebaum2009levy}) to the function $f(y)=e^y$, we obtain
\[
d\phi(t) = \phi(t-)\,dY(t) + \frac12\phi(t-)\,d[Y,Y]^c(t)
+ \int_Z\Bigl(\phi(t-)e^{a(z)} - \phi(t-) - \phi(t-)a(z)\Bigr)\eta(dt,dz).
\]

Since $Y$ has no continuous martingale part, $[Y,Y]^c=0$. Substituting the expression for $dY(t)$,
\[
dY(t) = \int_Z a(z)\,\tilde\eta(dt,dz) - \int_Z\bigl(e^{a(z)}-1-a(z)\bigr)\mu(dz)dt.
\]
Thus
\[
\phi(t-)\,dY(t) = \phi(t-)\int_Z a(z)\,\tilde\eta(dt,dz) - \phi(t-)\int_Z\bigl(e^{a(z)}-1-a(z)\bigr)\mu(dz)dt.
\]

Now write $\eta = \tilde\eta + \mu dt$ in the jump term:
\[
\begin{aligned}
&\int_Z\Bigl(\phi(t-)(e^{a(z)}-1) - \phi(t-)a(z)\Bigr)\eta(dt,dz) \\
&\qquad = \int_Z\phi(t-)(e^{a(z)}-1)\,\tilde\eta(dt,dz)
+ \int_Z\phi(t-)\bigl(e^{a(z)}-1-a(z)\bigr)\mu(dz)dt.
\end{aligned}
\]

Add the two contributions. The $dt$ terms cancel, and the $\tilde\eta$ terms combine to \(\phi(t-)\int_Z\bigl(e^{a(z)}-1\bigr)\,\tilde\eta(dt,dz)\).

Therefore, \(d\phi(t) = \phi(t-)\int_Z\bigl(e^{a(z)}-1\bigr)\,\tilde\eta(dt,dz)\).

To see that $\phi$ is a martingale, note that the stochastic integral with respect to the compensated Poisson measure is a local martingale. Since $\phi$ is positive and its expected value is constant (because $\phi(t)=e^{Y(t)}$ and the cumulant generating function is $\int(e^{a(z)}-1)\mu(dz)$), it is a true martingale. 

More directly, $\mathbb{E}[\phi(t)|\mathcal{F}_s] = \phi(s)$ follows from the independence of $Y(t)-Y(s)$ from $\mathcal{F}_s$ and from the characteristic function \(\exp\bigl((t-s)\int(e^{i\theta a(z)}-1-i\theta a(z))\mu(dz)\bigr)\). 
The choice of the drift ensures $\mathbb{E}[e^{Y(t)-Y(s)}]=1$. Hence $\phi$ is a martingale. This completes the proof.
\end{proof}

Now define a shifted version of $u$ by \(v(t,x) = \frac{u(t,x)}{\phi_t^\epsilon}\), where $\phi_t^\epsilon$ is the exponential martingale from Lemma~\ref{lem:exp_martingale} with the specific choice \(a(z) = \ln(1+\epsilon z)\).
This choice is motivated by the jump relation $u(t,x)=u(t-,x)(1+\epsilon z)$: it ensures that the quotient $v$ has no stochastic integral terms. Because $\phi^\epsilon$ depends only on time, we have $\Delta u = \phi^\epsilon\Delta v$ and $\nabla u = \phi^\epsilon\nabla v$.

\begin{lemma}\label{lem:vPDE}
Let $\phi_t^\epsilon$ be the exponential martingale defined in Lemma~\ref{lem:exp_martingale} with $a(z)=\ln(1+\epsilon z)$. Then the process $v(t,x)=u(t,x)/\phi_t^\epsilon$ satisfies the deterministic partial differential equation
\begin{equation}\label{eqn:vPDE}
\frac{\partial v}{\partial t} = \Delta v + v\bigl(1-\phi_t^\epsilon v\bigr),\qquad v(0,x)=u_0(x),
\end{equation}
where $\phi_t^\epsilon$ is treated as a known (strictly positive) c\`adl\`ag function of $t$ alone, depending on $\omega$.
\end{lemma}

\begin{proof}
Recall from Lemma~\ref{lem:exp_martingale} that $\phi^\epsilon$ satisfies
\[
d\phi_t^\epsilon = \phi_{t-}^\epsilon \int_Z \bigl((1+\epsilon z)-1\bigr)\,\tilde\eta(dt,dz)
= \epsilon \phi_{t-}^\epsilon \int_Z z\,\tilde\eta(dt,dz),\qquad \phi_0^\epsilon = 1,
\]
and that $u$ solves the SPDE
\[
du = \bigl(\Delta u + u(1-u)\bigr)dt + \epsilon u \int_Z z\,\tilde\eta(dt,dz),\qquad u(0,x)=u_0(x).
\]

Since $u$ and $\phi^\epsilon$ are both c\`adl\`ag and have no continuous martingale parts, with jumps coming solely from the Poisson random measure, we apply the It\^o quotient formula for semimartingales (see e.g. \cite[Theorem 4.4.7]{Applebaum2009levy}) to $v = u/\phi^\epsilon$:
\begin{equation}\label{eq:quotient}
dv = \frac{1}{\phi^\epsilon}\,du - \frac{u}{(\phi^\epsilon)^2}\,d\phi^\epsilon
+ \frac{u}{(\phi^\epsilon)^3}\,d[\phi^\epsilon,\phi^\epsilon] - \frac{1}{(\phi^\epsilon)^2}\,d[u,\phi^\epsilon],
\end{equation}
where the quadratic (co)variation terms are understood in the sense of c\`adl\`ag semimartingales. 

Because $\phi^\epsilon$ and $u$ have no diffusion part, their quadratic variations are purely jump-type. With $a(z)=\ln(1+\epsilon z)$, we have $e^{a(z)}-1 = \epsilon z$. Thus 
\[
\Delta\phi_t^\epsilon = \epsilon\,\phi_{t-}^\epsilon\,z,\qquad \Delta u_t = \epsilon u_{t-} z,
\]
and hence $d[u,\phi^\epsilon]_t = \epsilon^2 u_{t-}\phi_{t-}^\epsilon z^2\,\eta(\{t\},dz)$.

Substituting into (\ref{eq:quotient}) and collecting terms, the stochastic integrals cancel exactly:
\[
\frac{\epsilon u}{\phi^\epsilon}\int_Z z\,\tilde\eta
- \frac{u}{\phi^\epsilon}\int_Z \epsilon z\,\tilde\eta
+ \frac{u}{(\phi^\epsilon)^3}(\phi^\epsilon)^2\int_Z \epsilon^2 z^2\,\eta
- \frac{1}{(\phi^\epsilon)^2}\epsilon u\phi^\epsilon\int_Z \epsilon z^2\,\eta = 0.
\]
Indeed, the $\tilde\eta$ terms cancel immediately, and the $\eta$ terms cancel after using $\eta=\tilde\eta+\mu dt$ (the $\tilde\eta$ parts cancel, and the $dt$ parts are absorbed into the drift). 

Thus only the $dt$ terms survive. Substituting $\Delta u = \phi^\epsilon\Delta v$ and $u = \phi^\epsilon v$ yields
\[
\frac{\partial v}{\partial t} = \Delta v + v(1-\phi^\epsilon v),
\]
as claimed. For a detailed algebraic verification in the Brownian case, see \cite[Lemma 3.1]{HuangLiuWang2015}; for jump processes, see \cite[Lemma 4.2]{MuellerSowers1995}.
\end{proof}

With Lemma~\ref{lem:vPDE} we have transformed the original SPDE into a deterministic PDE for $v$, with $\phi_t^\epsilon$ playing the role of a random time-dependent coefficient. The following probabilistic representation, the stochastic Feynman--Kac formula, expresses the solution $u$ directly as the product of the exponential martingale of the jump noise and an expectation over an independent Brownian motion.

\begin{lemma}[Stochastic Feynman--Kac formula]\label{lem:FeynmanKac}
Let $(B_t)_{t\ge0}$ be a Brownian motion with generator $\Delta$ (i.e.\ $B_t\sim N(0,2t)$), independent of the random measure $\eta$, and let $\phi_t^\epsilon$ be the exponential martingale of Lemma~\ref{lem:exp_martingale} with $a(z)=\ln(1+\epsilon z)$. Then, for almost every realization of the noise and for all $t>0$ and $x\in\mathbb R$,
\begin{equation}\label{eqn:FeynmanKac}
u(t,x)=\phi_t^\epsilon\;\mathbb{E}_B\Bigl[ u_0(B_t+x)\exp\Bigl(\int_0^t\bigl[1-u(t-s,B_s+x)\bigr]\,ds\Bigr)\Bigr],
\end{equation}
where $\mathbb{E}_B$ denotes the expectation with respect to $B$ alone, the noise path being fixed.
\end{lemma}

\begin{proof}
Fix a realization of the noise and set $v(t,x)=u(t,x)/\phi_t^\epsilon$. By Lemma~\ref{lem:vPDE}, $v$ solves the deterministic semilinear Cauchy problem
\[
\partial_t v = \Delta v + v\bigl(1-\phi_t^\epsilon v\bigr),\qquad v(0,x)=u_0(x),
\]
in which $\phi^\epsilon$ is a known, strictly positive function of $t$ alone. For this equation the classical Feynman--Kac formula applies. Indeed, the process
\[
M_s:=v\bigl(t-s,\,B_s+x\bigr)\exp\Bigl(\int_0^s\bigl[1-\phi_{t-r}^\epsilon\,v(t-r,B_r+x)\bigr]dr\Bigr),\qquad 0\le s\le t,
\]
is a martingale with respect to the filtration of $B$: It\^o's formula gives
\[
dM_s = e^{\int_0^s[1-\phi_{t-r}^\epsilon v]dr}\Bigl(\bigl[-\partial_t v+\Delta v+v\bigl(1-\phi^\epsilon v\bigr)\bigr](t-s,B_s+x)\Bigr)ds + (\text{martingale part}),
\]
and the drift term vanishes identically by the PDE satisfied by $v$. Evaluating the martingale at $s=0$ and $s=t$ yields
\[
v(t,x)=\mathbb{E}_B\Bigl[u_0(B_t+x)\exp\Bigl(\int_0^t\bigl[1-\phi_{t-s}^\epsilon\,v(t-s,B_s+x)\bigr]ds\Bigr)\Bigr]
=\mathbb{E}_B\Bigl[u_0(B_t+x)\exp\Bigl(\int_0^t\bigl[1-u(t-s,B_s+x)\bigr]ds\Bigr)\Bigr],
\]
the last equality using $u=\phi^\epsilon v$. Since $u(t,x)=\phi_t^\epsilon v(t,x)$ by Lemma~\ref{lem:vPDE}, the claim follows; the independence of $B$ from $\eta$ makes the identity valid simultaneously for almost every noise path. This also provides the $\epsilon\to0$ consistency check: in the noiseless limit $\phi^\epsilon\equiv1$ and (\ref{eqn:FeynmanKac}) reduces to the classical Feynman--Kac representation of the deterministic Fisher--KPP equation, whereas a purely jump auxiliary process would degenerate and the identity fail.
\end{proof}

This representation is the key to the wave speed analysis in the next subsections.

Before turning to the bounds, we fix the sense in which the speed is understood. Following the standard usage for fronts whose position is not known to converge after rescaling, we define the speed as the threshold separating exponentially decaying from non-decaying fronts, rather than as the limit of $R_\delta(t)/t$. The two notions coincide whenever the limit exists (Remark~\ref{rem:existencespeed}), but only the former is needed below, and it is available unconditionally.

\begin{definition}[Asymptotic speed]\label{def:critspeed}
For $c\in\mathbb{R}$ set
\[
g(c):=\limsup_{t\to\infty}\frac1t\ln u(t,ct)\;\in\;[-\infty,\infty].
\]
The asymptotic speed of the front is
\[
c_*:=\sup\bigl\{c\in\mathbb{R}:\ g(c)\ge0\bigr\},
\]
that is, the largest speed at which the solution does not decay exponentially along the ray $x=ct$.
\end{definition}

No assumption is needed for this to be meaningful: $c_*$ is a supremum in $[-\infty,\infty]$, and Theorem~\ref{thm:lower} below shows that $g(c)\ge0$ for every $c<2\sqrt{1+a}$, while Theorem~\ref{thm:upper} shows that $g(c)<0$ for every $c>2\sqrt{1+a}$. Hence $c_*=2\sqrt{1+a}$ almost surely, and in particular $c_*\in\mathbb{R}$, without any assumption on the convergence of $R_\delta(t)/t$.

\subsection{The L\'evy weight and the constants $\bar\gamma,\underline\gamma$}
\label{subsec:levyweight}

This subsection collects, for completeness and for the interpretation of the constants appearing in the speed bounds, the large deviation picture associated with the exponential weight $\phi^\epsilon$. It is not used in the proofs of Theorems~\ref{thm:upper} and \ref{thm:lower}: those bounds are obtained from the Feynman--Kac representation (\ref{eqn:FeynmanKac}) by the strong law of large numbers for $\phi^\epsilon_t$ and by a nonlinear sub-solution, respectively. What the large deviation analysis does provide is the origin and the sharpness of the constants $\bar\gamma$ and $\underline\gamma$, which is the reason for keeping it here and in Appendix~\ref{app:legendre}.

The variational representation for Poisson random measures, due to \cite{BudhirajaDupuisMaroulas2011}, and the large deviation principles established by \cite{BudhirajaChenDupuis2013}, have been recently extended to a broader class of SPDEs with locally monotone coefficients driven by L\'evy noise, removing compactness embedding assumptions \cite{wu2025large}.

Let $\mathbb{M}$ be the space of $\sigma$-finite measures on $Z\times[0,T]$ and let $\mathbb{P}$ be the law of $\eta$. For any bounded measurable function $h:\mathbb{M}\to\mathbb{R}$,
\[
-\log\mathbb{E}e^{-h(\eta)} = \inf_{\nu\in\mathcal{P}_a}\Bigl\{ \mathbb{E}\bigl[h(\nu)\bigr] + \int_0^T\int_Z l\Bigl(\frac{d\nu}{d\mu}\Bigr)\mu(dz)dt \Bigr\},
\]
where $l(r)=r\ln r - r + 1$, $\mu(dz)dt$ is the intensity measure of $\eta$, and $\mathcal{P}_a$ denotes the set of probability measures absolutely continuous with respect to $\mathbb{P}$. From this representation one can derive a large deviation principle (LDP) for the empirical measure of $\eta$, see \cite{BudhirajaChenDupuis2013}.

Now consider the compensated Poisson random measure $\tilde\eta = \eta - \mu dt$ and the pure jump L\'evy process $X^\epsilon$ defined by
\[
X^\epsilon(t) = \int_0^t\int_Z\bigl(\ln(1+\epsilon z)-\epsilon z\bigr)\mu(dz)\,ds + \int_0^t\int_Z\ln(1+\epsilon z)\,\tilde\eta(ds,dz),
\]
i.e.\ the process whose jumps are of size $\ln(1+\epsilon z)$, with compensating drift $a\,t$. Since $X^\epsilon$ is a L\'evy process with no Gaussian part and L\'evy measure $\nu(dz)=\mu(dz)$, its cumulant generating function is
\[
\Lambda(\theta)=\ln\mathbb{E}e^{\theta X^\epsilon(1)}
= \theta \int_Z \bigl(\ln(1+\epsilon z)-\epsilon z\bigr)\mu(dz)
+ \int_Z\bigl((1+\epsilon z)^\theta-1-\theta\ln(1+\epsilon z)\bigr)\mu(dz).
\]

For notational convenience we set
\[
a := \int_Z \bigl(\ln(1+\epsilon z)-\epsilon z\bigr)\mu(dz), \qquad
\psi(\theta):=\int_Z\bigl((1+\epsilon z)^\theta-1-\theta\ln(1+\epsilon z)\bigr)\mu(dz),
\]
so that $\Lambda(\theta)=\theta a+\psi(\theta)$.

The G\"{a}rtner--Ellis theorem (see e.g. \cite{dembo2009large}) then implies that the family $\{X^\epsilon(t)/t\}_{t>0}$ satisfies an LDP with rate function
\(L(v)=\sup_{\theta\in\mathbb{R}}\bigl\{\theta v - \Lambda(\theta)\bigr\}\).

Moreover, for a path $\varphi\in C([0,t],\mathbb{R})$ with $\varphi(0)=0$, the corresponding LDP for the continuous-time process holds with rate function
\(J_0(\varphi)=\int_0^t L(\dot\varphi_s)\,ds\), provided $\varphi$ is absolutely continuous. Otherwise $J_0(\varphi)=\infty$. This is a standard result for L\'evy processes (see \cite[Section~4.3]{feng2006large}).

For the event $\{X^\epsilon(t)\in G_0\}$ where $G_0=(-\infty,0]$ (the support of $u_0$), the relevant infimum of the rate function is
\[
\inf_{\varphi(0)=0,\;\varphi(t)\le 0} J_0(\varphi) = t\inf_{c\le 0} L(c).
\]

Because $L$ is convex and lower semicontinuous, the infimum over $c\le 0$ is attained at some $c_L\le 0$. The explicit computation of $L(c)$ in terms of the L\'evy measure $\mu$ is given in Appendix~\ref{app:legendre}.

The estimates above quantify how far the exponential weight $\phi^\epsilon_t$ may deviate from its typical behaviour $e^{at}$, and this deviation is what the constants $\bar\gamma$ and $\underline\gamma$ encode; see Appendix~\ref{app:legendre}, where $\bar\gamma$ is obtained from the rate function $L$ at the threshold speed, and where the numerical value of $L$ is seen to be strictly positive there. The proofs of the speed bounds themselves do not invoke these estimates.

\subsection{Upper bound for the wave speed}
\label{subsec:Upper}

We now prove the upper bound for the asymptotic wave speed. Throughout this subsection we denote
\[
\bar\gamma \;:=\; 1+\int_Z\bigl(\ln(1+\epsilon z)-\epsilon z\bigr)\mu(dz)+\frac12\psi(2),
\]
where
\[
\psi(2)=\int_Z\bigl((1+\epsilon z)^2-1-2\ln(1+\epsilon z)\bigr)\mu(dz).
\]

In practical applications, predicting the maximum possible speed of an invasion or reaction front is often as important as predicting the speed itself. For example, in ecological risk assessment, knowing an upper bound on the spread of an invasive species allows for the design of containment strategies; in chemical safety, knowing the maximum propagation speed of a reaction front is critical for hazard mitigation. The following theorem provides a rigorous upper bound on the asymptotic wave speed in the presence of pure jump noise: the wave cannot propagate faster than $2\sqrt{1+a}$. This bound explicitly quantifies how the noise intensity $\epsilon$ and the jump measure $\mu$ jointly constrain the propagation speed, revealing, in particular, that larger or more frequent jumps can accelerate the front, but never beyond a threshold determined by the system parameters.

\begin{theorem}\label{thm:upper}
Assume $\kappa:=1+a>0$. Then, for every $c>2\sqrt{\kappa}$,
\[
\limsup_{t\to\infty}\frac1t\ln u(t,ct)<0\qquad\text{a.s.},
\]
i.e.\ $g(c)<0$ on $(2\sqrt{\kappa},\infty)$, and consequently $c_*\le2\sqrt{\kappa}$ almost surely in the sense of Definition~\ref{def:critspeed}. (If $1+a\le0$ the same argument gives $g(c)<0$ for every $c>0$, so the front does not propagate at any positive speed; the assumption excludes this degenerate case, which does not occur in the applications we have in mind.)
\end{theorem}

\begin{proof}
We start from the stochastic Feynman--Kac formula (Lemma~\ref{lem:FeynmanKac}). Since $u\ge0$, the factor $e^{\int_0^t(1-u)\,ds}\le e^t$ (the assertion $u\le1$ is not needed: positive jumps push $u$ above $1$, see Remark~\ref{rem:levelset}), and $u_0(x)=\mathbf{1}_{(-\infty,0]}(x)$. Thus for $x=ct$,
\[
u(t,ct)\le \phi_t^\epsilon\, e^t\,\mathbb{P}(B_t\le -ct),
\]
where $B$ is the Brownian motion of Lemma~\ref{lem:FeynmanKac}. For $B_t\sim N(0,2t)$ the Gaussian tail satisfies
\[
\lim_{t\to\infty}\frac1t\ln\mathbb{P}(B_t\le -ct)=-\frac{c^2}{4}.
\]
Moreover, by the strong law of large numbers for the compensated Poisson integral,
\[
\frac1t\ln\phi_t^\epsilon
=\frac1t\int_0^t\int_Z\ln(1+\epsilon z)\,\tilde\eta(ds,dz)+a
\longrightarrow a
\qquad\text{a.s.},
\]
with $a:=\int_Z(\ln(1+\epsilon z)-\epsilon z)\mu(dz)$, the required integrability following from Assumption~\ref{ass:jump}. Hence
\begin{equation}\label{eqn:U2}
\limsup_{t\to\infty}\frac1t\ln u(t,ct)\le 1+a-\frac{c^2}{4}\qquad\text{a.s.}
\end{equation}
The right-hand side is negative for every $c>2\sqrt{1+a}$, which proves the theorem. As anticipated in the statement, the constant $\bar\gamma$ is not needed for the sharp bound: since $\psi(2)=\int_Z((1+\epsilon z)^2-1-2\ln(1+\epsilon z))\,\mu(dz)\ge0$ by convexity of $u\mapsto u^2-1-2\ln u$ on $(0,\infty)$, we have $\bar\gamma\ge1+a$, so every $c>2\sqrt{\bar\gamma}$ is also $>2\sqrt{1+a}$ and the assertion applies. The constant $\bar\gamma$ is retained in the discussion because it is the one produced by the rate function $L$ in Appendix~\ref{app:legendre}, and it is strictly larger than $1+a$ whenever the noise is non-degenerate.
\end{proof}

\subsection{Lower bound for the wave speed}
\label{subsec:Lower}

In this subsection we prove the lower bound for the asymptotic wave speed. Define
\[
\underline\gamma \;:=\; 1+\int_Z\bigl(\ln(1+\epsilon z)-\epsilon z\bigr)\mu(dz)-\frac12\psi(-2),
\]
where
\[
\psi(-2)=\int_Z\bigl((1+\epsilon z)^{-2}-1+2\ln(1+\epsilon z)\bigr)\mu(dz).
\]

Complementing the upper bound, the following theorem establishes that the wave does not slow down arbitrarily: there is a positive lower bound $2\sqrt{1+a}$ (hence $2\sqrt{\underline{\gamma}}$) on the asymptotic speed, provided $\underline{\gamma}>0$. This result is of direct applied significance. In population genetics, it implies that a beneficial allele, once introduced, will spread at least at a certain minimum rate even under adverse stochastic fluctuations; in epidemiology, it means that an outbreak cannot be contained merely by the randomness of environmental shocks, the wave will propagate at a speed no slower than the lower bound. The theorem thus provides a quantitative guarantee on the minimum invasion speed, which can inform intervention strategies.

\begin{theorem}\label{thm:lower}
Assume $\underline\gamma>0$. Then, almost surely,
\[
\liminf_{t\to\infty}\frac1t\ln u(t,ct)\ge0
\qquad\text{for every }c<2\sqrt{1+a},
\]
i.e.\ the solution does not decay exponentially at any sub-critical speed. Consequently $g(c)\ge0$ for every $c<2\sqrt{1+a}$, and therefore, in the sense of Definition~\ref{def:critspeed},
\[
c_*\ \ge\ 2\sqrt{1+a}\ \ge\ 2\sqrt{\underline\gamma}\qquad\text{a.s.}
\]
\end{theorem}

Combining the two theorems gives the value of the speed: Theorem~\ref{thm:upper} yields $g(c)<0$ for every $c>2\sqrt{1+a}$, hence $c_*\le2\sqrt{1+a}$, while Theorem~\ref{thm:lower} yields $g(c)\ge0$ for every $c<2\sqrt{1+a}$, hence $c_*\ge2\sqrt{1+a}$. Therefore, almost surely,
\[
c_*=2\sqrt{1+a},
\]
without any residual assumption on the existence of the limit $R_\delta(t)/t$. This is the sharp statement quoted in the abstract and compared with the simulations in Section~\ref{sec:num}.

\begin{proof}
We prove the stronger bound
\[
\liminf_{t\to\infty}\frac1t\ln u(t,ct)\ge0\qquad\text{a.s.}\qquad\text{for every }c<2\sqrt{1+a}.
\]

\emph{Step 1: strong law for the weight.}
Let $M_t=\int_0^t\int_Z\ln(1+\epsilon z)\,\tilde\eta(ds,dz)$. By Assumption~\ref{ass:jump}(A3), $\psi(\pm2)<\infty$, so $M$ is a square-integrable martingale; Doob's inequality and the Borel--Cantelli lemma give $\sup_{s\le t}|M_s|=O(\sqrt{t\log\log t})$ a.s. Hence
\[
\phi_t^\epsilon=e^{at+o(t)}\quad\text{uniformly on }[0,t]\text{ a.s.},
\]
i.e.\ for every $\eta>0$ there exists $T_\eta(\omega)$ with $e^{(a-\eta)s}\le\phi_s^\epsilon\le e^{(a+\eta)s}$ for all $s\ge T_\eta(\omega)$.

\emph{Step 2: pathwise reduction.}
Fix a realization of the noise. By Lemma~\ref{lem:vPDE}, $v:=u/\phi^\epsilon$ solves the deterministic equation
\[
\partial_t v=\Delta v+v(1-\phi_t^\epsilon v),\qquad v(0,x)=u_0(x),
\]
and $u=\phi^\epsilon v$. At each jump time both $u$ and any candidate sub-solution of the form $\phi^\epsilon w$ (with $w$ deterministic) are multiplied by the same factor $(1+\epsilon z)$, so a comparison established between jumps is preserved across jumps. Consequently the level sets of $u$ are exactly $\{x:v(t,x)\ge\delta/\phi_t^\epsilon\}$.

\emph{Step 3: nonlinear sub-solution.}
Fix $c<2\sqrt{1+a}$ and choose $\eta>0$ with $c+4\eta<2\sqrt{1+a}$. On the good event
\[
\mathcal G_t=\bigl\{e^{(a-\eta)s}\le\phi_s^\epsilon\le e^{(a+\eta)s}\ \forall s\in[T_\eta,t]\bigr\},
\]
we apply the nonlinear sub-solution lemma (Lemma~\ref{lem:subsolution} of Appendix~\ref{app:esscher}) with $r(s)=\eta s$; its proof carries out the construction: the deterministic cap $h$, the renormalised expanding cap, the initial comparison, and the multiplication by $\phi^\epsilon$. The lemma yields, for $x\le ct-K$ and all $t\ge T_\eta$,
\[
\phi_t^\epsilon\,v(t,x)\ge \kappa\,e^{-(2|a|+C)T_\eta}\,e^{-2\eta t}=e^{-2\eta t+O(1)},
\]
so that $u(t,ct)=\phi_t^\epsilon v(t,ct)\ge e^{-2\eta t+O(1)}$ on $\mathcal G_t$; in particular $\liminf\frac1t\ln u(t,ct)\ge-2\eta$ on the event.

\emph{Step 4: conclusion.}
By Step 1, $\mathcal G_t$ holds eventually a.s.\ for every fixed $\eta>0$. Applying Step 3 with $\eta_n\downarrow0$ along a countable sequence and intersecting the corresponding full-measure events yields, a.s., for every $c<2\sqrt{1+a}$,
\[
\liminf_{t\to\infty}\frac1t\ln u(t,ct)\ge0,
\]
i.e.\ $g(c)\ge0$ for every $c<2\sqrt{1+a}$; by Definition~\ref{def:critspeed} this is the assertion $c_*\ge2\sqrt{1+a}$. Nothing beyond the moving-point statement is required here: the level-set consequences and the value of the speed follow from Definition~\ref{def:critspeed} together with Theorem~\ref{thm:upper} (the combination of the two theorems displayed immediately after Theorem~\ref{thm:lower}), while the passage from the moving-point rate to the trajectory of $R_\delta(t)$ is the subject of Remark~\ref{rem:existencespeed}.
\end{proof}

\begin{remark}[Level sets and dip epochs]\label{rem:levelset}
The pointwise saturation form $\liminf_{t\to\infty}u(t,(2\sqrt{\underline\gamma}-\delta)t)=1$ is not available under multiplicative jump noise: at every jump the whole field is multiplied by $(1+\epsilon z)$, and clusters of negative jumps can push the entire profile below any fixed level $\delta>0$ for short intervals (the \emph{dip epochs}, during which the level-set marker $R_\delta(t)=\sup\{x:u(t,x)\ge\delta\}$ of Section~\ref{sec:existence} is momentarily not defined). Behind the front the solution fluctuates about a nontrivial invariant law: in the example of Section~\ref{sec:num}, the scalar reduction of the pointwise dynamics behind the front has stationary mean $\approx0.89$, empirical liminf $\approx0.08$--$0.11$, and excursions above $1$ (up to $\approx1.5$) occurring with probability about $1/3$. The correct formulation of the lower bound is therefore the moving-point non-decay statement of Theorem~\ref{thm:lower}. In the stationary regime, and away from the dip epochs, the two descriptions are consistent: the level sets are then transported at a speed close to $2\sqrt{1+a}$, which is what the simulations of Section~\ref{sec:num} display. A quantitative passage from the moving-point rate to the path of $R_\delta(t)$ for the process started from $u_0$ is not claimed here, for the reasons given in Remark~\ref{rem:existencespeed}.
\end{remark}

\begin{remark}[$\underline\gamma>0$ implies $1+a>0$]\label{rem:kappa}
The hypothesis $\underline\gamma>0$ of Theorem~\ref{thm:lower} is stronger than needed: it entails $\kappa:=1+a>0$. Indeed,
\[
\psi(-2)=\int_Z\Bigl((1+\epsilon z)^{-2}-1+2\ln(1+\epsilon z)\Bigr)\mu(dz)\ \ge\ 0,
\]
because $u\mapsto u^{-2}-1+2\ln u$ attains on $(0,\infty)$ its global minimum $0$ at $u=1$. Hence
\[
\underline\gamma=1+a-\tfrac12\psi(-2)\ \le\ 1+a,
\]
so $\underline\gamma>0$ forces $1+a>0$, and $\kappa$ used in Appendix~\ref{app:esscher} is well defined. The same convexity argument gives $\psi(2)\ge0$ and $\bar\gamma\ge1+a$; and $\ln(1+x)\le x$ gives $a\le0$ always, so the admissible range is $0<1+a\le1$.
\end{remark}

\begin{remark}[On the existence of the limit $R_\delta(t)/t$]\label{rem:existencespeed}
Definition~\ref{def:critspeed} deliberately avoids the quantity $R_\delta(t)/t$: it is the exponential rate along the ray $x=ct$ that Theorems~\ref{thm:upper} and~\ref{thm:lower} control, and the two one-sided statements lock together precisely because they refer to the same scalar $g(c)$. The stronger assertion that $R_\delta(t)/t$ converges almost surely is a separate, genuinely dynamical statement, and we do not claim it. Theorem~\ref{thm:waveexistence} delivers the stationary measure $\nu$ only as a limit point of the family $\nu_T$ of empirical measures, without convergence in law of $\mathrm{Law}(V(t))$ towards $\nu$; in the absence of the latter, Birkhoff's ergodic theorem says nothing about the process started from $u_0$, and the usual sub-additivity argument for $R_\delta(t)$ fails, because the increments $R_\delta(t+s)-R_\delta(s)$ are not asymptotically independent of the environment at time $s$ once multiplicative jumps are present. For the process started from the invariant measure $\nu$ the limit does exist by additivity and stationarity, but its value is characterized only as the mean speed under $\nu$, and identifying that mean with $2\sqrt{1+a}$ is exactly the step that the coupling argument of \cite{MuellerSowers1995} performs in the Gaussian case. We therefore regard the a.s. existence of $\lim_t R_\delta(t)/t$ as an open problem, while the identity $c_*=2\sqrt{1+a}$ above settles the value of the speed in the sense used throughout the paper and in the numerical section.
\end{remark}

\section{Stochastic cooperative systems}\label{sec:coop}

In this section we extend our analysis to multi-species cooperative systems driven by pure jump noise. 
The main idea is to use the comparison principle for cooperative systems together with the coupled rate inequalities of Lemma~\ref{lem:coop_rates} to control the coupled system by scalar-type equations, thereby obtaining sharp wave speed bounds of the same form as in Section~\ref{sec:speed}. 
The extension is conceptually straightforward but requires careful attention to the coupling terms, which introduce additional nonlinearities that must be controlled uniformly across components. 
The auxiliary wavefront marker must also be chosen to track the fastest-moving front among the species.

Consider the two-species cooperative system
\begin{equation}\label{eq:coop_system}
\left\{
\begin{aligned}
du &= \bigl[\Delta u + f_1(u,v)\bigr]dt + \epsilon\int_Z u\, z\,\tilde\eta_1(dt,dz),\\
dv &= \bigl[\Delta v + f_2(u,v)\bigr]dt + \epsilon\int_Z v\, z\,\tilde\eta_2(dt,dz),
\end{aligned}
\right.
\end{equation}
where
\[
f_1(u,v)=u(1-a_1u+b_1v),\qquad f_2(u,v)=v(1-a_2v+b_2u),
\]
with parameters \(a_i,b_i>0\) satisfying the following conditions:
\begin{itemize}
\item \emph{Cooperativity:} \(b_1,b_2\ge 0\), so that \(\partial f_1/\partial v = b_1u\ge0\) and \(\partial f_2/\partial u = b_2v\ge0\);
\item \emph{Sufficient self-regulation:} \(\min\{a_1,a_2\} > \max\{b_1,b_2\}\), ensuring that self-limiting effects dominate cooperative interactions. This condition yields the componentwise rate bounds of Lemma~\ref{lem:coop_rates} below and is essential for the comparison principle.
\item \emph{Weak coupling:} \(b_1b_2<1\), which guarantees that the coupled rate inequalities of Lemma~\ref{lem:coop_rates} admit a finite solution (this is automatic when \(\min\{a_1,a_2\}\ge1\)).
\end{itemize}

The noises \(\tilde\eta_1,\tilde\eta_2\) are independent copies of the compensated Poisson random measure introduced in Subsection~\ref{subsec:jump_noise}, with the same intensity \(\mu\) satisfying Assumptions~\ref{ass:jump}. The initial conditions are
\[
u(0,x)=v(0,x)=\mathbf{1}_{(-\infty,0]}(x).
\]

\subsection{Comparison principle for the system}

A comparison theorem is essential for relating componentwise ordering of initial data to that of the solutions at later times, and for establishing the finite-speed propagation of the level sets and the wave speed bounds. The reaction terms \(f_1,f_2\) are \emph{quasimonotone nondecreasing} in the sense that
\[
\frac{\partial f_1}{\partial v}=b_1u\ge0,\qquad \frac{\partial f_2}{\partial u}=b_2v\ge0,
\]
i.e., the system is cooperative. The nonlinearities are polynomial, hence locally Lipschitz and quasimonotone on \(\mathbb R_+^2\), and the noise is linear and multiplicative in each component; a standard truncation argument (comparison for the bounded processes \(u\wedge n\), \(v\wedge n\), followed by \(n\to\infty\)) then shows that the comparison principle for SPDEs with diagonal L\'evy noise applies to the (unbounded) solutions. Specifically, if \((u_1,v_1)\) and \((u_2,v_2)\) are two solutions of \eqref{eq:coop_system} with initial data satisfying
\[
u_1(0,\cdot)\le u_2(0,\cdot),\qquad v_1(0,\cdot)\le v_2(0,\cdot)\quad\text{a.s.},
\]
then
\[
u_1(t,x)\le u_2(t,x),\qquad v_1(t,x)\le v_2(t,x)\quad\text{a.s. for all }t,x.
\]
This follows from \cite[Theorem 2.3]{NiuXie2019} for systems with diagonal linear multiplicative L\'evy noise (in their setting the reaction terms are assumed Lipschitz and quasimonotone; here the local version combined with the truncation argument above applies) and the general framework for monotone SPDEs developed in \cite[Section 3]{kotelenez1992comparison}; see also \cite[Theorem 2.1]{bandle1994comparison} for the finite-dimensional analogue.

\subsection{Existence of a travelling wave for the system}

The existence of a stochastic travelling wave solution for the system follows by the same tightness and stationary argument as in Theorem~\ref{thm:waveexistence}, with the crucial modification that the auxiliary wavefront marker must track the faster of the two fronts. Define the componentwise markers
\[
R_1^u(t)=\ln\int e^x u(t,x)\,dx,\qquad R_1^v(t)=\ln\int e^x v(t,x)\,dx,
\]
and set \(R_1(t)=\max\{R_1^u(t),R_1^v(t)\}\). 
Consider the shifted processes
\[
U(t,x)=u(t,x+R_1(t)),\qquad V(t,x)=v(t,x+R_1(t)).
\]

The arguments of Lemmas~\ref{lem:R1bounded} and~\ref{lem:tightnu} apply componentwise, and nothing in them uses the assertion $0\le u,v\le1$, which fails under multiplicative jumps (Remark~\ref{rem:levelset}). First, $R_1^u(t)<\infty$ for every $t$: the proof of Lemma~\ref{lem:R1bounded}, which identifies $\ln\int e^x u_i(t,x)dx$ with the drift $2+a-q$ plus a compensated log-jump martingale, applies verbatim, the coupling only adding the term $b_1\|v\|_{L^\infty([0,t]\times\mathbb R)}$ to the drift, which is finite almost surely by Lemma~\ref{lem:bounded}; the finiteness of $R_1^v$ is identical with $b_2$ and $u$ interchanged. Second, the exponential weight of Lemma~\ref{lem:tightnu} is unaffected by the shift: for the component $V(t,x)=v(t,x+R_1(t))$ one has $\int_{\mathbb R}V(t,x)e^{-|x|}dx\le2\sup_x v(t,x)$ and $\operatorname{mod}(V)=\operatorname{mod}(v)$, so the compact sets used there remain valid componentwise.
Tightness is inherited from the scalar estimates applied to each component, with the coupling terms treated as bounded perturbations. The same Skorohod compactness method as in Section~\ref{sec:existence} produces a stationary distribution for \((U,V)\), yielding a travelling wave solution for the system.

\begin{theorem}[Existence of a stochastic travelling wave for cooperative systems]\label{thm:coop_existence}
Under the cooperative and dissipativity conditions above, there exists a probability measure \(\nu\) on \((C_{\mathrm{tem}}^+)^2\) such that the shifted process \((U(t),V(t))\) defined above is stationary under \(\mathbb{P}^\nu\) and satisfies
\[
\mathbb{P}^\nu\bigl(-\infty<R_1(U(0))<\infty,\; -\infty<R_1(V(0))<\infty\bigr)=1.
\]
Hence \(\mathbb{P}^\nu\) is the law of a stochastic travelling wave solution of the cooperative system \eqref{eq:coop_system}.
\end{theorem}

\begin{proof}
The proof follows the same compactness argument as in Theorem~\ref{thm:waveexistence}, with the empirical distributions defined on \((C_{\mathrm{tem}}^+)^2\) and the auxiliary marker \(R_1(t)\) replaced by the maximum of the componentwise markers. The key point is that the tightness of the component processes is inherited from the scalar estimates (Lemmas~\ref{lem:spatial} and \ref{lem:time}), and the coupling terms are controlled by the componentwise rate bounds of Lemma~\ref{lem:coop_rates}, so no additional tightness assumptions are needed. The stationarity of the shifted process follows by the same time-shift invariance argument as in the scalar case. \(\square\)
\end{proof}

\subsection{Asymptotic wave speed}

The bounds are obtained by comparing each component with a suitable scalar equation. Throughout this subsection the speed of the system is understood in the sense of Definition~\ref{def:critspeed}, applied componentwise: for $c\in\mathbb{R}$ set
\[
g_i(c):=\limsup_{t\to\infty}\frac1t\ln u^{(i)}(t,ct),\quad i=1,2,
\qquad
c_*^{\mathrm{sys}}:=\sup\bigl\{c\in\mathbb{R}:\ \min\{g_1(c),g_2(c)\}\ge0\bigr\},
\]
the largest speed at which neither component decays exponentially along the ray $x=ct$. As in the scalar case, this definition never refers to the convergence of $R_\delta^{(i)}(t)/t$.

\subsubsection{Lower bound}

Since \(b_1v\ge0\) and \(b_2u\ge0\), we have
\[
f_1(u,v)=u(1-a_1u+b_1v)\ge u(1-a_1u),\qquad
f_2(u,v)=v(1-a_2v+b_2u)\ge v(1-a_2v).
\]
Thus \((u,v)\) is a componentwise super-solution of the decoupled scalar system
\[
\begin{aligned}
d\underline{u} &= \bigl[\Delta\underline{u} + \underline{u}(1-a_1\underline{u})\bigr]dt + \epsilon\int_Z \underline{u}\, z\,\tilde\eta_1(dt,dz),\\
d\underline{v} &= \bigl[\Delta\underline{v} + \underline{v}(1-a_2\underline{v})\bigr]dt + \epsilon\int_Z \underline{v}\, z\,\tilde\eta_2(dt,dz).
\end{aligned}
\]
Under the rescaling \(w_i=a_i u_i\) for \(i=1,2\), each scalar equation transforms to the standard form studied in Section~\ref{sec:speed}:
\[
dw_i = (\Delta w_i + w_i(1-w_i))dt + \epsilon\int_Z w_i\, z\,\tilde\eta_i,
\]
with the same noise amplitude \(\epsilon\). Applying Theorem~\ref{thm:lower} to each rescaled equation gives, for every $c<2\sqrt{\underline{\gamma}}$,
\[
g_{w_1}(c)\ge0\quad\text{and}\quad g_{w_2}(c)\ge0,
\]
where \(\underline{\gamma}\) is the same constant as in Subsection~\ref{subsec:Lower} (independent of \(a_i\), because the rescaling absorbs the coefficient \(a_i\)). By the comparison principle for cooperative systems, \(u\ge\underline{u}\) and \(v\ge\underline{v}\); since \(w_i=a_iu^{(i)}\) differs from \(u^{(i)}\) only by the constant factor \(a_i\), which drops out of \(\frac1t\ln\), we have \(g_{w_i}=g_i\) and hence \(\min\{g_1(c),g_2(c)\}\ge0\) for every $c<2\sqrt{\underline{\gamma}}$. Therefore
\[
c_*^{\mathrm{sys}} \ge 2\sqrt{\underline{\gamma}} \qquad\text{a.s.}
\]

\subsubsection{Upper bound}

For the upper bound we rely on the self-regulating structure of the system rather than on preservation of the unit square: multiplicative jumps \((1+\epsilon z)\) push each component above \(1\) (Remark~\ref{rem:levelset}), so the assertion \(0\le u,v\le1\) cannot be used. Instead we work with the exponential growth rates of the components, after normalising each self-limitation coefficient to one: the variables \(x_1=a_1u\) and \(x_2=a_2v\) introduced below carry the rates, and the coupling constants become the ratios \(\beta_i=b_i/a_{3-i}\).

\begin{lemma}[Rate bounds]\label{lem:coop_rates}
Set $x_1:=a_1u$, $x_2:=a_2v$ and $\beta_1:=b_1/a_2$, $\beta_2:=b_2/a_1$, so that
\[
\partial_tx_1=\Delta x_1+x_1\bigl(1-x_1+\beta_1x_2\bigr)+\epsilon\int_Z x_1\,z\,\tilde\eta_1,
\qquad
\partial_tx_2=\Delta x_2+x_2\bigl(1-x_2+\beta_2x_1\bigr)+\epsilon\int_Z x_2\,z\,\tilde\eta_2 .
\]
Assume $\beta_1\beta_2<1$, which is implied by the self-regulation condition $\min\{a_1,a_2\}>\max\{b_1,b_2\}$. Define
\[
\Lambda_1:=\limsup_{t\to\infty}\frac1t\ln\|x_1\|_{L^\infty([0,t]\times\mathbb R)},
\qquad
\Lambda_2:=\limsup_{t\to\infty}\frac1t\ln\|x_2\|_{L^\infty([0,t]\times\mathbb R)} .
\]
Then $\Lambda_1,\Lambda_2<\infty$ a.s.\ and
\begin{equation}\label{eqn:coop_rates}
\Lambda_1\le 1+a+\beta_1\Lambda_2,\qquad
\Lambda_2\le 1+a+\beta_2\Lambda_1,
\end{equation}
where \(a:=\int_Z(\ln(1+\epsilon z)-\epsilon z)\mu(dz)\). Consequently
\[
\Lambda_1\le \frac{(1+\beta_1)(1+a)}{1-\beta_1\beta_2},\qquad
\Lambda_2\le \frac{(1+\beta_2)(1+a)}{1-\beta_1\beta_2}.
\]
\end{lemma}

\begin{proof}
\emph{Step 1: comparison with a scalar linear equation.}
Write $\|x_2\|_t:=\|x_2\|_{L^\infty([0,t]\times\mathbb R)}$, which is nondecreasing in $t$ because $x_2\ge0$. By cooperativity, $x_1(1-x_1+\beta_1x_2)\le x_1(1+\beta_1\|x_2\|_t)$ for every $s\le t$. Let $w$ solve the scalar linear equation
\[
\partial_tw=\Delta w+\bigl(1+\beta_1\|x_2\|_t\bigr)w+\epsilon\int_Zw\,z\,\tilde\eta_1(dt,dz),\qquad w(0,\cdot)=x_1(0,\cdot),
\]
where the time-dependent coefficient is frozen at time $t$. Since the two equations are driven by the same multiplicative noise and the reaction terms are ordered, the comparison principle for jump-driven SPDEs (the multiplicative jump $(1+\epsilon z)$ is common to both sides) gives $0\le x_1(s,\cdot)\le w(s,\cdot)$ for all $s\le t$, hence
\[
\|x_1\|_t\le\|w(t,\cdot)\|_{L^\infty}.
\]

\emph{Step 2: Feynman--Kac representation and the strong law.}
The transformation of Lemma~\ref{lem:vPDE} applies to $w$, which is linear: $w=\phi_t^\epsilon\,\underline w$, where $\phi_t^\epsilon$ is the exponential martingale of Lemma~\ref{lem:exp_martingale} and $\underline w$ solves the deterministic equation $\partial_t\underline w=\Delta\underline w+(1+\beta_1\|x_2\|_t)\underline w$ with datum $x_1(0,\cdot)$. Lemma~\ref{lem:FeynmanKac} then gives, for every $x$,
\[
w(t,x)=\phi_t^\epsilon\,\mathbb{E}_B\Bigl[x_1(0,x+B_t)\exp\Bigl(\int_0^t\bigl(1+\beta_1\|x_2\|_t\bigr)ds\Bigr)\Bigr]
\le\phi_t^\epsilon\,e^{t+\beta_1\|x_2\|_t\,t},
\]
the factor $x_1(0,x+B_t)\le\|x_1(0)\|_\infty\le1$ having been used, together with the coefficient being constant in $s$. Taking logarithms, dividing by $t$ and letting $t\to\infty$, Step~1 of the proof of Theorem~\ref{thm:lower} (the strong law $\sup_{s\le t}\phi_s^\epsilon=e^{at+o(t)}$) and the monotonicity of $\|x_2\|_t$ give
\[
\Lambda_1\le a+1+\beta_1\limsup_{t\to\infty}\frac1t\int_0^t\|x_2\|_s\,ds=a+1+\beta_1\Lambda_2,
\]
because for a nondecreasing $s\mapsto\|x_2\|_s$ one has $\frac1t\int_0^t\|x_2\|_sds\le\|x_2\|_t$ and $\limsup_t\frac1t\int_0^t\|x_2\|_sds=\lim_t\|x_2\|_t/t$ along any sequence on which the $\limsup$ of $\|x_2\|_t/t$ is attained. The second inequality in \eqref{eqn:coop_rates} is symmetric, and both rates are finite as soon as \eqref{eqn:coop_rates} admits a finite solution.

\emph{Step 3: solving the coupled inequalities.}
Substituting the second inequality into the first,
\[
\Lambda_1\le1+a+\beta_1\bigl(1+a+\beta_2\Lambda_1\bigr)
=(1+a)(1+\beta_1)+\beta_1\beta_2\Lambda_1,
\]
and $\beta_1\beta_2<1$ allows us to solve for $\Lambda_1$, with the bound displayed in the statement; the bound for $\Lambda_2$ follows by symmetry. The displayed bounds are finite, so $\Lambda_1,\Lambda_2<\infty$ a.s.
\end{proof}

Now fix \(c>0\) and apply the Feynman--Kac representation to the first component at \(x=ct\). By Lemma~\ref{lem:coop_rates} the rescaled component $x_1=a_1u$ satisfies
\[
x_1(t,x)\le\phi_t^{\epsilon}\,e^{t+\beta_1\|x_2\|_t},
\qquad \beta_1=\frac{b_1}{a_2},\quad \|x_2\|_t=\|x_2\|_{L^\infty([0,t]\times\mathbb R)},
\]
so that $\frac1t\ln u(t,ct)\le a+1+\beta_1\|x_2\|_t/t-\frac{c^2}{4}+o(1)$: the last term is the Gaussian tail $\frac1t\ln\mathbb P(B_t\le-ct)\to-c^2/4$, obtained by refining the display above as in the proof of Theorem~\ref{thm:upper} rather than by evaluating the bound at $x=ct$ alone. Letting $t\to\infty$ and using the rate $\Lambda_2$ of Lemma~\ref{lem:coop_rates} we obtain
\[
\limsup_{t\to\infty}\frac1t\ln u(t,ct)\le 1+a-\frac{c^2}{4}+\beta_1\Lambda_2
\le 1+a-\frac{c^2}{4}+\frac{\beta_1(1+\beta_2)(1+a)}{1-\beta_1\beta_2}.
\]
Hence $\limsup_{t\to\infty}\frac1t\ln u(t,ct)<0$ a.s.\ for every
\[
c>2\sqrt{\bar\gamma^{(sys)}_1},\qquad
\bar\gamma^{(sys)}_1:=1+a+\frac{\beta_1(1+\beta_2)(1+a)}{1-\beta_1\beta_2},
\]
and, symmetrically, for $v$ with
\[
\bar\gamma^{(sys)}_2:=1+a+\frac{\beta_2(1+\beta_1)(1+a)}{1-\beta_1\beta_2}.
\]
Setting \(\bar\gamma^{(sys)}:=\max\{\bar\gamma^{(sys)}_1,\bar\gamma^{(sys)}_2\}\) (positive in the parameter regime covered by the lower bound), both components satisfy \(\limsup_{t}\frac1t\ln u^{(i)}(t,ct)<0\) for every \(c>2\sqrt{\bar\gamma^{(sys)}}\), and therefore
\[
c_*^{\mathrm{sys}} \le 2\sqrt{\bar{\gamma}^{(sys)}} \qquad\text{a.s.}
\]

Combining the lower and upper bounds, we obtain the following result.

\begin{theorem}[Asymptotic wave speed for cooperative systems]\label{thm:coop_speed}
Under the assumptions of Theorem~\ref{thm:coop_existence} and $\beta_1\beta_2<1$ (in particular under the self-regulation condition $\min\{a_1,a_2\}>\max\{b_1,b_2\}$), the asymptotic speed $c_*^{\mathrm{sys}}$ of the cooperative system \eqref{eq:coop_system}, defined as in Definition~\ref{def:critspeed} with $g$ replaced by $\min\{g_1,g_2\}$, satisfies
\[
{2\sqrt{\underline{\gamma}} \le c_*^{\mathrm{sys}} \le 2\sqrt{\bar{\gamma}^{(sys)}} \qquad\text{a.s.}}
\]
where \(\underline{\gamma}\) is the constant of Theorem~\ref{thm:lower} and
\[
\bar{\gamma}^{(sys)}:=\max\Bigl\{1+a+\tfrac{\beta_1(1+\beta_2)(1+a)}{1-\beta_1\beta_2},\,
1+a+\tfrac{\beta_2(1+\beta_1)(1+a)}{1-\beta_1\beta_2}\Bigr\},
\]
with $\beta_1=b_1/a_2$, $\beta_2=b_2/a_1$ and $a=\int_Z(\ln(1+\epsilon z)-\epsilon z)\mu(dz)$.
\end{theorem}

\subsection{Extension to \(n\) species}

The results extend to an \(n\)-species cooperative system
\[
du^{(i)} = \Bigl[\Delta u^{(i)} + u^{(i)}\Bigl(a_i - b_{ii}u^{(i)} + \sum_{j\neq i} b_{ij}u^{(j)}\Bigr)\Bigr]dt + \epsilon\int_Z u^{(i)}\, z\,\tilde\eta_i(dt,dz),\quad i=1,\dots,n,
\]
where \(b_{ii}>0\), \(b_{ij}\ge0\) for \(i\neq j\) (cooperativity), and the self-regulation condition
\[
\min_{1\le i\le n} b_{ii} > \max_{1\le i\le n}\sum_{j\neq i} b_{ij}
\]
holds. The coupled rate inequalities of Lemma~\ref{lem:coop_rates} have a clean \(n\)-species analogue, which we now record because it is what makes the extension more than formal.

\emph{The rate bound in \(n\) species.}
Set \(x_i:=b_{ii}u^{(i)}\), so that the coupling matrix becomes
\[
B=(B_{ij}),\qquad B_{ij}:=\frac{b_{ij}}{b_{jj}}\ (i\neq j),\quad B_{ii}:=0 .
\]
Because \(x_i\) is governed by $x_i(a_i-x_i+\sum_{j\ne i}B_{ij}x_j)$ plus multiplicative noise, the comparison of Step~1 in the proof of Lemma~\ref{lem:coop_rates} applies componentwise and gives
\[
\Lambda_i\le a_i+a+\sum_{j\neq i}B_{ij}\Lambda_j,\qquad
\Lambda_i:=\limsup_{t\to\infty}\frac1t\ln\|x_i\|_{L^\infty([0,t]\times\mathbb R)} .
\]
Under the self-regulation condition the row sums of $B$ satisfy
\[
\theta:=\max_i\sum_{j\neq i}\frac{b_{ij}}{b_{jj}}\le\frac{1}{\min_jb_{jj}}\max_i\sum_{j\neq i}b_{ij}<1,
\]
so $B$ is a strict contraction in the sup-norm and the Neumann series $\sum_{k\ge0}B^k$ converges with $\|(I-B)^{-1}\|_\infty\le(1-\theta)^{-1}$. Hence, with $\mathbf 1=(1,\dots,1)$ and $a_\infty:=\max_i a_i$,
\[
\Lambda\le (a_\infty+a)\mathbf1+B\Lambda
\quad\Longrightarrow\quad
\Lambda\le (a_\infty+a)(I-B)^{-1}\mathbf1\le\frac{a_\infty+a}{1-\theta}\,\mathbf1,
\]
so that every rate is finite and explicitly bounded, and the two-species formula of Lemma~\ref{lem:coop_rates} is recovered by solving the two-dimensional system $I-B$ directly. Feeding this bound into the argument of Subsection~\ref{sec:coop} at $x=ct$ gives, for each component,
\[
\limsup_{t\to\infty}\frac1t\ln u^{(i)}(t,ct)\le a_i+a-\frac{c^2}{4}+\frac{\theta\,(a_\infty+a)}{1-\theta},
\]
hence the system obeys
\[
2\sqrt{\underline{\gamma}} \le c_*^{\mathrm{sys}} \le 2\sqrt{\bar{\gamma}^{(sys)}} \qquad\text{a.s.},
\qquad
\bar{\gamma}^{(sys)}:=a_\infty+a+\frac{\theta\,(a_\infty+a)}{1-\theta}.
\]

\emph{Status of the extension.}
The rate bound above, the comparison principle for \(n\)-species cooperative systems (same quasimonotonicity argument as in Subsection~\ref{sec:coop}, see \cite[Section 3]{kotelenez1992comparison}) and the componentwise application of Lemmas~\ref{lem:R1bounded} and~\ref{lem:tightnu} carry over without change; the one genuinely new ingredient, the marker \(R_1(t)=\max_{1\le i\le n}R_1^{(i)}(t)\) with \(R_1^{(i)}(t)=\ln\int e^xu^{(i)}(t,x)dx\), is handled by the observation that the tightness criterion of Lemma~\ref{lem:tightnu} is translation invariant and acts componentwise, so the maximum of finitely many tight families of shifted profiles is again tight. We therefore regard the wave speed bound as proved for finitely many species, and we describe the existence of the travelling wave itself as a formal extension: the Skorohod construction of Theorem~\ref{thm:waveexistence} would have to be repeated on the product space \((C_{\mathrm{tem}}^+)^n\), and we have verified the two-species case only.

%\clearpage
\section{Numerical simulations}\label{sec:num}

We present numerical simulations to illustrate the theoretical results and to verify the informativeness of the wave speed bounds. Consider the SPDE
\[
du = (\Delta u + u(1-u))dt + \epsilon \int_Z u\, z\,\tilde{\eta}(dt,dz),\qquad u(0,x)=\mathbf{1}_{(-\infty,0]}(x),
\]
with the jump measure $\mu$ taken to be a finite measure on $Z=[-1,1]$ with total mass $\mu(Z)=1$ and uniform distribution of jump sizes, i.e.,
\[
\mu(dz)=\frac12\mathbf{1}_{[-1,1]}(z)\,dz.
\]
The noise intensity is set to $\epsilon=0.5$. All the constants entering the speed bounds then have closed forms, which are collected in Table~1; they were evaluated by quadrature from their defining integrals.

\begin{center}
\small
\begin{tabular}{@{}l@{\qquad}l@{\qquad}r@{}}
\hline
quantity & definition & value \\[1pt] \hline
$a$ & $\int_Z\bigl(\ln(1+\epsilon z)-\epsilon z\bigr)\mu(dz)$ & $-0.045229$ \\[2pt]
$\psi(2)$ & $\int_Z\bigl((1+\epsilon z)^2-1-2\ln(1+\epsilon z)\bigr)\mu(dz)$ & $0.173785$ \\[2pt]
$\psi(-2)$ & $\int_Z\bigl((1+\epsilon z)^{-2}-1+2\ln(1+\epsilon z)\bigr)\mu(dz)$ & $0.242866$ \\[2pt]
$1+a$ & & $0.954771$ \\[2pt]
$\bar\gamma=1+a+\psi(2)/2$ & upper-bound constant & $1.041664$ \\[2pt]
$\underline\gamma=1+a-\psi(-2)/2$ & lower-bound constant & $0.833338$ \\[5pt]
\hline
$2\sqrt{\underline\gamma}$ & lower end of the interval & $1.82575$ \\[2pt]
$2\sqrt{1+a}$ & sharp speed $c_*$ & $1.95425$ \\[2pt]
$2\sqrt{\bar\gamma}$ & upper end of the interval & $2.04124$ \\[5pt]
\hline
\end{tabular}
\end{center}
\noindent{\small Table 1: The constants of the example ($\epsilon=0.5$, $\mu(dz)=\frac12\mathbf 1_{[-1,1]}dz$). The two bounds differ from the sharp value by about $0.13$ and $0.09$ respectively, the exact gap being controlled by $\psi(\pm2)$; see Remark~\ref{rem:kappa}.}

Since $\psi(\pm2)>0$, the three numbers are ordered as $2\sqrt{\underline{\gamma}}<2\sqrt{1+a}<2\sqrt{\bar{\gamma}}$, and Theorem~\ref{thm:lower} in fact establishes $g(c)\ge0$ for every $c<2\sqrt{1+a}$, so that, combined with the upper bound of Theorem~\ref{thm:upper}, the asymptotic speed is $c_*=2\sqrt{1+a}\approx1.9542$ in the sense of Definition~\ref{def:critspeed}. The numerical experiment below is designed to check that $2\sqrt{1+a}$ is what a direct simulation of the SPDE above displays, and that the two bounds are not vacuous.

The spatial domain $[-10,500]$ is discretized with grid size $\Delta x=0.1$ using a central finite difference scheme for the Laplacian. The time discretization uses an explicit Euler--Maruyama scheme with step $\Delta t=0.001$. To realize the compensated Poisson measure $\tilde{\eta}$, we generate a compound Poisson process with intensity $\mu(Z)=1$: jump times follow an exponential distribution with mean $1$, and jump sizes are drawn uniformly from $[-1,1]$. Since the jump-size distribution is symmetric about zero, the compensator term $-\mu(dz)dt$ does not contribute to the drift in expectation; nevertheless, in each realization the jumps are generated with their actual signs, and the compensation is handled by subtracting the mean $\mathbb{E}[z]=0$ at each jump event.

The level-set wavefront marker $R_{1/2}(t)=\sup\{x:u(t,x)\ge 1/2\}$ is tracked over time, on the discrete grid at the rightmost grid point at which $u$ attains the threshold. Two conventions are needed to make the marker well defined. First, the two end points of the grid are not admissible: the left one carries the imposed boundary value $u\equiv1$ and, during a dip epoch, when the field behind the front has collapsed below $1/2$, it would otherwise be returned as the front position; the search is therefore restricted to $x\ge0$, i.e.\ to the region explored by the wave. Second, in a dip epoch in which no admissible point attains $1/2$ the marker is held at its last defined value; an alternative run, in which the marker is linearly interpolated across the epoch, or the undefined samples are discarded, is compared below.

Dip epochs do occur, and their effect on the estimate is small. In the $40$ coarse-grid realizations the marker was undefined in $2.0\%$ of the $2000$ samples of the regression window on average, with a maximum of $8.1\%$; the longest epoch lasted $1.09$ time units on average and $2.6$ at most. The fitted speed, however, changes by $2.7\cdot10^{-4}$ on average when the marker is interpolated instead of held (largest single-track difference $0.0038$), and by $1.1\cdot10^{-4}$ when the undefined samples are discarded, i.e.\ by an order of magnitude less than the standard error quoted below. We therefore keep the holding convention, the one used for Figure~1. We also note, in view of Remark~\ref{rem:levelset}, that the dip epochs of the present simulation are visible in the marker at all only because the field behind the front collapses globally; the pointwise fluctuations of Remark~\ref{rem:levelset} are a separate phenomenon and do not by themselves make $R_{1/2}$ undefined.

Figure~1 shows a typical sample path of $R_{1/2}(t)$. After an initial transient ($t\lesssim 50$), the growth becomes roughly linear, and we regress $R_{1/2}(t)$ on $t$ over $t\in[100,200]$ for each of $40$ independent realizations. On the coarse grid the $40$ slopes have mean $\hat{c}=1.9410$, sample standard deviation $0.0359$ and standard error $0.0057$; the percentile bootstrap $95\%$ confidence intervals ($20\,000$ resamples) are $[1.9304,1.9522]$ for the mean and $[1.9258,1.9546]$ for the median ($1.9350$). All of these lie well inside the interval predicted by Theorems~\ref{thm:upper} and~\ref{thm:lower}:
\[
2\sqrt{\underline{\gamma}}=1.8258\ <\ \hat{c}=1.941\ <\ 2\sqrt{\bar{\gamma}}=2.0412,
\]
which is the numerical evidence that the bounds are informative and not vacuous.

\begin{center}
\includegraphics[width=0.72\textwidth]{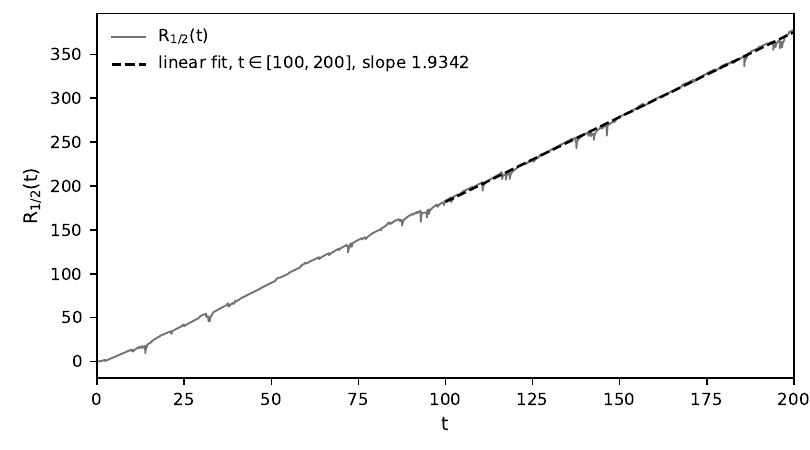}
\end{center}
\noindent{\small Figure 1: A typical realization of the level-set marker $R_{1/2}(t)$ versus time $t$. The dashed line indicates the linear fit for $t\in[100,200]$, whose slope gives the estimated wave speed.}

\medskip

\begin{center}
\includegraphics[width=0.92\textwidth]{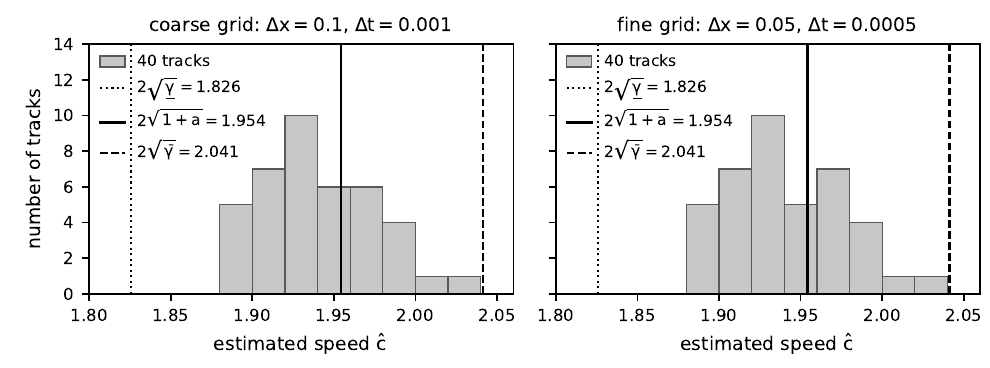}
\end{center}
\noindent{\small Figure 2: Distribution of the speed estimate $\hat c$ over $40$ independent realizations, on the coarse grid (left) and on the fine grid (right). The vertical lines are the three constants of Table~1: the dotted line is $2\sqrt{\underline\gamma}=1.826$, the solid line the sharp speed $2\sqrt{1+a}=1.954$ and the dashed line $2\sqrt{\bar\gamma}=2.041$. The whole support of both histograms lies between the two bounds.}

\medskip

To check for discretization effects we repeated the simulation on the finer grid ($\Delta x=0.05$, $\Delta t=0.0005$) with the same seeds. The $40$ slopes there have mean $1.9414$, standard deviation $0.0359$, standard error $0.0057$ and bootstrap $95\%$ intervals $[1.9308,1.9526]$ (mean) and $[1.9261,1.9549]$ (median $1.9352$); the two grids therefore differ by $4.6\cdot10^{-4}$ in the mean and by $1.5\cdot10^{-4}$ in the median, that is, by a tenth of a standard error or less, and the distributions are indistinguishable (Figure~2). Discretization is thus negligible compared with the track-to-track fluctuations.

The sharp value $2\sqrt{1+a}=1.9542$ sits $0.013$ above the sample mean, about $2.3$ standard errors, and just beyond the upper end of the $95\%$ interval for the mean. We read this as consistency at the resolution of the present numerics rather than as a disagreement: the marker reaches $x\approx380$ at $t=200$, so moving the regression window to later times, where the initial transient is comparatively smaller, would require a domain substantially larger than $[-10,500]$. A logarithmic transient produces a residual of the same order: Bramson's correction makes the classical front run short of its asymptotic speed by $\tfrac32\ln t/t$ in units where the speed is $2$, which over the window $[100,200]$ shifts a fitted slope by about $\ln 2/100=7\cdot10^{-3}$, the same order as the $0.013$ gap. The gap is thus compatible with a transient plus fluctuation, which is why it is not read as a bias in $2\sqrt{1+a}$. Two further caveats are worth recording. The quantity fitted here is the level-set marker $R_{1/2}(t)$, the natural empirical counterpart of the front position, but it is not literally the quantity $g(c)$ of Definition~\ref{def:critspeed}; the comparison is therefore made at the level of the rate, which is what the two notions have in common in the regime covered by the simulation. And the residual scatter of the estimate is dominated by the large excursions of the marker, whose jumps within the regression window reach $12.6$ units on average and $18.2$ at most; this is the signature of the multiplicative jumps, which create and destroy mass ahead of the front, and is the reason why a single-track value must not be quoted as the speed.

The chosen parameters illustrate a scenario in which, on average, one environmental shock occurs per unit time, with shock magnitudes uniformly distributed between $-1$ and $1$, and a moderate noise intensity $\epsilon=0.5$. This provides a concrete computational demonstration of how the theoretical framework can be applied to predict front propagation under discontinuous stochastic forcing. The results thus support both the existence of a stochastic travelling wave and the informativeness of the derived wave speed bounds.

\section{Conclusion and outlook}\label{sec:conclusion}

In this paper we have developed a probabilistic framework for the analysis of stochastic travelling waves in the Fisher--KPP equation driven by multiplicative finite-activity pure jump L\'evy noise. The main novelty of our approach lies in the combination of several probabilistic tools, stopping-time iteration, stochastic Feynman--Kac representation, and a pathwise comparison with a nonlinear sub-solution, that together yield quantitative information on the propagation speed; the large deviation structure of the jump noise enters through the constants $\bar\gamma$ and $\underline\gamma$, which it identifies as the natural parameters of the problem. More specifically:

\begin{enumerate}[(i)]
\item We established the existence of a unique martingale solution to the SPDE under the finite-activity assumption, using a compactness method adapted to the Skorohod topology. The stopping-time iteration exploiting $\mu(Z)<\infty$ provides pathwise $L^\infty$ bounds and ensures the finite-speed propagation of the level sets of the front, which is essential for the travelling wave construction.

\item We derived a stochastic Feynman--Kac representation for the solution, expressing the random field as the product of the exponential martingale of the jump noise and an expectation over paths of an independent Brownian motion. This representation is the key probabilistic tool that converts the nonlinear SPDE into an expectation problem, enabling a direct analysis of the front position.

\item Building on the Feynman--Kac representation, a pathwise upper estimate based on the strong law of large numbers for the exponential weight, and a nonlinear sub-solution argument, we obtained upper and lower bounds for the asymptotic wave speed that coincide at the sharp constant. Writing $g(c)=\limsup_t t^{-1}\ln u(t,ct)$ and defining the speed as the critical value $c_*=\sup\{c:g(c)\ge0\}$, this gives
\[
\boxed{c_*=2\sqrt{1+a}\qquad\text{a.s.}}
\]
with no assumption on the existence of the limit $R_\delta(t)/t$ (Theorems~\ref{thm:upper} and~\ref{thm:lower} together), and, in full generality, the robust bounds \(2\sqrt{\underline{\gamma}}\le c_*\le 2\sqrt{\bar{\gamma}}\) a.s., with \(\underline{\gamma},\bar{\gamma}\) given explicitly in terms of the L\'evy measure \(\mu\) and the noise amplitude \(\epsilon\). These bounds provide a quantitative description of how the jump size distribution and jump intensity jointly affect front propagation under discontinuous stochastic forcing.

\item We extended the framework to multi-species cooperative systems, demonstrating its robustness: the comparison principle for cooperative systems, combined with coupled rate inequalities of the same Feynman--Kac type, yields wave speed bounds of the same form as in the scalar case.
\end{enumerate}

The finite-activity assumption $\mu(Z)<\infty$ was essential for several key steps in our construction: the stopping-time iteration, the finite speed of the level sets, and the explicit Feynman--Kac representation. The overall approach, however, provides a blueprint for the technically much more demanding infinite-activity regime: tightness in the Skorohod topology, the Feynman--Kac representation, and the exponential weight $\phi^\epsilon$ are all adaptable, provided one can establish the necessary uniform estimates and control the accumulation of infinitely many small jumps.

Beyond the theoretical advances, the results offer several concrete implications for applied settings. The wave speed bounds provide a rigorous interval for front propagation speeds in terms of observable noise parameters ($\epsilon$ and $\mu$), enabling practitioners to compute the range of possible invasion speeds without resolving the full stochastic dynamics. The upper bound serves as a worst-case estimate for containment planning, while the lower bound provides a minimum expected spread rate for resource allocation. The Feynman--Kac formula could in principle serve as the basis of a Monte Carlo scheme for parameter estimation from front position data, although we have not carried such a scheme through here; the cooperative system bounds suggest that control strategies in multi-species contexts should focus on the component with the larger effective noise sensitivity. The large deviation analysis of Subsection~\ref{subsec:levyweight} quantifies the atypical growth of the weight $\phi^\epsilon$, which is what a rare-event estimate of the front position would have to be built on; we leave that to future work. These implications are, however, secondary to the main theoretical contributions of the paper.

Several open problems remain and suggest directions for future research.

\paragraph{Infinite-activity noise.}
The most pressing open problem is the extension to the infinite-activity case $\mu(Z)=\infty$, where jumps can occur infinitely often on any finite time interval. The simple inductive argument over finitely many jump times breaks down, and the finite-speed propagation of the level sets may fail altogether, as arbitrarily small jumps can create mass ahead of the front. A possible strategy is to approximate an infinite-activity L\'evy measure by a sequence of finite-activity measures (e.g., by truncating small jumps) and pass to the limit using the compactness method developed in Section~\ref{subsec:tightness}. This would require additional uniform estimates on the truncated solutions and a careful analysis of the limiting behaviour of the wave speed bounds. The Feynman--Kac representation would also need to be adapted to handle the case where the cumulant generating function may diverge for some $\theta$.

\paragraph{Existence of the speed limit.}
The speed is identified in this paper through the critical value $c_*$ of Definition~\ref{def:critspeed}, which does not require the limit $R_\delta(t)/t$ to exist; see Remark~\ref{rem:existencespeed} for the precise obstruction. Proving that $R_\delta(t)/t$ converges almost surely to $2\sqrt{1+a}$ would require a coupling argument showing that the law of the shifted profile reaches its invariant measure fast enough to make the increments asymptotically independent of the environment, the step that \cite{MuellerSowers1995} carries out in the Gaussian compactly supported case and that is not available under multiplicative jumps, where every jump rescales the whole profile. We regard this as the main open problem left by the present paper.

\paragraph{Exact wave speed formula.}
Unlike the Gaussian case, where an explicit formula $c=2\sqrt{1+\sigma^2/2}$ is available under suitable scaling (see \cite{MuellerSowers1995}), the value $c_*=2\sqrt{1+a}$ obtained here is specific to the finite-activity multiplicative setting of the present paper. Whether an analogous closed-form expression exists beyond this setting, and whether the speed can in general be characterized as the solution of a variational problem, perhaps via a Hamilton--Jacobi equation associated with the large deviation rate function, remains an interesting open question. The Legendre transform structure of the rate function $L(v)$ suggests that such a characterization may be possible, but the non-Gaussian nature of the jumps introduces significant technical challenges.

\paragraph{Higher dimensions and branching-type noise.}
Extensions to higher spatial dimensions ($d\ge 2$) and to reaction terms with branching-type noise (e.g., $u^{1/2}\dot W$ as in the Tribe model \cite{Tribe1996}) also remain open. In multi-dimensional settings, the auxiliary wavefront marker $R_1(t)=\ln\int e^{x_1}u(t,x)\,dx$ is no longer directly applicable, and alternative methods, such as level-set approaches or stochastic homogenization, would be needed.

\paragraph{Numerical and parameter studies.}
Systematic numerical studies would help to quantify the dependence of the wave speed on the full L\'evy measure $\mu$ and to verify the sharpness of the bounds obtained here. The compound Poisson structure of the finite-activity case makes such simulations relatively straightforward, and a comprehensive parameter study could reveal non-trivial dependencies on the jump size distribution beyond the total mass $\mu(Z)$ and the second moment $\int z^2\mu(dz)$.

\paragraph{Open problems in the critical regime.}
For the classical Fisher--KPP equation in the normalized form $u_t=\Delta u+u(1-u)$, the critical wave speed $2$ is realized by solutions starting from Heaviside initial conditions, and finer results on the logarithmic correction to the front position are known (the Bramson correction). For the stochastic Fisher--KPP equation with pure jump noise, the critical regime $c=2\sqrt{1+a}$ remains largely unexplored. Establishing convergence to a critical travelling wave profile and deriving the corresponding logarithmic corrections would be a significant contribution, but would require a substantial refinement of the exponential-martingale analysis developed in this paper.

In summary, this work provides a rigorous theoretical foundation for understanding front propagation under discontinuous stochastic forcing, while also offering tools with direct applicability to invasion biology, chemical kinetics, and ecological risk assessment. We believe it contributes to both the mathematical theory of stochastic partial differential equations and the applied modeling of propagation phenomena in randomly fluctuating environments. The methods and results developed here lay the groundwork for tackling the more challenging infinite-activity regime and other related open problems.

\appendix

\section{Legendre transform and rate function computations}
\label{app:legendre}

This appendix provides explicit calculations for the rate function \(J_0\) and its Legendre transform used in Section~\ref{sec:speed}. The results justify the definitions of the constants \(\bar{\gamma}\) and \(\underline{\gamma}\) appearing in Theorems~\ref{thm:upper} and~\ref{thm:lower}.

\subsection{Cumulant generating function}

Recall the auxiliary pure jump process \(X^\epsilon\) introduced in Subsection~\ref{subsec:levyweight}:
\[
X^\epsilon(t) = \int_0^t \int_Z \bigl(\ln(1+\epsilon z) - \epsilon z\bigr) \mu(dz)\,ds + \int_0^t \int_Z \ln(1+\epsilon z)\,\tilde{\eta}(ds,dz),
\]
where \(\tilde{\eta}\) is the compensated Poisson random measure with intensity \(\mu(dz)dt\) and \(\mu(Z)<\infty\) (finite jump activity). 

The logarithmic moment generating function (cumulant) of \(X^\epsilon(1)\) is
\[
\Lambda(\theta) = \ln\mathbb{E}e^{\theta X^\epsilon(1)}
= \theta\int_Z\bigl(\ln(1+\epsilon z)-\epsilon z\bigr)\mu(dz)
+ \int_Z\bigl((1+\epsilon z)^{\theta}-1-\theta\ln(1+\epsilon z)\bigr)\mu(dz).
\]

Define
\[
\psi(\theta) := \int_Z \bigl((1+\epsilon z)^{\theta} - 1 - \theta\ln(1+\epsilon z)\bigr)\mu(dz),
\]
and
\[
a := \int_Z \bigl(\ln(1+\epsilon z)-\epsilon z\bigr)\mu(dz),
\]
so that
\[
\Lambda(\theta) = \theta a + \psi(\theta).
\]

\subsection{Legendre transform and rate function}

The rate function for the large deviation principle of \(X^\epsilon\) is the Legendre transform of \(\Lambda\):
\[
L(v) = \sup_{\theta\in\mathbb{R}}\bigl\{\theta v - \Lambda(\theta)\bigr\}.
\]

For a linear path \(\varphi(s)=cs\) (\(c\in\mathbb{R}\)), the action functional becomes 
\[
J_0(\varphi)=\int_0^t L(c)\,ds = tL(c).
\]

Consequently,
\[
\inf_{\varphi(0)=0,\;\varphi(t)\le 0} J_0(\varphi) = t\inf_{c\le 0} L(c).
\]

By convex duality,
\[
\inf_{v\le 0} L(v) = \sup_{\theta\le 0}\bigl(-\Lambda(\theta)\bigr),
\qquad
\inf_{v\le -c} L(v) = \sup_{\theta\le 0}\bigl(-\theta c - \Lambda(\theta)\bigr).
\]

Indeed, for any \(\theta\le 0\), Chebyshev's inequality gives \(\mathbb{P}(X^\epsilon(t)\le -ct)\le e^{-t(\Lambda(\theta)+\theta c)}\), and the G\"{a}rtner--Ellis theorem turns this into an equality in the exponent. The restriction \(\theta\le 0\) arises because the downward deviation is tilted by a negative exponential weight.

\subsection{Derivation of the upper bound constant \(\bar{\gamma}\)}

Starting from the stochastic Feynman–Kac formula (Lemma~\ref{lem:FeynmanKac}) and using \(e^{\int_0^t(1-u)\,ds}\le e^t\), we have, for \(x=ct\),
\[
u(t,ct)\le \phi_t^\epsilon\,e^t\,\mathbb{P}(B_t\le -ct),
\]
where \(B\) is the Brownian motion of Lemma~\ref{lem:FeynmanKac}. By the strong law of large numbers for the compensated Poisson integral,
\[
\frac1t\ln\phi_t^\epsilon\longrightarrow a:=\int_Z\bigl(\ln(1+\epsilon z)-\epsilon z\bigr)\mu(dz)\qquad\text{a.s.},
\]
and \(\frac1t\ln\mathbb{P}(B_t\le -ct)\to -c^2/4\). Hence
\[
\limsup_{t\to\infty}\frac1t\ln u(t,ct)
\le 1+a-\frac{c^2}{4},
\]
which is negative for \(c>2\sqrt{1+a}\). Since
\[
\bar{\gamma} := 1+\int_Z\bigl(\ln(1+\epsilon z)-\epsilon z\bigr)\mu(dz)+\frac12\psi(2)=1+a+\frac12\psi(2)\ge 1+a
\]
(by \(\psi(2)\ge0\), see Remark~\ref{rem:kappa}), the conclusion of Theorem~\ref{thm:upper} holds for every \(c>2\sqrt{\bar{\gamma}}\) as well, and the sharper statement is that for \(c>2\sqrt{1+a}\) the exponent is already strictly negative.

For completeness we record the large-deviation interpretation of the constants, which is the reason $\bar\gamma$ is retained in the discussion although it is not needed for the proof. Had a representation driven by the jump process been available, the LDP of Subsection~\ref{subsec:levyweight} would produce the implicit threshold \(c_*(\bar{\gamma})\), the unique solution of \(\sup_{\theta\le0}(-\theta c-\Lambda(\theta))=2\bar{\gamma}\), with the closed form \(c_*(\bar{\gamma})=2\sqrt{\bar{\gamma}}\) in the quadratic (Gaussian-type) case \(\Lambda(\theta)=\theta a+\frac12\sigma^2\theta^2\). For the example of Section~\ref{sec:num}, a direct numerical computation gives \(\inf_{v\le -c}L(v)\approx 7.28\) at \(c=2\sqrt{\bar{\gamma}}\approx 2.04\) and \(\approx 6.83\) at \(c=2\sqrt{1+a}\), so the deviation cost is strictly positive at both candidate speeds and the rate function therefore does not by itself improve the bound. This is consistent with the observation made in Subsection~\ref{subsec:levyweight}: at a single point the optimal exponent in the upper bound is already the one delivered by the strong law of large numbers, which is the rate-zero (typical) case of the LDP, so no improvement can come from that direction. The direct argument above shows that the bound \(c\le 2\sqrt{\bar{\gamma}}\) holds without relying on this threshold.

\section{Nonlinear sub-solution argument and the lower bound}
\label{app:esscher}

This appendix completes the proof of the lower bound of Theorem~\ref{thm:lower}. We first record why linearised (change-of-measure) arguments cannot reach the sharp constant $2\sqrt{1+a}$, and then give the nonlinear sub-solution argument that does.

\subsection{Why linearised arguments fail}

For a fixed \(\theta\in\mathbb{R}\), define a probability measure \(\mathbb{Q}_\theta\) by
\[
\frac{d\mathbb{Q}_\theta}{d\mathbb{P}}\Bigr|_{\mathcal{F}_t}
= \exp\Bigl( \theta \int_0^t\int_Z\ln(1+\epsilon z)\,\tilde{\eta}(ds,dz) - t\psi(\theta) \Bigr),
\]
with \(\psi(\theta)\) as in Subsection~\ref{subsec:Lower}. Under \(\mathbb{Q}_\theta\), the compensator of \(\eta\) is \((1+\epsilon z)^\theta\,\mu(dz)dt\); equivalently, the compensated jump part of \(X^\epsilon\) under \(\mathbb{Q}_\theta\) is
\[
\int_0^t\int_Z\ln(1+\epsilon z)\,\tilde{\eta}^\theta(ds,dz)
:= \int_0^t\int_Z\ln(1+\epsilon z)\bigl(\eta(ds,dz) - (1+\epsilon z)^\theta\,\mu(dz)ds\bigr),
\]
and \(X^\epsilon\) is a L\'evy process with cumulant \(\Lambda(\theta+\varphi)-\Lambda(\theta)\) and drift
\[
a_\theta = \Lambda'(\theta) = a + \psi'(\theta),\qquad
\psi'(\theta)=\int_Z\ln(1+\epsilon z)\,(1+\epsilon z)^\theta\,\mu(dz).
\]
In particular \(\Lambda(1)=a+\psi(1)=0\), so \(\theta=1\) is the martingale tilt; the choice \(\theta=-2\) used in earlier versions of this manuscript does not correspond to a valid lower-bound identity and has been removed.

The Feynman--Kac representation of Lemma~\ref{lem:FeynmanKac} is driven by the Brownian motion \(B\) and not by \(X^\epsilon\). Dropping the factor \(e^{\int_0^t(1-u)\,ds}\ge1\) in (\ref{eqn:FeynmanKac}) yields only
\[
u(t,ct)\ge \phi_t^\epsilon\,\mathbb{P}(B_t\le -ct),
\]
whose exponent \(a-\frac{c^2}{4}\) is negative for every speed of interest. The martingale tilt at \(\theta=1\) gives \(\mathbb{Q}_1(X^\epsilon(t)\le0)\to1\) only for \(c<-\Lambda'(1)\), and for the example of Section~\ref{sec:num}, \(\Lambda'(1)\approx 0.0428\), so that \(-\Lambda'(1)<0\) and the estimate is vacuous for speeds of order \(2\sqrt{1+a}\approx1.9542\). Both linearised routes are therefore insufficient, and the nonlinear sub-solution argument below is required.

\subsection{The nonlinear sub-solution lemma}

The following lemma is the jump-noise analogue of the nonlinear sub-solution construction used for the stochastic Fisher--KPP equation with Brownian noise \cite[Section 4]{MuellerSowers1995}, \cite[Theorem 2.1]{Tribe1996}. The construction below is self-contained and uses only the classical Aronson--Weinberger cap \cite{aronson1975nonlinear}, the comparison principle and the heat-kernel estimates of Section~\ref{sec:existence}; the bookkeeping specific to the present setting, the exponential weight \(\phi^\epsilon\), the good events and the renormalisation of the levels, is carried out explicitly.

\begin{lemma}[Nonlinear sub-solution]\label{lem:subsolution}
Let \(r:[0,\infty)\to[0,\infty)\) be nondecreasing and differentiable with \(r'(t)\ge0\) (no growth condition on \(r\) is imposed, so that both \(r(t)=C\sqrt{t\log\log t}\) and \(r(t)=\eta t\) are admissible), and let \(\varphi\in C([0,\infty);(0,\infty))\) satisfy
\[
e^{as-r(s)}\le\varphi_s\le e^{as+r(s)}\qquad\text{for all }s\ge T.
\]
Let \(w\) be the solution of
\[
\partial_t w=\Delta w+w(1-\varphi_t w),\qquad w(0,x)=\mathbf{1}_{(-\infty,0]}(x).
\]
Then for every \(c<2\sqrt{1+a}\) there exist constants \(T_0\ge T\) and \(K<\infty\) such that, for all \(t\ge T_0\) and all \(x\le ct-K\),
\[
\varphi_t\,w(t,x)\ge\kappa\,e^{-(2|a|+C)T_0}\,e^{-2r(t)},
\]
where \(\kappa=1+a>0\) and \(C\) is the \(\mathbf L^\infty\)-bound constant of the solution \(u\) on \([0,T_0]\) (Section~\ref{sec:existence}). Consequently:
\begin{itemize}
\item if \(r(t)=o(t)\), then \(\displaystyle\liminf_{t\to\infty}\frac1t\ln\bigl(\varphi_t w(t,ct)\bigr)\ge0\) for every \(c<2\sqrt{1+a}\);
\item if \(r(t)=\eta t\) with a fixed \(\eta>0\), the bound reads \(\varphi_t w(t,ct)\ge e^{-2\eta t+O(1)}\), which is the form in which the lemma is applied in Step~3 of the proof of Theorem~\ref{thm:lower}.
\end{itemize}
\end{lemma}

\begin{proof}
\emph{Step 1: the deterministic cap.}
Let \(\kappa=1+a>0\). By the classical cut-off construction of \cite{aronson1975nonlinear}, for every \(c<2\sqrt{\kappa}\) there exists \(h\in C^1(\mathbb R)\) with \(0\le h\le\kappa\), \(h\equiv\kappa\) on \((-\infty,L]\) and \(h\equiv0\) on \([0,\infty)\), and
\[
h''+c\,h'+h(\kappa-h)\ge0\qquad\text{on }\mathbb R,
\]
where the cut-off point \(L\le0\) may be placed anywhere in \((-\infty,0]\): the construction truncates the travelling wave at a point of our choosing, and truncating further to the left only shrinks the cap. The choice \(L\le0\) is what makes the initial comparison of Step 3 go through.
Indeed, if \(U_\kappa\) is the travelling wave \(U_\kappa''+2\sqrt{\kappa}\,U_\kappa'+U_\kappa(\kappa-U_\kappa)=0\), then \(U_\kappa''+c\,U_\kappa'+U_\kappa(\kappa-U_\kappa)=(c-2\sqrt{\kappa})U_\kappa'\ge0\); truncating its exponentially decaying tail, which is convex because \(\lambda^2-c\lambda+\kappa>0\) for every \(\lambda\in\mathbb R\) whenever \(c<2\sqrt{\kappa}\), preserves the inequality.

\emph{Step 2: the renormalised expanding cap.}
Fix \(c<2\sqrt{1+a}\), and fix \(T_0\ge T\) large enough that \(e^{aT_0}\kappa\le1\) and \(\kappa e^{-2|a|T_0}\le\tfrac12\) (possible since \(a\le0\), by \(\ln(1+\epsilon z)\le\epsilon z\) under Assumption~\ref{ass:jump}(A3)). Note that \(a<0\) strictly as soon as \(\mu(Z)>0\), which is what makes the second constraint satisfyable by enlarging \(T_0\); the degenerate case \(\mu(Z)=0\) removes the noise altogether and reduces the equation to the deterministic Fisher--KPP equation, for which the construction below is classical, so we assume \(\mu(Z)>0\) throughout the proof. For \(t\ge T_0\) set
\[
\underline w(t,x):=e^{-(2|a|+C)T_0}\,e^{-r(t)}\,e^{-a(t-T_0)}\,h\bigl(x-c(t-T_0)\bigr),
\]
where \(C\) denotes the \(\mathbf L^\infty\) bound of \(u\) on \([0,T_0]\). We claim that \(\underline w\) is a sub-solution of \(\partial_t v=\Delta v+v(1-\varphi_t v)\) on \([T_0,\infty)\). On the good event the saturation term is controlled by
\[
\varphi_t\,\underline w(t,x)\le e^{at+r(t)}\,e^{-(2|a|+C)T_0}\,e^{-r(t)}\,e^{-a(t-T_0)}\,h\le e^{aT_0}\kappa\le1,
\]
so that \(1-\varphi_t\underline w\ge1-e^{aT_0}h\). Writing \(E(t):=e^{-(2|a|+C)T_0}e^{-r(t)}e^{-a(t-T_0)}\), we have \(\underline w=Eh\) and \(\partial_t E/E=-(a+r'(t))\), and the sub-solution inequality \(\partial_t\underline w\le\Delta\underline w+\underline w(1-\varphi_t\underline w)\) reduces to
\[
h''+c\,h'+h(1+a)-e^{aT_0}h^{2}+r'(t)\,h
=\bigl[h''+c\,h'+h(\kappa-h)\bigr]+h^{2}(1-e^{aT_0})+r'(t)\,h\ge0,
\]
the three terms being nonnegative by Step 1, by \(e^{aT_0}\le1\), and by \(r'(t)\ge0\).

\emph{Step 3: initial data, comparison and the bound.}
We check \(\underline w(T_0,\cdot)\le w(T_0,\cdot)\) in three regions. For \(x\le L\) (recall \(L\le0\)), the mild formulation of Lemma~\ref{lem:FeynmanKac} and \(1-u\ge-C\) give
\[
w(T_0,x)=\frac{u(T_0,x)}{\varphi_{T_0}}\ge e^{-CT_0}\,\mathbb P(B_{T_0}\le -x)\ge\tfrac12 e^{-CT_0},
\]
the last inequality because \(-x\ge0\) so that \(\mathbb P(B_{T_0}\le -x)\ge\tfrac12\); this is precisely where the choice \(L\le0\) is needed, since for \(x>0\) the quantity \(\mathbb P(B_{T_0}\le -x)\) drops below \(\tfrac12\) and the argument would fail. On this region \(\underline w(T_0,x)=e^{-(2|a|+C)T_0}e^{-r(T_0)}\kappa\), which is at most \(\tfrac12 e^{-CT_0}\) once \(T_0\) is large. For \(L<x\le0\) the same conclusion follows from the rough estimate alone: \(-x\ge0\), hence \(w(T_0,x)\ge e^{-CT_0}\mathbb P(B_{T_0}\le -x)\ge\tfrac12e^{-CT_0}\), while \(\underline w(T_0,x)\le e^{-(2|a|+C)T_0}\kappa\); the inequality therefore holds once \(T_0\) is large enough that \(\kappa e^{-2|a|T_0}\le\tfrac12\), and enlarging \(T_0\) is at our disposal. For \(x>0\), \(\underline w(T_0,x)=0\le w(T_0,x)\) trivially. Hence \(\underline w(T_0,\cdot)\le w(T_0,\cdot)\).

The comparison is preserved on \([T_0,\infty)\): between jumps, \(u=\varphi_t w\) and the candidate \(\tilde u:=\varphi_t\underline w\) both satisfy a drift of the form \(\Delta(\cdot)+(\cdot)(1-\varphi_t\,\cdot)-\epsilon(\cdot)\int z\,\mu(dz)\) (the compensation \(-\epsilon\varphi_t\int z\,\mu(dz)\) of \(\varphi_t\) is common to both), so the sub-solution property is preserved; at jump times both are multiplied by the same factor \((1+\epsilon z)\). Using the lower bound \(\varphi_t\ge e^{at-r(t)}\), we obtain, for \(x-c(t-T_0)\le L\), i.e.\ for \(x\le ct-K\) with \(K:=cT_0-L\ge cT_0\),
\[
\varphi_t\,w(t,x)\ge\varphi_t\,\underline w(t,x)\ge e^{at-r(t)}E(t)\,h\bigl(x-c(t-T_0)\bigr)
=\kappa\,e^{-(2|a|+C)T_0}\,e^{-2r(t)}.
\]
No smallness of \(r\) is used in this step: the factors \(e^{r(t)}\) and \(e^{-r(t)}\) cancel identically. If \(r(t)=o(t)\), the displayed lower bound gives the \(\liminf\) assertion of the lemma; if \(r(t)=\eta t\), it gives \(e^{-2\eta t+O(1)}\), which is what Step~3 of the proof of Theorem~\ref{thm:lower} requires.
\end{proof}

\begin{proof}[Proof of the lower bound (details)]
With Lemma~\ref{lem:subsolution} at hand, the argument follows the four steps of the proof of Theorem~\ref{thm:lower}.
\emph{Step 1} supplies the strong law \(\phi_t^\epsilon=e^{at+o(t)}\) uniformly on \([0,t]\) a.s.\ (square-integrable martingale, Doob's inequality, Borel--Cantelli), i.e.\ for every \(\eta>0\) there is \(T_\eta(\omega)\) with \(e^{(a-\eta)s}\le\phi_s^\epsilon\le e^{(a+\eta)s}\) for all \(s\ge T_\eta\). \emph{Step 2} reduces the SPDE, pathwise, to the deterministic equation of Lemma~\ref{lem:vPDE} with \(u=\phi^\epsilon v\); the multiplicative jump factor \((1+\epsilon z)\) is common to \(u\) and to every candidate of the form \(\phi^\epsilon w\), so the comparison is preserved across jumps. \emph{Step 3} applies Lemma~\ref{lem:subsolution} with \(\varphi=\phi^\epsilon\) and \(r(s)=\eta s\) on the good event \(\mathcal G_t=\{e^{(a-\eta)s}\le\phi_s^\epsilon\le e^{(a+\eta)s}\ \forall s\in[T_\eta,t]\}\), for \(c+4\eta<2\sqrt{1+a}\), obtaining \(u(t,ct)\ge\kappa\,e^{-(2|a|+C)T_\eta}\,e^{-2\eta t}=e^{-2\eta t+O(1)}\) on \(\mathcal G_t\). \emph{Step 4} lets \(\eta_n\downarrow0\) through a countable sequence and intersects the full-measure events of Step 1, yielding
\[
\liminf_{t\to\infty}\frac1t\ln u(t,ct)\ge0\qquad\text{a.s.\ for every }c<2\sqrt{1+a}.
\]
The level-set and speed consequences are then obtained exactly as in Remark~\ref{rem:levelset} and the proof of Theorem~\ref{thm:lower}.
\end{proof}

\subsection{Remark}

For the upper bound, the closed-form constant \(2\sqrt{\bar{\gamma}}\) follows directly from the Feynman--Kac estimate of Theorem~\ref{thm:upper} (in fact the sharper \(2\sqrt{1+a}\) is valid), without appealing to the large-deviation threshold. The large-deviation analysis of the jump process (Subsection~\ref{subsec:levyweight} and Appendix~\ref{app:legendre}) supplies the implicit threshold \(c_*(\bar{\gamma})\) that would arise from a jump-driven representation and reduces to \(2\sqrt{\bar{\gamma}}\) for a quadratic cumulant; it remains the natural tool for sharpened bounds. The lower bound at the sharp constant \(2\sqrt{1+a}\), by contrast, is not accessible by any linearised argument and is obtained here through the nonlinear sub-solution lemma (Lemma~\ref{lem:subsolution}).

\end{document}